\documentclass[a4paper, 12pt]{article}
\pdfoutput=1

\usepackage[margin=25mm]{geometry}

\usepackage[utf8]{inputenc}
\usepackage[T1]{fontenc}
\usepackage[british]{babel}
\usepackage[style=british]{csquotes}
\usepackage[en-GB]{datetime2}
\DTMlangsetup[en-GB]{ord=omit}
\usepackage{xcolor}

\usepackage{amsthm, thmtools}
\usepackage{mathtools}
\usepackage{titlesec} 
\usepackage{hyperref}
\hypersetup{
    bookmarksnumbered,
    colorlinks,
    linkcolor={cyan!50!blue},
    citecolor={cyan!50!blue},
    urlcolor={cyan!50!blue}
}
\usepackage[capitalise,nameinlink]{cleveref}

\usepackage{tikz}
\usepackage{tikz-cd}

\usepackage[nopatch=footnote]{microtype}
\usepackage[lining,tabular,scale=.9]{FiraSans}
\usepackage[tt=false,semibold]{libertine}
\usepackage[libertine,timesmathacc,smallerops]{newtxmath}
\usepackage[lining,scaled=.8]{FiraMono}
\usepackage[cal=boondox,frak=euler]{mathalpha}

\DeclareFontEncoding{MDA}{}{}
\DeclareSymbolFont{mathdesignA}{MDA}{mdput}{m}{n}
\SetSymbolFont{mathdesignA}{bold}{MDA}{mdput}{b}{n}
\DeclareSymbolFontAlphabet{\mathbb}{mathdesignA}

\tikzcdset{
    arrow style=tikz,
    arrows={/tikz/line width=.5pt},
    diagrams={>={Straight Barb[scale=0.8]}}
}

\DeclareFontFamily{OMX}{MnSymbolE}{}
\DeclareSymbolFont{MnLargeSymbols}{OMX}{MnSymbolE}{m}{n}
\SetSymbolFont{MnLargeSymbols}{bold}{OMX}{MnSymbolE}{b}{n}
\DeclareFontShape{OMX}{MnSymbolE}{m}{n}{
    <-6>  MnSymbolE5
   <6-7>  MnSymbolE6
   <7-8>  MnSymbolE7
   <8-9>  MnSymbolE8
   <9-10> MnSymbolE9
  <10-12> MnSymbolE10
  <12->   MnSymbolE12
}{}
\DeclareFontShape{OMX}{MnSymbolE}{b}{n}{
    <-6>  MnSymbolE-Bold5
   <6-7>  MnSymbolE-Bold6
   <7-8>  MnSymbolE-Bold7
   <8-9>  MnSymbolE-Bold8
   <9-10> MnSymbolE-Bold9
  <10-12> MnSymbolE-Bold10
  <12->   MnSymbolE-Bold12
}{}
\DeclareMathDelimiter{[}{\mathopen}{MnLargeSymbols}{'000}{MnLargeSymbols}{'000}
\DeclareMathDelimiter{]}{\mathclose}{MnLargeSymbols}{'005}{MnLargeSymbols}{'005}
\DeclareMathDelimiter{\llbr}{\mathopen}{MnLargeSymbols}{'102}{MnLargeSymbols}{'102}
\DeclareMathDelimiter{\rrbr}{\mathclose}{MnLargeSymbols}{'107}{MnLargeSymbols}{'107}

\usepackage{bm}
\usepackage{cases}
\usepackage{accents}

\usepackage{setspace}
\usepackage{subdepth}

\newcommand{\initlengths}{%
    \setlength{\abovedisplayshortskip}{3pt plus 9pt minus 3pt}%
    \setlength{\belowdisplayshortskip}{9pt plus 9pt minus 9pt}%
    \setlength{\abovedisplayskip}{9pt plus 9pt minus 9pt}%
    \setlength{\belowdisplayskip}{9pt plus 9pt minus 9pt}%
    \tolerance 500
}

\usepackage{titling}
\pretitle{\vspace{-18pt}\setstretch{1.05}\begin{center}\LARGE\libertineDisplay\fontdimen2\font=0.3em} 
\posttitle{\end{center}\vspace{6pt}}

\numberwithin{paragraph}{subsection}
\newcommand{\parafont}{\bfseries}
\newcommand{\parasep}{9pt plus 3pt minus 3pt}

\titleformat{\section}{\Large\libertineDisplay}{\thesection}{1em}{}
\titleformat{\subsection}{\large\firamedium\boldmath}{\thesubsection}{1em}{}
\titleformat{\paragraph}[runin]{\parafont}{\theparagraph.}{.33em}{\normalfont\bfseries\boldmath}
\titlespacing{\paragraph}{0pt}{\parasep}{.5em}

\usepackage{needspace}

\usepackage{tocloft}

\renewenvironment{abstract}{%
    \centering\begin{minipage}{.85\textwidth}%
    \setlength{\parindent}{1.5em}%
    \centerline{\large\firamedium\abstractname}%
    \par\vspace{12pt}%
}{\end{minipage}\par\vspace{3pt}}

\newcommand{\authorinforule}{\noindent\rule{0.38\textwidth}{0.4pt}}

\newlength{\authorwidth}
\newcommand{\authorinfo}[3]{%
    \setlength{\leftskip}{1.5em}
    \setlength{\parindent}{0em}
    \setstretch{1}
    \par%
    {\small%
    \makebox[\authorwidth][l]{#1}%
    \texttt{#2}%
    \\
    #3.}
    \vspace{6pt}\par
}

\makeatletter
\renewcommand{\Hy@numberline}[1]{#1. }
\newcommand*{\@parabookmark}{%
  \pdfbookmark[3]{%
    \theparagraph
    \ifx\@currentlabelname\@empty
    \else
      .\space\@currentlabelname%
    \fi
  }{\theparagraph}
}
\newcommand*{\@defbookmark}{%
  \pdfbookmark[3]{%
    \theparagraph.\space\thmt@thmname
    \ifx\@currentlabelname\@empty
    \else
      .\space\@currentlabelname%
    \fi
  }{\theparagraph}
}
\newcommand*{\@thmbookmark}{%
  \pdfbookmark[3]{%
    \theparagraph.\space\thmt@thmname
    \ifx\@currentlabelname\@empty
    \else
      \space(\@currentlabelname)%
    \fi
  }{\theparagraph}
}
\newcommand*{\parabookmark}{\@parabookmark}
\newcommand*{\defbookmark}{\@defbookmark}
\newcommand*{\thmbookmark}{\@thmbookmark}
\newcommand*{\suppressparabookmarks}{%
  \renewcommand*{\parabookmark}{}%
  \renewcommand*{\defbookmark}{}%
  \renewcommand*{\thmbookmark}{}%
}
\newcommand*{\resumeparabookmarks}{%
  \renewcommand*{\parabookmark}{\@parabookmark}%
  \renewcommand*{\defbookmark}{\@defbookmark}%
  \renewcommand*{\thmbookmark}{\@thmbookmark}%
}

\declaretheoremstyle[
    spaceabove=\parasep, spacebelow=\parasep,
    postheadspace=.5em,
    postheadhook=\thmbookmark,
    headfont=\normalfont\bfseries,
    headpunct={},
    headformat={\NUMBER.\@\ \NAME\NOTE.},
    notefont=\normalfont,
    notebraces={(}{)},
    bodyfont=\itshape,
]{theorem}
\declaretheoremstyle[
    spaceabove=\parasep, spacebelow=\parasep,
    postheadspace=.5em,
    headfont=\normalfont\bfseries,
    headpunct={},
    headformat={\NAME.\@\NOTE},
    notefont=\normalfont\bfseries\boldmath,
    notebraces={}{.},
    bodyfont=\itshape,
]{theorem*}
\declaretheoremstyle[
    spaceabove=\parasep, spacebelow=\parasep,
    postheadspace=.5em,
    postheadhook=\defbookmark,
    headfont=\normalfont\bfseries,
    headpunct={},
    headformat={\NUMBER.\@\ \NAME.\@\NOTE},
    notefont=\normalfont\bfseries\boldmath,
    notebraces={}{.},
]{definition}
\declaretheoremstyle[
    spaceabove=\parasep, spacebelow=\parasep,
    postheadspace=.5em,
    postheadhook=\parabookmark,
    headfont=\normalfont\bfseries,
    headpunct={},
    headformat={\NUMBER.\@\NOTE},
    notefont=\normalfont\bfseries\boldmath,
    notebraces={}{.},
]{para}

\renewenvironment{proof}[1][\proofname]{\par
    \pushQED{\qed}%
    \normalfont\trivlist
    \item[\hskip\labelsep\bfseries #1\@addpunct{.}]\ignorespaces
}{%
    \popQED\endtrivlist\@endpefalse
}
\makeatother

\declaretheorem[sibling=paragraph, style=para, refname={\S,\S\S}]{para}
\declaretheorem[sibling=paragraph, style=theorem, name=Theorem]{theorem}
\declaretheorem[sibling=paragraph, style=theorem, name=Lemma]{lemma}
\declaretheorem[sibling=paragraph, style=theorem, name=Corollary]{corollary}
\declaretheorem[sibling=paragraph, style=theorem, name=Proposition]{proposition}

\declaretheorem[numbered=no, style=theorem*, name=Theorem]{theorem*}
\declaretheorem[numbered=no, style=theorem*, name=Lemma]{lemma*}
\declaretheorem[sibling=paragraph, style=definition, name=Definition]{definition}
\declaretheorem[sibling=paragraph, style=definition, name=Example]{example}

\numberwithin{equation}{paragraph}

\crefdefaultlabelformat{#2{\upshape#1}#3}
\crefname{figure}{Figure}{Figures}

\crefformat{equation}{#2{\upshape(#1)}#3}
\crefrangeformat{equation}{{\upshape#3(#1)#4--#5(#2)#6}}
\crefmultiformat{equation}{#2{\upshape(#1)}#3}{ and #2{\upshape(#1)}#3}{, #2{\upshape(#1)}#3}{, and~#2{\upshape(#1)}#3}

\crefformat{enumi}{#2{\upshape#1}#3}
\crefrangeformat{enumi}{{\upshape#3#1#4--#5#2#6}}
\crefmultiformat{enumi}{#2{\upshape#1}#3}{ and #2{\upshape#1}#3}{, #2{\upshape#1}#3}{, and~#2{\upshape#1}#3}

\crefformat{section}{#2{\upshape\S#1}#3}
\crefrangeformat{section}{{\upshape#3\S\S#1#4--#5#2#6}}
\crefmultiformat{section}{{\upshape#2\S\S#1#3}}{ and {\upshape#2#1#3}}{, {\upshape#2#1#3}}{, and~{\upshape#2#1#3}}
\crefformat{subsection}{#2{\upshape\S#1}#3}
\crefrangeformat{subsection}{{\upshape#3\S\S#1#4--#5#2#6}}
\crefmultiformat{subsection}{{\upshape#2\S\S#1#3}}{ and {\upshape#2#1#3}}{, {\upshape#2#1#3}}{, and~{\upshape#2#1#3}}
\crefformat{para}{#2{\upshape\S#1}#3}
\crefrangeformat{para}{{\upshape#3\S\S#1#4--#5#2#6}}
\crefmultiformat{para}{{\upshape#2\S\S#1#3}}{ and {\upshape#2#1#3}}{, {\upshape#2#1#3}}{, and~{\upshape#2#1#3}}

\usepackage{enumitem}
\setlist{noitemsep, listparindent=\parindent}
\setlist[enumerate]{label=\textnormal{(\roman*)}}

\usepackage{caption}
\DeclareCaptionFont{fira}{\firamedium\firaproportional}
\newcommand{\preparebibliography}{
    \phantomsection
    \addcontentsline{toc}{section}{References}
    \sloppy
    \setstretch{1.1}
    \renewcommand*{\bibfont}{\normalfont\small}
    \setlength{\bibitemsep}{0.25\baselineskip}
}

\usepackage[
    backend=biber,
    style=oxnum,
    giveninits,
    maxbibnames=5,
    maxcitenames=5,
    sorting=nyt
]{biblatex}

\DefineBibliographyStrings{english}{
  withafterword = {with an appendix by},
}

\DeclareFieldFormat[article]{version}{\IfDecimal{#1}{version~}{}#1}
\renewbibmacro*{series+number}{%
    \printfield{series}\isdot%
    \setunit*{\addcomma\addspace}%
    \printfield{number}}

\DeclareFieldFormat{pages}{%
  \iffieldundef{bookpagination}%
    {#1}%
    {\mkpageprefix[bookpagination]{#1}}%
}

\newcommand{\colim}{\operatornamewithlimits{colim}}
\newcommand{\crk}{\operatorname{crk}}
\newcommand{\dFilt}{\prescript{\smash{\mathrm{d}}}{}{\mathrm{Filt}}}
\newcommand{\dGrad}{\prescript{\smash{\mathrm{d}}}{}{\mathrm{Grad}}}
\newcommand{\dMap}{\prescript{\smash{\mathrm{d}}}{}{\mathrm{Map}}}
\newcommand{\dmathcal}[1]{\prescript{\smash{\mathrm{d}}\mspace{-2mu}}{}{}\mathcal{#1}}
\newcommand{\Filt}{\mathrm{Filt}}

\newcommand{\Gm}{\mathbb{G}_\mathrm{m}}
\newcommand{\Grad}{\mathrm{Grad}}
\newcommand{\Hom}{\mathrm{Hom}}
\newcommand{\id}{\mathrm{id}}

\newcommand{\longsimto}{\mathrel{\overset{\smash{\raisebox{-.8ex}{$\sim$}}\mspace{3mu}}{\longrightarrow}}}

\newcommand{\Map}{\mathrm{Map}}

\newcommand{\Q}{\mathbb{Q}}
\newcommand{\rk}{\operatorname{rk}}
\newcommand{\simto}{\mathrel{\overset{\smash{\raisebox{-.8ex}{$\sim$}}\mspace{3mu}}{\to}}}

\newcommand{\Spec}{\operatorname{Spec}}
\newcommand{\Z}{\mathbb{Z}}

\renewcommand{\geq}{\geqslant}
\renewcommand{\leq}{\leqslant}

\newcommand{\bA}{{\mathbb A}}

\newcommand{\bG}{{\mathbb G}}

\newcommand{\bN}{{\mathbb N}}

\newcommand{\bP}{{\mathbb P}}
\newcommand{\bQ}{{\mathbb Q}}
\newcommand{\bR}{{\mathbb R}}

\newcommand{\bZ}{{\mathbb Z}}

\newcommand{\cD}{{\mathcal D}}

\newcommand{\cM}{{\mathcal M}}

\newcommand{\cU}{{\mathcal U}}
\newcommand{\cV}{{\mathcal V}}

\newcommand{\cX}{{\mathcal X}}
\newcommand{\cY}{{\mathcal Y}}
\newcommand{\cZ}{{\mathcal Z}}

\newcommand{\restr}{\mathord{\downarrow}} 

\title{Intrinsic Donaldson--Thomas theory\\I. Component lattices of stacks}
\author{\clap{Chenjing Bu\quad Daniel Halpern-Leistner\quad Andrés Ibáñez Núñez\quad Tasuki Kinjo}}
\date{}

\begin{document}

\initlengths

\maketitle

\begin{abstract}
    This is the first paper in a series on
\emph{intrinsic Donaldson--Thomas theory},
where we develop a new framework for enumerative geometry
that allows the generalization of constructions and results
from linear moduli stacks to general non-linear algebraic stacks.

In this paper, we introduce the \emph{component lattice} of an algebraic stack.
This is a key object in our theory,
defined using the formalism of stacks of graded and filtered points.
It provides the combinatorial data
needed to formulate various results in enumerative geometry,
such as decomposition-type theorems and wall-crossing formulae.
Later papers in the series will focus on extending Donaldson--Thomas theory
to the non-linear case,
and we expect that our approach will be useful for
extending many other flavours of enumerative invariants
beyond the linear case as well.

This paper proves several foundational results of our framework.
The first is the \emph{constancy theorem},
which states that the isomorphism types of connected components
of the stacks of graded and filtered points
stay constant within chambers in the component lattice.
The second is the \emph{finiteness theorem},
which provides a criterion for the finiteness of
the number of possible isomorphism types of these components.
The third is the \emph{associativity theorem},
generalizing the structure of Hall algebras
from linear stacks to general stacks.

We also discuss some applications of these results
outside Donaldson--Thomas theory,
including a construction of stacks of real-weighted filtrations,
and a generalization of the semistable reduction theorem
to real-weighted filtrations.

\end{abstract}

\clearpage

{
    \hypersetup{linkcolor=black}
    \setstretch{1.2}
    \tableofcontents
}

\clearpage
\suppressparabookmarks
\section{Introduction}

\addtocounter{subsection}{1}

\begin{para}[Overview]
    \label{para-intro-overview}
    This is the first paper in a series
    \cite{part-ii,part-iii}
    on \emph{intrinsic Donaldson--Thomas theory},
    a new framework for studying
    the enumerative geometry of general algebraic stacks.

    The theory of \emph{Donaldson--Thomas invariants}
    has been a central topic in enumerative geometry.
    It was initiated in the works of \textcite{donaldson-thomas-1998,thomas-2000-dt},
    and further developed by
    \textcite{joyce-song-2012-dt},
    \textcite{kontsevich-soibelman-motivic-dt},
    and many others.
    These invariants are defined for moduli stacks of objects
    in $3$-Calabi--Yau abelian categories,
    and count semistable objects in the category in a certain sense.

    One major difficulty in constructing enumerative invariants
    is in dealing with \emph{stacky points} in the moduli stack,
    or points corresponding to objects with
    an automorphism group of positive dimension.
    Many existing enumerative theories
    are limited to the case of abelian categories
    primarily because of this difficulty.
    For example, this is the case for
    the aforementioned Donaldson--Thomas theory,
    the theory of homological enumerative invariants of
    \textcite{joyce-wall-crossing},
    the theory of \emph{Vafa--Witten invariants}
    of \textcite{tanaka-thomas-2020-i,tanaka-thomas-2018-ii},
    and the theory of \emph{DT4 invariants}
    originating from the works of
    \textcite{cao-leung-dt4},
    \textcite{borisov-joyce-2017},
    \textcite{oh-thomas-2023-i,oh-thomas-ii},
    and others.

    More recently, \textcite{bu-self-dual-1,bu-self-dual-2}
    developed an approach to addressing this issue
    for orthogonal and symplectic enumerative problems.
    However, it has not been clear how to deal with this difficulty
    beyond these special cases.

    The goal of this series of papers is to develop a new framework
    for enumerative theories of general algebraic stacks,
    which is not limited to the linear case of abelian categories.
    While we expect that our framework will be useful for
    extending many other types of enumerative theories beyond the linear case,
    this series will focus on the construction of
    \emph{Donaldson--Thomas invariants}
    for general $(-1)$-shifted symplectic stacks,
    which are intrinsic to the stack,
    and depend on a certain notion of stability data attached to it.

    In this first paper, we lay out the foundations of our framework,
    involving a detailed study of the stacks of \emph{graded and filtered points}
    of an algebraic stack,
    introduced by \textcite{halpern-leistner-instability},
    which provide a way to generalize
    filtrations in an abelian category to an arbitrary stack.
    We show that the set of connected components of these stacks,
    called the \emph{component lattice},
    has interesting combinatorial structures,
    and study its interactions with the geometry of the stack.
    We establish some fundamental properties of the component lattice,
    which will be used in the following papers to define and study
    the intrinsic Donaldson--Thomas invariants,
    and we expect them to be useful for other types of enumerative theories as well.

    In a closely related work by
    \textcite{bu-davison-ibanez-nunez-kinjo-padurariu},
    the component lattice plays a crucial role in
    extending aspects of cohomological Donaldson--Thomas theory
    beyond the linear case,
    leading to an explicit formula for
    the decomposition theorem of the good moduli space morphism
    for a large class of stacks,
    generalizing the \emph{cohomological integrality theorem} of
    \textcite{davison-meinhardt-2020}.

    In the remainder of this introduction, we give a brief summary
    of the main ideas involved in this paper,
    as well as the main results and applications.
\end{para}

\begin{para}[Graded and filtered points]
    \label{para-intro-grad-filt}
    We briefly introduce the \emph{graded and filtered points}
    of an algebraic stack,
    following \textcite{halpern-leistner-instability},
    which generalize the notion of filtrations in an abelian category
    to an arbitrary algebraic stack.

    Given an algebraic stack~$\mathcal{X}$,
    its \emph{stack of graded points} and \emph{stack of filtered points}
    are defined as the mapping stacks
    \begin{align*}
        \mathrm{Grad} (\mathcal{X})
        & =
        \mathrm{Map} (* / \mathbb{G}_\mathrm{m}, \mathcal{X}) \ ,
        \\
        \mathrm{Filt} (\mathcal{X})
        & =
        \mathrm{Map} (\mathbb{A}^1 / \mathbb{G}_\mathrm{m}, \mathcal{X}) \ ,
    \end{align*}
    which are again algebraic stacks under mild conditions on~$\mathcal{X}$.

    If~$\mathcal{X}$ is a moduli stack of objects
    in an abelian category~$\mathcal{A}$,
    then $\mathrm{Grad} (\mathcal{X})$ is the stack of
    $\mathbb{Z}$-graded objects in~$\mathcal{A}$,
    and $\mathrm{Filt} (\mathcal{X})$ is the stack of
    $\mathbb{Z}$-filtrations in~$\mathcal{A}$.

    On the other hand, the stacks $\mathrm{Grad} (\mathcal{X})$ and $\mathrm{Filt} (\mathcal{X})$
    have explicit descriptions
    when~$\mathcal{X}$ is a quotient stack.
    For simplicity, consider a quotient stack~$V / G$,
    where~$G$ is a reductive group over an algebraically closed field~$K$,
    and~$V$ is a finite-dimensional $G$-representation.
    Then, as in \textcite[Theorem~1.4.8]{halpern-leistner-instability},
    we have
    \begin{align*}
        \mathrm{Grad} (V / G)
        & =
        \coprod_{\lambda \colon \mathbb{G}_\mathrm{m} \to G}
        V^{\lambda, 0} / L_\lambda \ ,
        \\
        \mathrm{Filt} (V / G)
        & =
        \coprod_{\lambda \colon \mathbb{G}_\mathrm{m} \to G}
        V^{\lambda, +} / P_\lambda \ ,
    \end{align*}
    where~$\lambda$ runs over conjugacy classes of \emph{cocharacters} of~$G$,
    that is, maps~$\mathbb{G}_\mathrm{m} \to G$.
    Here, $V^{\lambda, 0}, V^{\lambda, +} \subset V$
    are the parts where~$\mathbb{G}_\mathrm{m}$
    acts with zero and non-negative weights through~$\lambda$, respectively,
    and $L_\lambda, P_\lambda \subset G$
    are the Levi and parabolic subgroups corresponding to~$\lambda$.

    In particular, if~$T \subset G$ is a maximal torus,
    then the connected components of the stacks
    $\mathrm{Grad} (V / G)$ and $\mathrm{Filt} (V / G)$
    are indexed by the set~$\Lambda_T / W$,
    where~$\Lambda_T = \mathrm{Hom} (\mathbb{G}_\mathrm{m}, T)$
    is the cocharacter lattice of~$T$,
    and~$W$ is the Weyl group.
\end{para}

\begin{para}[The component lattice]
    \label{para-intro-component-lattice}
    We define the \emph{component lattice}
    of a stack~$\mathcal{X}$ to be the set
    $\mathrm{CL} (\mathcal{X}) = \uppi_0 (\mathrm{Grad} (\mathcal{X}))
    \simeq \uppi_0 (\mathrm{Filt} (\mathcal{X}))$,
    equipped with a certain extra structure.
    For example, for a quotient stack~$V / G$ as above,
    we have $\mathrm{CL} (V / G) \simeq \Lambda_T / W$.

    We explain some key ideas that will appear in this paper
    in a concrete example.
    Consider the group~$G = \mathrm{GL} (2)$
    and its standard $2$-dimensional representation $V = \mathbb{A}^2$,
    and consider the quotient stack~$V / G$.
    The stacks $\mathrm{Grad} (V / G)$ and $\mathrm{Filt} (V / G)$
    are depicted in \cref{fig-intro-component-lattice}.
    The component lattice is $\mathrm{CL} (V / G) \simeq \mathbb{Z}^2 / W$,
    where $W \simeq \mathbb{Z}/2\mathbb{Z}$ acts by reflection across the diagonal line.
    In the picture, each dot corresponds to a connected component of
    $\mathrm{Grad} (V / G)$ or $\mathrm{Filt} (V / G)$,
    with their isomorphism types as indicated,
    and each dot is identified with its reflection across the diagonal line.
    The notations are as follows:
    $T \subset G$ is the maximal torus consisting of diagonal matrices,
    $B \subset G$ is the Borel subgroup of upper triangular matrices,
    $\bar{B} \subset G$ is the opposite Borel subgroup of lower triangular matrices,
    and~$V_1, V_2 \subset V$ are the one-dimensional subspaces containing
    the standard basis vectors.
    Note that we have isomorphisms $V_1 / B \simeq V_2 / \bar{B}$, etc.,
    allowing the identification of dots across the diagonal line.

    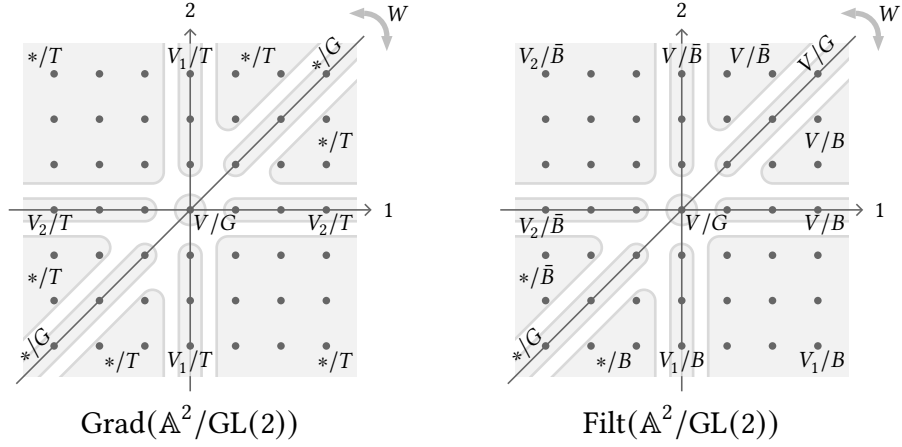
\begin{figure}[t]
        \centering
        \begin{tikzpicture}[>={Straight Barb[scale=.8]}]
            \def\drawbase{
                \def\fillcolor{black!5}

                \begin{scope}
                    \clip (-2.2, -2.2) rectangle (2.2, 2.2);
                    \draw[line width=1, black!15, fill=\fillcolor]
                        (.35, -2.4) -- (.35, -.5) arc (180:90:.15) -- (2.4, -.35) -- (2.4, -2.4) -- cycle
                        (-.35, 2.4) -- (-.35, .5) arc (0:-90:.15) -- (-2.4, .35) -- (-2.4, 2.4) -- cycle
                        (2.4, .35) -- (1.207107, .35) arc (-90:-225:.15) -- (2.4, 1.905025) -- cycle
                        (-2.4, -.35) -- (-1.207107, -.35) arc (90:-45:.15) -- (-2.4, -1.905025) -- cycle
                        (-.35, -2.4) -- (-.35, -1.207107) arc (0:135:.15) -- (-1.905025, -2.4) -- cycle
                        (.35, 2.4) -- (.35, 1.207107) arc (-180:-45:.15) -- (1.905025, 2.4) -- cycle
                        (2.4, .15) -- (.6, .15) arc (90:270:.15) -- (2.4, -.15) -- cycle
                        (-2.4, -.15) -- (-.6, -.15) arc (-90:90:.15) -- (-2.4, .15) -- cycle
                        (-.15, -2.4) -- (-.15, -.6) arc (180:0:.15) -- (.15, -2.4) -- cycle
                        (.15, 2.4) -- (.15, .6) arc (0:-180:.15) -- (-.15, 2.4) -- cycle
                        (2.187868, 2.4) -- (.493934, .706066) arc (135:315:.15) -- (2.4, 2.187868) -- cycle
                        (-2.187868, -2.4) -- (-.493934, -.706066) arc (-45:135:.15) -- (-2.4, -2.187868) -- cycle
                        (0, 0) circle (.2);
                \end{scope}
                \foreach \x in {-3, ..., 3} {
                    \foreach \y in {-3, ..., 3} {
                        \fill[black!60] (0.6 * \x, 0.6 * \y) circle (0.05);
                    }
                }
                \draw[line width=.6, black!60] (-2.3, -2.3) -- (2.3, 2.3);
                \draw[line width=.6, black!60, ->] (-2.4, 0) -- (2.4, 0)
                    node[right, black] {$\scriptstyle 1$};
                \draw[line width=.6, black!60, ->] (0, -2.4) -- (0, 2.4)
                    node[above, black] {$\scriptstyle 2$};

                \draw[line width=2, black!25] (2.3, 2.6) arc (90:0:.3);
                \fill[black!25] (2.3, 2.525) -- (2.3, 2.675) -- (2.1, 2.6)
                    (2.525, 2.3) -- (2.675, 2.3) -- (2.6, 2.1);
                \node[inner sep=0] at (2.75, 2.6) {$\scriptstyle W$};

                \fill[\fillcolor, opacity=.75]
                    (-.1, 2.2) rectangle (.1, 1.85)
                    (-.1, -2.2) rectangle (.1, -1.85);
            }

            \begin{scope}
                \drawbase
                \node[inner sep=0] at (0, -2.8) {$\mathrm{Grad} (\mathbb{A}^2 / \mathrm{GL} (2))$};

                \node[anchor=north west, inner sep=1] at (-2.2, 2.2) {$\scriptstyle * / T$};
                \node[anchor=north, inner sep=1] at (0, 2.2) {$\scriptstyle V_1 / T$};
                \node[anchor=north, inner sep=1] at (.9, 2.2) {$\scriptstyle * / T$};
                \node[rotate=45] at (1.8, 2.05) {$\scriptstyle * / G$};
                \node[anchor=east, inner sep=1] at (2.2, .9) {$\scriptstyle * / T$};
                \node[anchor=north east, inner sep=1] at (2.2, 0) {$\scriptstyle V_2 / T$};
                \node[anchor=south east, inner sep=1] at (2.2, -2.2) {$\scriptstyle * / T$};
                \node[anchor=south, inner sep=1] at (0, -2.2) {$\scriptstyle V_1 / T$};
                \node[anchor=south, inner sep=1] at (-.9, -2.2) {$\scriptstyle * / T$};
                \node[rotate=45] at (-2.05, -1.8) {$\scriptstyle * / G$};
                \node[anchor=west, inner sep=1] at (-2.2, -.9) {$\scriptstyle * / T$};
                \node[anchor=north west, inner sep=1] at (-2.2, 0) {$\scriptstyle V_2 / T$};
                \node[anchor=north west, inner sep=1] at (0, 0) {$\scriptstyle V / G$};
            \end{scope}

            \begin{scope}[shift={(6.5, 0)}]
                \drawbase
                \node[inner sep=0] at (0, -2.8) {$\mathrm{Filt} (\mathbb{A}^2 / \mathrm{GL} (2))$};

                \node[anchor=north west, inner sep=1] at (-2.2, 2.2) {$\scriptstyle V_2 / \smash{\bar{B}}$};
                \node[anchor=north, inner sep=1] at (0, 2.2) {$\scriptstyle V / \smash{\bar{B}}$};
                \node[anchor=north, inner sep=1] at (.9, 2.2) {$\scriptstyle V / \smash{\bar{B}}$};
                \node[rotate=45] at (1.8, 2.05) {$\scriptstyle V / G$};
                \node[anchor=east, inner sep=1] at (2.2, .9) {$\scriptstyle V / B$};
                \node[anchor=north east, inner sep=1] at (2.2, 0) {$\scriptstyle V / B$};
                \node[anchor=south east, inner sep=1] at (2.2, -2.2) {$\scriptstyle V_1 / B$};
                \node[anchor=south, inner sep=1] at (0, -2.2) {$\scriptstyle V_1 / B$};
                \node[anchor=south, inner sep=1] at (-.9, -2.2) {$\scriptstyle * / B$};
                \node[rotate=45] at (-2.05, -1.8) {$\scriptstyle * / G$};
                \node[anchor=west, inner sep=1] at (-2.2, -.9) {$\scriptstyle * / \smash{\bar{B}}$};
                \node[anchor=north west, inner sep=1] at (-2.2, 0) {$\scriptstyle V_2 / \bar{B}$};
                \node[anchor=north west, inner sep=1] at (0, 0) {$\scriptstyle V / G$};
            \end{scope}
        \end{tikzpicture}

        \caption{Graded and filtered points of~$\mathbb{A}^2 / \mathrm{GL} (2)$}
        \label{fig-intro-component-lattice}
    \end{figure}

    An important observation here
    is that the component lattice
    admits a wall-and-chamber structure,
    such that the isomorphism types of the components of
    $\mathrm{Grad} (\mathcal{X})$ and $\mathrm{Filt} (\mathcal{X})$
    stay constant within each region.
    In the above example, the walls are the three lines in the picture,
    dividing the lattice into $13$~regions as indicated,
    or $8$~regions after taking the quotient by~$W$.
\end{para}

\begin{para}[The constancy theorem]
    Our first main result,
    the \emph{constancy theorem},
    \cref{thm-constancy},
    is a generalization of the above observation
    to arbitrary algebraic stacks,
    and states that such a wall-and-chamber structure
    on the component lattice exists for
    general algebraic stacks.

    More precisely,
    we show that the component lattice of a stack~$\mathcal{X}$
    is covered by a collection of sublattices
    called \emph{special faces},
    which we introduce in \cref{subsec-faces}.
    These are roughly the walls
    of the wall-and-chamber structure mentioned above.
    For example, in \cref{fig-intro-component-lattice},
    the special faces are the whole plane,
    the coordinate axes, the diagonal line, and the origin.

    On each special face,
    we define a hyperplane arrangement called the \emph{cotangent arrangement}
    in \cref{subsec-cotangent-weights}
    under mild conditions,
    and the constancy theorem states that the isomorphism types of the components of
    $\mathrm{Grad} (\mathcal{X})$ and $\mathrm{Filt} (\mathcal{X})$
    stay constant within each chamber of the cotangent arrangement,
    as long as we do not enter a lower-dimensional special face.

    When~$\mathcal{X}$ is the moduli stack of objects in an abelian category,
    this result amounts to the following fact:
    For a map of sets $\gamma \colon \mathbb{Z} \to \uppi_0 (\mathcal{X})$
    such that $\gamma (n) = 0$ for all but finitely many~$n$,
    which we think of as choosing classes for the graded pieces
    of a $\mathbb{Z}$-graded object or a $\mathbb{Z}$-filtration,
    the moduli stacks of $\mathbb{Z}$-graded objects and $\mathbb{Z}$-filtrations
    with the chosen classes for the graded pieces,
    denoted by~$\mathcal{X}_\gamma$ and~$\mathcal{X}_\gamma^+$,
    only depend on the non-zero classes that appear in~$\gamma$
    and their order,
    but not on the exact gradings of these classes.
\end{para}

\needspace{2\baselineskip}
\begin{para}[The finiteness theorem]
    Our second main result is the \emph{finiteness theorem},
    \cref{thm-finiteness},
    which states that a quasi-compact stack~$\mathcal{X}$
    with \emph{quasi-compact graded points}
    has finitely many special faces.
    Here, having quasi-compact graded points
    is a very mild condition.

    By the constancy theorem,
    this will imply that for such a stack~$\mathcal{X}$,
    all the connected components of
    $\Grad (\mathcal{X})$ and $\Filt (\mathcal{X})$
    can only take finitely many isomorphism types.
\end{para}

\begin{para}[Hall categories and associativity]
    \label{para-intro-associativity}
    Our third main result is the \emph{associativity theorem},
    \cref{thm-associativity},
    which states that the commutative diagram
    \begin{equation*}
        \begin{tikzcd}[column sep=1em, row sep={3.5em,between origins}]
            \Filt (\Filt (\mathcal{X}))
            \ar[r] \ar[d]
            & \Filt (\Grad (\mathcal{X}))
            \mathrlap{{} \simeq \Grad (\Filt (\mathcal{X}))}
            \ar[d]
            \\
            \Filt (\mathcal{X})
            \ar[r]
            & \Grad (\mathcal{X})
        \end{tikzcd}
        \hspace{5em}
    \end{equation*}
    is cartesian on some connected components of
    $\Filt (\Filt (\mathcal{X}))$,
    although this is not true for all components.
    Here, the horizontal arrows take the associated graded objects,
    and the vertical arrows take the total objects of a filtration.
    The condition on which components this diagram is cartesian
    can be described using the wall-and-chamber structure.

    More precisely, the associativity theorem states that we have an isomorphism
    \begin{equation*}
        \mathcal{X}_{\smash{\sigma_1 \uparrow \sigma_2}}^+ \longsimto
        \mathcal{X}_{\sigma_1}^+ \underset{\mathcal{X}_{\alpha_1}}{\times}
        \mathcal{X}_{\sigma_2}^+ \ ,
    \end{equation*}
    where $\sigma_1$ is a chamber in a special face~$\alpha_1$
    in the component lattice of~$\mathcal{X}$,
    and $\mathcal{X}_{\sigma_1}^+$ denotes the component
    in $\Filt (\mathcal{X})$
    corresponding to the points in the cones~$\sigma_1$,
    which stays constant within the cone by the constancy theorem;
    $\mathcal{X}_{\alpha_1}$~denotes the corresponding component
    in $\Grad (\mathcal{X})$;
    $\sigma_2$~is a cone containing~$\alpha_1$,
    and $\mathcal{X}_{\sigma_2}^+$ can be seen as a component
    in $\Filt (\mathcal{X}_{\alpha_1})$,
    so a component in $\Filt (\Grad (\mathcal{X}))$;
    and $\mathcal{X}_{\smash{\sigma_1 \uparrow \sigma_2}}^+$
    can be identified with a component of
    $\Filt (\Filt (\mathcal{X}))$.

    In the case when~$\mathcal{X}$ is the moduli stack
    of objects in an abelian category,
    the associativity theorem is the statement that
    for classes $\gamma_1, \gamma_2, \gamma_3 \in \uppi_0 (\mathcal{X})$,
    we have cartesian diagrams
    \begin{equation*}
        \begin{tikzcd}[column sep=2em, row sep={3.5em,between origins}]
            \mathcal{X}_{\gamma_1, \gamma_2, \gamma_3}^+
            \ar[r] \ar[d]
            \ar[dr, phantom, "\ulcorner", pos=.2]
            & \mathcal{X}_{\gamma_1} \times \mathcal{X}_{\gamma_2, \gamma_3}^+
            \ar[d]
            \\
            \mathcal{X}_{\gamma_1, \gamma_2 + \gamma_3}^+
            \ar[r]
            & \mathcal{X}_{\gamma_1} \times \mathcal{X}_{\gamma_2 + \gamma_3}
            \vphantom{^+} \rlap{\ ,}
        \end{tikzcd}
        \qquad
        \begin{tikzcd}[column sep=2em, row sep={3.5em,between origins}]
            \mathcal{X}_{\gamma_1, \gamma_2, \gamma_3}^+
            \ar[r] \ar[d]
            \ar[dr, phantom, "\ulcorner", pos=.2]
            & \mathcal{X}_{\gamma_1, \gamma_2}^+ \times \mathcal{X}_{\gamma_3}
            \ar[d]
            \\
            \mathcal{X}_{\gamma_1 + \gamma_2, \gamma_3}^+
            \ar[r]
            & \mathcal{X}_{\gamma_1 + \gamma_2} \times \mathcal{X}_{\gamma_3}
            \vphantom{^+} \rlap{\ ,}
        \end{tikzcd}
    \end{equation*}
    where $\mathcal{X}_{\gamma_1, \gamma_2, \gamma_3}^+$
    is the moduli stack of $3$-step filtrations
    with the graded pieces of the given classes in the given order,
    and similarly for the other notations~$\mathcal{X}_{\smash{(\ldots)}}^+$.

    In this abelian category case, these cartesian diagrams
    are important in the construction of
    various types of \emph{Hall algebras},
    and are used to prove the associativity of Hall algebras.
    Different flavours of Hall algebras include
    the original \emph{Hall algebras} studied by
    \textcite{ringel-1990-hall-1,ringel-1990-hall-2},
    \emph{motivic Hall algebras}
    constructed by \textcite{joyce-2007-configurations-ii},
    and \emph{cohomological Hall algebras}
    constructed by \textcite{kontsevich-soibelman-2011-coha},
    and later extended by \textcite{kinjo-park-safronov-coha}.

    For a general stack~$\mathcal{X}$,
    under mild conditions,
    we define a \emph{Hall category} of~$\mathcal{X}$
    using the component lattice,
    which captures the discrete data that is relevant for enumerative theories.
    The associativity theorem can be formulated as a functorial property
    related to the Hall category,
    as in \cref{thm-functorial-assoc}.
    This construction produces functors from the Hall category,
    which we regard as analogues of Hall algebras for general stacks.
\end{para}

\begin{para}[Real filtrations]
    \label{para-intro-real-filtrations}
    Our main results also have several applications
    outside Donaldson--Thomas theory,
    which we discuss here and in \cref{para-intro-semistable-reduction} below.

    In \cref{S:real_filtrations},
    we construct stacks of \emph{real filtrations} in a given stack~$\mathcal{X}$.
    These stacks parametrize filtrations in~$\mathcal{X}$
    that are weighted by real numbers,
    rather than by integers.
    These generalized filtrations naturally occur in many moduli problems,
    such as in K-stability, as in
    \textcite{chen-sun-wang-2018}
    and \textcite{blum-liu-xu-zhuang-2023-kahler-ricci},
    and in Bridgeland stability, as in
    \textcite[\S 6.5]{halpern-leistner-instability}, for example.

    We also define a version of stacks of filtrations
    associated to closed convex cones in~$\mathbb{R}^n$,
    with the cone~$\mathbb{R}_{\geq 0} \subset \mathbb{R}$
    corresponding to the stack of real filtrations mentioned above.
\end{para}

\begin{para}[Semistable reduction]
    \label{para-intro-semistable-reduction}
    Another related application of our main results
    generalizes the \emph{semistable reduction theorem} of
    \textcite[Theorem~B]{alper-halpern-leistner-heinloth-2023-moduli},
    which provides a general method of proving
    the properness of good moduli spaces in many moduli problems,
    given a \emph{$\Theta$-stratification}
    in the sense of \textcite[Definition~2.1.2]{halpern-leistner-instability}.

    In \cref{thm-semistable-reduction},
    we extend this result from the case of usual $\Theta$-stratifications
    to the case of \emph{real $\Theta$-stratifications},
    defined using stacks of real filtrations as above,
    and also to the case of \emph{cone $\Theta$-stratifications},
    also called \emph{higher $\Theta$-stratifications}
    in \textcite{odaka2024stability}.
    Our proof of this result takes a different approach than that of the closely related theorem \cite[Corollary~3.10, Theorem~3.18]{odaka2024stability}, and involves the constancy and finiteness theorems. Roughly speaking, we show that the real-weighted Harder--Narasimhan filtrations can be perturbed to $\mathbb{Z}$-weighted Harder--Narasimhan filtrations without affecting the underlying geometry of the $\Theta$-stratification, and then the original semistable reduction theorem applies. We expect this result to see applications in the theory of K-stability.
\end{para}

\begin{para}[Future directions]
    In the following papers \cite{part-ii,part-iii}
    of this series,
    based on the results of this paper,
    we will construct intrinsic Donaldson--Thomas invariants
    for general $(-1)$-shifted symplectic Artin stacks,
    extending the constructions of
    \textcite{joyce-song-2012-dt},
    \textcite{kontsevich-soibelman-motivic-dt},
    as well as the motivic enumerative invariants of
    \textcite{joyce-2006-configurations-i,joyce-2007-configurations-ii,joyce-2007-configurations-iii,joyce-2008-configurations-iv,joyce-2007-stack-functions}.

    We also hope that the formalism developed in this paper
    will be useful for extending other types of enumerative theories
    from the linear case of abelian categories to the non-linear case of more general stacks,
    such as Joyce's homological enumerative invariants
    \cite{joyce-wall-crossing},
    the theory of \emph{Vafa--Witten invariants}
    of \textcite{tanaka-thomas-2020-i,tanaka-thomas-2018-ii},
    and the theory of \emph{DT4~invariants}
    in \cite{cao-leung-dt4,borisov-joyce-2017,oh-thomas-2023-i,oh-thomas-ii},
    as mentioned in \cref{para-intro-overview}.
\end{para}

\needspace{2\baselineskip}
\begin{para}[Acknowledgements]
    We thank
    Dominic Joyce
    for helpful discussions related to this work,
    and for his comments on an earlier draft of this paper.
    We also thank
    Johan de~Jong,
    Yuji Odaka,
    Yukinobu Toda,
    Chenyang Xu,
    and Ziquan Zhuang,
    for useful discussions.

    C.~Bu would like to thank
    the Mathematical Institute, University of Oxford,
    for its support during the preparation of this paper.
    D.~Halpern-Leistner was supported by the NSF grants DMS-2052936 and DMS-1945478.
    A.~Ibáñez Núñez would like to thank the Mathematical Institute of the University of Oxford, the Newton Institute in Cambridge and the Department of Mathematics at Columbia University for their support.
    T.~Kinjo was supported by JSPS KAKENHI Grant Number 23K19007.
\end{para}

\begin{para}[Conventions]
    \label{assumption-stack-basic}
    Throughout this paper, we work under the following conventions:

    \begin{itemize}
        \item
            We work over a base algebraic space~$S$
            which is quasi-separated and locally noetherian.

        \item
            All schemes, algebraic spaces,
            and algebraic stacks, if not specified otherwise,
            are defined over~$S$, and are assumed to be
            quasi-separated and locally of finite type over~$S$.
            In the case of stacks, we further assume that
            they have affine stabilizers and separated inertia,
            unless otherwise stated.
            We denote by $\mathsf{St}_S$ the $2$-category of
            such algebraic stacks over~$S$.

        \item
            For a group algebraic space~$G$ acting on an algebraic space~$X$ over a base~$S$,
            we denote by $X/G$ the quotient stack, without the customary brackets.
    \end{itemize}
\end{para}

\resumeparabookmarks
\section{Formal lattices}

\subsection{Formal lattices}

\begin{para}
    We introduce the notion of \emph{formal lattices},
    which will be used to capture the combinatorial structure
    of the component lattice introduced in \cref{para-intro-component-lattice}.
\end{para}

\begin{para}[Formal lattices]
    \label{para-formal-lattices}
    Let $R$ be a commutative ring,
    which we will mostly take to be~$\mathbb{Z}$ or~$\mathbb{Q}$.
    Let $\mathsf{Lat} (R)$
    denote the category of free $R$-modules of finite rank,
    or \emph{$R$-lattices}.

    A \emph{formal $R$-lattice} is a functor
    \[
        X \colon \mathsf{Lat} (R)^\mathrm{op} \longrightarrow \mathsf{Set} \ .
    \]
    Denote by $\mathsf{fLat} (R) = \mathsf{Fun} (\mathsf{Lat} (R)^\mathrm{op}, \mathsf{Set})$
    the category of formal $R$-lattices.

    The \emph{underlying set} of a formal $R$-lattice~$X$
    is the set $|X| = X (R)$.
\end{para}

\begin{para}[Examples]
    \label{para-formal-lattices-examples}
    We give the following examples of formal lattices:

    \begin{enumerate}
        \item
            Any $R$-module~$M$ can be seen as a formal $R$-lattice
            via the Yoneda embedding, that is, by setting
            $M (F) = \Hom_R (F, M)$ for every free $R$-module~$F$ of finite rank.
            These formal $R$-lattices are precisely those
            that preserve finite products.
    \end{enumerate}
    The category of formal lattices has all limits and colimits,
    giving the following examples:

    \begin{enumerate}[resume]
        \item
            For a family of formal $R$-lattices $(X_i)_{i \in I}$,
            the disjoint union $\coprod_{i \in I} X_i$
            and the product $\prod_{i \in I} X_i$
            are again formal $R$-lattices.
        \item
            We are also allowed to glue formal $R$-lattices along arbitrary maps.
            For instance, one can glue two copies of $R$ along the origin,
            obtaining a formal $R$-lattice that looks like the set
            $\{ (r, s) \in R^2 \mid r = 0 \text{ or } s = 0 \}$.
        \item
            \label{item-eg-formal-lattice-quotient}
            If $X$ is a formal $R$-lattice,
            and $G$ is a group acting on $X$,
            then we can form the quotient $X / G$
            as a colimit in $\mathsf{fLat} (R)$,
            which assigns $(X/G)(F) = X(F)/G$
            for any $R$-lattice~$F$.
        \item
            To any abstract simplicial complex $(S,\Delta)$, meaning a finite
            set $S$ and a collection $\Delta$ of subsets such that any subset
            of a set $T \in \Delta$ also lies in $\Delta$, one can associate a
            formal lattice $X = \colim_{T \in \Delta} R^T$.
    \end{enumerate}
\end{para}

\begin{para}[Faces]
    \label{para-faces}
    Let $X$ be a formal $R$-lattice.
    The \emph{category of faces} of $X$ is the category $\mathsf{Face} (X)$ defined as follows:

    \begin{itemize}
        \item
            An object is a pair $(F, \alpha)$, called a \emph{face} of $X$,
            where $F$ is a free $R$-module of finite rank,
            and $\alpha \colon F \to X$ is a morphism in $\mathsf{fLat} (R)$.
        \item
            A morphism $(F', \alpha') \to (F, \alpha)$ is a morphism $f \colon F' \to F$
            in $\mathsf{Lat} (R)$ such that $\alpha' = \alpha \circ f$.
    \end{itemize}
    We sometimes abbreviate the pair $(F, \alpha)$ as~$\alpha$.

    Equivalently, $\mathsf{Face} (X)$ is the comma category
    $(\mathsf{Lat} (R) \downarrow X)$ taken in $\mathsf{fLat} (R)$.

    We say that a map $f$ between faces of $X$ is \emph{injective} or \emph{surjective} if
    the underlying map between $R$-lattices is injective or surjective, respectively.

    A face $\alpha \colon F \to X$ is \emph{non-degenerate},
    if it does not factor through
    a free $R$-module of rank less than that of~$F$.
    The \emph{dimension} of such a non-degenerate face is the rank of~$F$.
    Denote by
    \begin{equation*}
        \mathsf{Face}^{\mathrm{nd}} (X) \subset \mathsf{Face} (X)
    \end{equation*}
    the full subcategory of non-degenerate faces.
    In particular, if~$R$ is a principal ideal domain,
    then all morphisms in $\mathsf{Face}^{\mathrm{nd}} (X)$
    are injective.

    A map $f \colon X \to Y$ of formal $R$-lattices
    is \emph{unramified} if it preserves non-degenerate faces.
\end{para}

\begin{example}
    Let~$K$ be a field,
    $F$ a finite-dimensional $K$-vector space,
    and $G$ a group acting on~$F$.
    Consider the formal $K$-lattice $X=F/G$ as in
    \cref{para-formal-lattices-examples}~\cref{item-eg-formal-lattice-quotient}.

    The category $\mathsf{Face}^{\mathrm{nd}} (F / G)$ of non-degenerate faces can
    be explicitly described as follows.
    Its objects are linear subspaces $E \subset F$.
    A morphism between objects $E, E' \subset F$ is a linear map $E\to E'$ that is induced by the action of some element $g \in G$.
\end{example}

\begin{para}[Rational cones]
    \label{para-cones}
    For a finite-dimensional $\mathbb{Q}$-vector space~$F$,
    a (\emph{polyhedral}) \emph{cone} in~$F$ is a subset of the form
    \begin{equation*}
        C = \{ a_1 v_1 + \cdots + a_n v_n \mid a_1, \dotsc, a_n \in \mathbb{Q}_{\geq 0} \}
    \end{equation*}
    for some elements~$v_1, \ldots, v_n \in F$.
    We say that~$C$ is of \emph{full dimension} if it spans~$F$.

    Define a \emph{rational cone} to be a commutative monoid~$C$
    which is isomorphic to a polyhedral cone
    in a finite-dimensional $\mathbb{Q}$-vector space~$F$.
    The vector space~$F$ can then be recovered
    as the groupification~$F_C$ of~$C$.
    Define the \emph{dimension} of~$C$ to be
    $\dim C = \dim F_C$.
    Write $\mathsf{Cone}_\mathbb{Q}$ for the category of rational cones,
    where the morphisms are monoid homomorphisms.

    For a formal $\mathbb{Q}$-lattice~$X$
    and a rational cone~$C$,
    define a \emph{morphism} $C \to X$
    to be a morphism of formal $\mathbb{Q}$-lattices $F_C \to X$,
    although we think of this as restricting the latter morphism
    to the cone $C \subset F_C$.
    In this way, $X$ also defines a functor
    $\mathsf{Cone}_{\smash{\mathbb{Q}}}^\mathrm{op} \to \mathsf{Set}$.

    Define the \emph{category of cones} in~$X$,
    denoted by $\mathsf{Cone} (X)$, as follows:

    \begin{itemize}
        \item
            An object is a pair~$(C, \sigma)$,
            where~$C$ is a rational cone,
            and $\sigma \colon C \to X$ is a morphism.
        \item
            A morphism $(C, \sigma) \to (C', \sigma')$
            is a morphism $f \colon C \to C'$
            such that $\sigma' \circ f = \sigma$.
    \end{itemize}
    We sometimes abbreviate $(C, \sigma)$ as~$\sigma$.
    There is a natural inclusion
    $\mathsf{Face} (\mathcal{X}) \hookrightarrow \mathsf{Cone} (\mathcal{X})$.

    For cones~$(C, \sigma), (C', \sigma')$ in~$X$,
    we use the notation $\sigma \subset \sigma'$
    to mean that the morphism~$\sigma$
    factors through an injective map $C \to C'$,
    and that we fix a choice of such a factorization.

    The \emph{span} of a cone $(C, \sigma)$ in~$X$
    is the face $\operatorname{span} (C, \sigma) = (F, \alpha)$,
    where $F = F_C$ is the groupification of~$C$,
    and~$\alpha \colon F \to X$ is the morphism that underlies~$\sigma$.
    This defines a functor
    $\operatorname{span} \colon \mathsf{Cone} (X) \to \mathsf{Face} (X)$,
    which is left adjoint to the inclusion
    $\mathsf{Face} (X) \hookrightarrow \mathsf{Cone} (X)$.

    Define the category of \emph{non-degenerate cones} in~$X$
    as the full subcategory
    \begin{equation*}
        \mathsf{Cone}^{\mathrm{nd}} (X) \subset \mathsf{Cone} (X)
    \end{equation*}
    consisting of objects~$(C, \sigma)$
    such that $\mathrm{span} (C, \sigma)$ is a non-degenerate face.
\end{para}

\begin{para}[Integral cones]
    \label{para-integral-cones}
    Most of \cref{para-cones}
    can also be adapted to formal $\mathbb{Z}$-lattices, as follows.

    Define an \emph{integral cone}
    to be a commutative monoid~$\Sigma$
    which is isomorphic to $\Lambda \cap C$
    for a rational cone~$C$
    and a $\mathbb{Z}$-lattice~$\Lambda \subset F_C$ of full rank,
    where the intersection is taken in~$F_C$.
    The lattice~$\Lambda$ can then be recovered
    as the groupification~$\Lambda_\Sigma$ of~$\Sigma$,
    and the cone~$C$ as the rationalization
    $\Sigma_\mathbb{Q} = \mathbb{Q}_{\geq 0} \cdot \Sigma \subset \Lambda_\Sigma \otimes \mathbb{Q}$.
    Define the \emph{rank} of~$\Sigma$ to be
    $\operatorname{rk} \Sigma = \operatorname{rk} \Lambda_\Sigma$.
    Write $\mathsf{Cone}_\mathbb{Z}$ for the category of integral cones,
    where the morphisms are monoid homomorphisms.

    Note that in the literature,
    integral cones are also called \emph{normal affine monoids};
    see, for example, \textcite[\S2.B]{gubeladze-winfried-2009-polytopes}

    Analogously to \cref{para-cones},
    for a formal $\mathbb{Z}$-lattice~$X$,
    we can define a notion of morphisms from integral cones to~$X$,
    a category $\mathsf{Cone} (X)$ of cones in~$X$,
    and a full subcategory
    $\mathsf{Cone}^{\mathrm{nd}} (X) \subset \mathsf{Cone} (X)$
    of non-degenerate cones in~$X$.
\end{para}

\begin{para}[Restriction and extension of scalars]
    \label{para-restriction}

    Let $\phi\colon R\to S$ be a homomorphism of commutative rings.
    The \emph{restriction of scalars} functor
    \[
        (-) \restr^S_R \colon
        \mathsf{fLat} (S) \longrightarrow \mathsf{fLat} (R)
    \]
    is defined by precomposition along
    $(-) \otimes_R S \colon \mathsf{Lat}(R) \to \mathsf{Lat}(S)$.
    Restriction of scalars has a left adjoint
    $(-) \otimes_R S \colon \mathsf{fLat}(R) \to \mathsf{fLat}(S)$,
    given by
    \[
        X \otimes_R S = \colim_{(F, \alpha) \in \mathsf{Face} (X)} {} (F \otimes_R S) \ .
    \]
    We call this left adjoint \emph{extension of scalars},
    and we will often abbreviate $X_S = X \otimes_R S$.
\end{para}

\begin{lemma}\label{lemma: extension of scalars for modules is usual tensor product}
    If\/ $M$ is an $R$-module, regarded as a formal $R$-lattice,
    then the extension of scalars $M \otimes_R S$
    agrees with the ordinary tensor product.
\end{lemma}

\begin{proof}
    Since $M$ is a module, the category $\mathsf{Face} (M)$ has
    coproducts, so it is a sifted category.
    This implies that
    $M\otimes_R S\colon \mathsf{fLat}(S)^\mathrm{op} \to \mathsf{Set}$
    preserves finite products, hence it is an $S$-module. Since this applies to
    any $R$-module $M$, the functor
    $(-)\otimes_R S\colon \mathsf{fLat}(R)\to \mathsf{fLat}(S)$
    restricts to a functor $\mathsf{Mod}(R)\to \mathsf{Mod}(S)$ which is a left
    adjoint to the restriction of scalars functor
    $(-)\restr^S_R\colon \mathsf{Mod}(S)\to \mathsf{Mod}(R)$
    and therefore must be the ordinary tensor product $(-)\otimes_R S$.
\end{proof}

\begin{para}[Geometric realization]
    \label{para-geometric-realization}
    For a formal $\mathbb{R}$-lattice~$X$,
    its \emph{geometric realization}
    is a topological space with the underlying set
    $|X| = X (\mathbb{R})$, defined by
    \begin{equation}
        \label{eq-def-geometric-realization}
        |X| = \colim_{(F,\alpha) \in \mathsf{Face}(X)} F \ ,
    \end{equation}
    equipped with the colimit topology of
    the euclidean topologies on the $\mathbb{R}$-vector spaces~$F$.

    In particular, for a formal $\mathbb{Z}$-~or $\mathbb{Q}$-lattice~$X$,
    writing~$X_\mathbb{R}$ for the extension of scalars to~$\mathbb{R}$
    as in \cref{para-restriction},
    we have a canonical homeomorphism
    \begin{equation}
        \label{eq-real-geometric-realization}
        |X_\mathbb{R}| \simeq \colim_{(F,\alpha) \in \mathsf{Face}(X)} {} (F \otimes \mathbb{R}) \ ,
    \end{equation}
    where each $F \otimes \mathbb{R}$ is equipped with the euclidean topology.
\end{para}

\subsection{Face and cone arrangements}

\begin{para}
    In this section, we introduce the notions of
    \emph{face arrangements} and \emph{cone arrangements}
    in a formal lattice,
    which will be used to describe the wall-and-chamber structure
    on the component lattice introduced in \cref{para-intro-component-lattice}.
    The former notion is a generalization of
    the notion of a hyperplane arrangement in a vector space.
\end{para}

\begin{para}[Hyperplane arrangements]
    \label{para-hyperplane-arrangements}
    To fix terminology, we make the following definitions.

    For a finite-dimensional $\mathbb{Q}$-vector space~$F$,
    a~\emph{hyperplane arrangement} in~$F$ is a finite set
    $\Phi = \{ H_1, \ldots, H_n \}$ of hyperplanes in~$F$,
    where each $H_i \subset F$ is a linear subspace of codimension~$1$.

    A~\emph{flat} of~$\Phi$ is a linear subspace of~$F$
    of the form $\bigcap_{i \in I} H_i$ for some subset $I \subset \{ 1, \ldots, n \}$.
    The whole vector space $F$ is always a flat, realized as the empty intersection.

    A~\emph{chamber} of~$\Phi$ is a cone $C \subset F$ of full dimension,
    such that the interior of its closure in~$F_\mathbb{R}$
    is a connected component of
    $F_\mathbb{R} \setminus \bigcup_{H \in \Phi} H_\mathbb{R}$.
\end{para}

\begin{para}[Face arrangements]
    \label{para-face-arrangement}
    Let $X$ be a formal $\mathbb{Q}$-lattice.

    A \emph{face arrangement}~$\Phi$ in~$X$
    is a full subcategory
    \begin{equation*}
        \Phi \subset \mathsf{Face}^\mathrm{nd} (X) \ ,
    \end{equation*}
    such that the inclusion functor
    $\Phi \hookrightarrow \mathsf{Face} (X)$
    admits a left adjoint
    $(-)^\Phi \colon \mathsf{Face} (X) \to \Phi$,
    called the \emph{closure} in~$\Phi$.

    In this case, for each face $\alpha \in \mathsf{Face} (X)$,
    the adjunction unit $\alpha \to \alpha^\Phi$
    is universal among morphisms from~$\alpha$
    to faces in~$\Phi$, so that~$\alpha^\Phi$
    can be seen as the smallest face in~$\Phi$ containing~$\alpha$.

    For example, a face arrangement
    in a finite-dimensional $\mathbb{Q}$-vector space~$F$
    is a collection of subspaces of~$F$,
    containing~$F$ and closed under intersection.

    For any face arrangement~$\Phi$ in~$X$,
    for each face~$(F, \alpha) \in \Phi$,
    the overcategory $\Phi_{/ (F, \alpha)}$
    is a face arrangement in~$F$.
\end{para}

\begin{para}[Cone arrangements]
    \label{para-cone-arrangement}
    Let~$X$ be a formal $\mathbb{Q}$-lattice.

    A \emph{cone arrangement} in~$X$
    is a full subcategory
    \begin{equation*}
        \Psi \subset \mathsf{Cone}^\mathrm{nd} (X) \ ,
    \end{equation*}
    such that the inclusion functor
    $\Psi \hookrightarrow \mathsf{Cone} (X)$
    admits a left adjoint
    $(-)^\Psi \colon \mathsf{Cone} (X) \to \Psi$,
    called the \emph{closure} in~$\Psi$.

    For example, every face arrangement in~$X$
    is also a cone arrangement in~$X$.
    A finite cone arrangement in a finite-dimensional $\mathbb{Q}$-vector space~$F$
    is a finite collection of polyhedral cones in~$F$,
    containing~$F$ and closed under intersection.

    For any cone arrangement~$\Psi$ in~$X$,
    for each cone~$(C, \sigma) \in \Psi$,
    the collection of subcones of~$C$ that belong to~$\Psi$
    is closed under finite intersections,
    which follows from the adjunction property.
\end{para}

\begin{para}[Examples]
    \label{para-face-arrangement-examples}
    We give some examples
    of face and cone arrangements.

    \begin{enumerate}
        \item
            \label{item-eg-hyperplane-arrangement}
            Let~$F$ be a finite-dimensional $\mathbb{Q}$-vector space,
            equipped with a hyperplane arrangement.
            Then its flats form a face arrangement in~$F$,
            and the cones which are
            intersections of closed half-spaces cut out by the hyperplanes
            form a cone arrangement in~$F$.

        \item
            Let $F$ be a finite-dimensional $\mathbb{Q}$-vector space,
            acted on by a group~$G$,
            and consider the formal $\mathbb{Q}$-lattice~$X = F / G$ defined in
            \cref{para-formal-lattices-examples}~\cref{item-eg-formal-lattice-quotient}.
            Then a face (or~cone) arrangement in~$X$ is the same data as
            a $G$-invariant face (or~cone) arrangement in~$F$,
            such that for any $g \in G$,
            the fixed subspace $F^g \subset F$ belongs to the arrangement.

        \item
            Let $X = \mathbb{Q} \cup_{\{ 0 \}} \mathbb{Q}$
            be the formal $\mathbb{Q}$-lattice
            glued from two copies of~$\mathbb{Q}$ along the origin.
            Then~$X$ has a unique face arrangement,
            consisting of the two maximal non-degenerate faces
            and their intersection.
            It has precisely $16$ different cone arrangements,
            where each of the four rays can be optionally included.
    \end{enumerate}
\end{para}

\section{The component lattice}

\subsection{Graded and filtered points}

\begin{para}
    For an algebraic stack~$\mathcal{X}$,
    following \textcite{halpern-leistner-instability},
    we introduce the \emph{stack of graded points}
    and the \emph{stack of filtered points} of~$\mathcal{X}$,
    denoted by $\Grad (\mathcal{X})$ and $\Filt (\mathcal{X})$, respectively.
    These stacks have seen many important applications
    in enumerative geometry.

    For example, if $\mathcal{X}$ is the moduli stack of objects
    in an abelian category~$\mathcal{A}$,
    then $\Grad (\mathcal{X})$ gives the moduli stack of
    $\mathbb{Z}$-graded objects in~$\mathcal{A}$,
    and $\Filt (\mathcal{X})$ gives the moduli stack of
    $\mathbb{Z}$-filtered objects in~$\mathcal{A}$.
\end{para}

\begin{para}[Graded and filtered points]
    \label{para-grad-filt}
    Let $\mathcal{X}$ be an algebraic stack over~$S$
    as in \cref{assumption-stack-basic},
    and let $n \geq 0$ be an integer.

    The \emph{stack of\/ $\mathbb{Z}^n$-graded points}
    and the \emph{stack of\/ $\mathbb{Z}^n$-filtered points} of~$\mathcal{X}$
    are the mapping stacks
    \begin{align*}
        \Grad^n (\mathcal{X})
        & = \mathrm{Map} ( \mathrm{B} \mathbb{G}_\mathrm{m}^n, \mathcal{X} ) \ ,
        \\
        \Filt^n (\mathcal{X})
        & = \mathrm{Map} ( \Theta^n, \mathcal{X}) \ ,
    \end{align*}
    where $\mathbb{G}_\mathrm{m}$ is defined over~$\Spec \Z$,
    and $\Theta = \mathbb{A}^1 / \mathbb{G}_\mathrm{m}$
    is the quotient stack of~$\mathbb{A}^1$ by the scaling action of~$\mathbb{G}_\mathrm{m}$,
    also defined over $\Spec \Z$.

    The stacks $\Grad^n (\mathcal{X})$ and $\Filt^n (\mathcal{X})$
    are again algebraic stacks over~$S$
    satisfying the conditions in \cref{assumption-stack-basic},
    by \textcite[Proposition~1.1.2]{halpern-leistner-instability}. Separatedness of the inertia of $\Grad^n (\mathcal{X})$ and $\Filt^n (\mathcal{X})$ follows from that of the inertia of $\mathcal X$ and representability of the maps $\mathrm{tot}\colon \Grad^n (\mathcal{X}) \to \mathcal{X}$ and $\mathrm{ev}_1\colon \Filt^n (\mathcal{X}) \to \mathcal{X}$, which holds by \cite[Lemma~1.1.5 and Proposition~1.1.13]{halpern-leistner-instability}.

    We write $\Grad (\mathcal{X}) = \Grad^1 (\mathcal{X})$
    and $\Filt (\mathcal{X}) = \Filt^1 (\mathcal{X})$,
    and call them the \emph{stack of graded points}
    and the \emph{stack of filtered points} of~$\mathcal{X}$, respectively. 
\end{para}

\begin{para}[Induced morphisms]
    \label{para-grad-filt-morphisms}
    Consider the morphisms
    \begin{equation*}
        \begin{tikzcd}
            \mathrm{B} \mathbb{G}_\mathrm{m}^n
            \ar[shift left=0.5ex, r, "0"]
            &
            \Theta^n
            \ar[shift left=0.5ex, l, "\mathrm{pr}"]
            &
            * \vphantom{^n} \ ,
            \ar[shift left=0.5ex, l, "1"]
            \ar[shift right=0.5ex, l, "0"']
            \ar[ll, bend right, start anchor=north west, end anchor=north east, looseness=.8]
        \end{tikzcd}
    \end{equation*}
    where~$*$ denotes $\Spec \Z$,
    the map~$\mathrm{pr}$ is induced by the projection $\mathbb{A}^n \to *$,
    and~$1$ denotes the inclusion as the point $(1, \ldots, 1)$.
    These induce morphisms of stacks
    \begin{equation*}
        \begin{tikzcd}
            \Grad^n (\mathcal{X})
            \ar[rr, bend left, start anchor=north east, end anchor=north west, looseness=.6, "\smash{\mathrm{tot}}"]
            \ar[r, shift right=0.5ex, "\mathrm{sf}"']
            &
            \Filt^n (\mathcal{X})
            \ar[l, shift right=0.5ex, "\mathrm{gr}"']
            \ar[r, shift left=0.5ex, "\mathrm{ev}_0"]
            \ar[r, shift right=0.5ex, "\mathrm{ev}_1"']
            &
            \mathcal{X} \rlap{ ,}
        \end{tikzcd}
    \end{equation*}
    where the notations `$\mathrm{gr}$', `$\mathrm{sf}$', and `$\mathrm{tot}$' stand for
    the \emph{associated graded point},
    the \emph{split filtration},
    and the \emph{total point}, respectively.

    We list the following properties of these morphisms, which follow from
    \textcite[Lemma~1.1.5, Proposition~1.1.13, and Lemma~1.3.8]{halpern-leistner-instability}:

    \begin{enumerate}[midpenalty=10000]
        \item
            The morphisms $\mathrm{tot}$
            and $\mathrm{ev}_1$ are representable.

        \item
            The morphism $\mathrm{gr}$ is quasi-compact.

        \item
            The morphisms $\mathrm{gr}$ and $\mathrm{sf}$ induce an isomorphism
            $\uppi_0 (\Grad^n (\mathcal{X})) \simeq \uppi_0 (\Filt^n (\mathcal{X}))$.
    \end{enumerate}
\end{para}

\begin{example}[Quotient stacks]
    \label{eg-grad-quotient-stack}
    The stacks of graded and filtered points of a quotient stack
    can be described explicitly, following
    \textcite[Theorems~1.4.7 and~1.4.8]{halpern-leistner-instability}.

    Let~$\mathcal{X} = U / G$ be a quotient stack,
    where~$U$ is an algebraic space over~$S$,
    acted on by a smooth affine group scheme~$G$ over~$S$,
    and assume one of the following:
    \begin{itemize}
        \item
            $G = \mathrm{GL} (m)$ for some~$m$, or
        \item
            $S = \Spec K$ for a field~$K$,
            and~$G$ has a split maximal torus.
    \end{itemize}
    For a morphism of group schemes
    $\lambda \colon \mathbb{G}_\mathrm{m}^n \to G$,
    define the \emph{Levi subgroup}
    and the \emph{parabolic subgroup} of~$G$
    associated to~$\lambda$ by
    \begin{align*}
        L_\lambda & =
        \{ g \in G \mid g = \lambda (t) \, g \, \lambda (t)^{-1} \text{ for all } t \} \ ,
        \\
        P_\lambda & =
        \{ g \in G \mid \lim_{t \to 0} \lambda (t) \, g \, \lambda (t)^{-1} \text{ exists} \} \ ,
    \end{align*}
    respectively. Define the
    \emph{fixed locus} and the \emph{attractor} associated to~$\lambda$ by
    \begin{align*}
        U^{\lambda} & =
        \Map^{\mathbb{G}_\mathrm{m}^n} (*, U) \ ,
        \\
        U^{\lambda, +} & =
        \Map^{\mathbb{G}_\mathrm{m}^n} (\mathbb{A}^n, U) \ ,
    \end{align*}
    where $\Map^{\mathbb{G}_\mathrm{m}^n} (-, -)$
    denotes the $\mathbb{G}_\mathrm{m}^n$-equivariant mapping space,
    and~$\mathbb{G}_\mathrm{m}^n$ acts on~$U$ via~$\lambda$,
    and on $\mathbb{A}^n$ by scaling each coordinate.
    These are again algebraic spaces over~$S$.
    See \textcite[\S1]{drinfeld-gaitsgory-2014-braden}
    or \textcite[\S1.4]{halpern-leistner-instability} for more details.

    The $G$-action on~$U$ induces a $P_\lambda$-action on~$U^{\lambda, +}$
    and an $L_\lambda$-action on~$U^{\lambda}$.
    Moreover, by \cite[Theorems~1.4.7 and~1.4.8]{halpern-leistner-instability},
    we have
    \begin{align}
        \Grad^n (\mathcal{X}) & \simeq
        \coprod_{\lambda \colon \mathbb{G}_\mathrm{m}^n \to G} U^{\lambda} / L_\lambda \ , \label{eq:Grad-quotient-stack}
        \\
        \Filt^n (\mathcal{X}) & \simeq
        \coprod_{\lambda \colon \mathbb{G}_\mathrm{m}^n \to G} U^{\lambda, +} / P_\lambda \ ,  \label{eq:Filt-quotient-stack}
    \end{align}
    where the disjoint unions are taken over all
    conjugacy classes of maps~$\lambda$.

    This example will be continued in \cref{eg-filt-cones-quotient-stack},
    where we discuss a similar description for
    a generalized version of filtered points.
\end{example}

\begin{para}[Coordinate-free notation]
    \label{para-grad-lambda}
    We introduce a coordinate-free notation
    for $\Grad^n (\mathcal{X})$.

    For a free $\mathbb{Z}$-module~$\Lambda$ of finite rank,
    let $T_\Lambda = \mathrm{Spec} \, \mathbb{Z} [\Lambda^\vee] \simeq \Gm^{\rk \Lambda}$
    be the torus with cocharacter lattice~$\Lambda$.
    Define the \emph{stack of $\Lambda^\vee$-graded points} of~$\mathcal{X}$ by
    \begin{equation*}
        \Grad^\Lambda (\mathcal{X}) =
        \mathrm{Map} (\mathrm{B} T_\Lambda, \mathcal{X}) \ .
    \end{equation*}
    This construction is contravariant in~$\Lambda$.
    We have an isomorphism
    $\Grad^\Lambda (\mathcal{X}) \simeq \Grad^{\rk \Lambda} (\mathcal{X})$
    upon choosing an isomorphism $\Lambda \simeq \mathbb{Z}^{\rk \Lambda}$.

    There is also a coordinate-free version of the stack
    $\Filt^n (\mathcal{X})$,
    which we discuss in \cref{subsec-filt-cones}.
\end{para}

\begin{para}[Rational graded points]
    \label{para-rational-graded-points}
    We introduce stacks of \emph{$\mathbb{Q}^n$-graded points} of~$\mathcal{X}$,
    denoted by $\Grad^n_\mathbb{Q} (\mathcal{X})$,
    as in \textcite[\S2.2]{ibanez-nunez-filtrations}.
    For example, if~$\mathcal{X}$ is the moduli stack of objects
    in an abelian category~$\mathcal{A}$,
    then $\Grad^n_\mathbb{Q} (\mathcal{X})$
    parametrizes \emph{$\mathbb{Q}^n$-graded objects} in~$\mathcal{A}$,
    that is, families $(x_v)_{v \in \mathbb{Q}^n}$ of objects in~$\mathcal{A}$
    such that $x_v = 0$ for all but finitely many~$v$.

    For a finite-dimensional $\mathbb{Q}$-vector space~$F$,
    define the \emph{stack of\/ $F^\vee$-graded points} of~$\mathcal{X}$ by
    \begin{equation*}
        \Grad^F (\mathcal{X}) =
        \colim_{\Lambda \subset F} \Grad^\Lambda (\mathcal{X}) \ ,
    \end{equation*}
    where the colimit is taken over all free
    $\mathbb{Z}$-submodules~$\Lambda \subset F$
    of full rank.
    We have a morphism
    $\mathrm{tot} \colon \Grad^F (\mathcal{X}) \to \mathcal{X}$,
    defined as the colimit of the morphisms
    $\mathrm{tot} \colon \Grad^\Lambda (\mathcal{X}) \to \mathcal{X}$.

    In particular, we write
    $\Grad_\mathbb{Q}^n (\mathcal{X}) = \Grad^{\mathbb{Q}^n} (\mathcal{X})$
    and $\Grad_\mathbb{Q} (\mathcal{X}) = \Grad^\mathbb{Q} (\mathcal{X})$.

    This construction does not produce essentially new stacks,
    since by \cref{lemma-grad-surj} below,
    every morphism in the colimit diagram is an open and closed immersion,
    so each component of $\Grad^n_\mathbb{Q} (\mathcal{X})$
    isomorphic to a component of $\Grad^n (\mathcal{X})$.
    In other words, $\Grad^n_\mathbb{Q} (\mathcal{X})$
    is just $\Grad^n (\mathcal{X})$ with each component duplicated many times.
    Our main reason for introducing these stacks is so that their components
    correspond to points in the \emph{rational component lattice}
    of~$\mathcal{X}$, as we define in \cref{para-component-lattice} below,
    which is a central object of study in this paper.
\end{para}

\begin{lemma}
    \label{lemma-grad-surj}
    Let $\mathcal{X}$ be a stack as in \cref{assumption-stack-basic}.
    Let $\Lambda_1,\Lambda_2$ be free $\Z$-modules of finite rank,
    and let $\pi \colon \Lambda_1 \to \Lambda_2$
    be a homomorphism whose image is of full rank.

    Then the induced morphism
    \begin{equation*}
        \pi^* \colon \Grad^{\Lambda_2} (\mathcal{X})
        \longrightarrow \Grad^{\Lambda_1} (\mathcal{X})
    \end{equation*}
    is an open and closed immersion.
\end{lemma}

\begin{proof}
    This is a special case of \cref{thm-filt-sigma-image} below.
\end{proof}

\begin{lemma}
    \label{lemma-grad-preserves-qc-morphisms}
    Let $\mathcal{X}$ and $\mathcal{Y}$ be stacks as in \cref{assumption-stack-basic},
    and let $f \colon \mathcal{Y} \to \mathcal{X}$ be a quasi-compact representable morphism.

    Then the induced morphisms
    $\Grad^\Lambda (\mathcal{Y}) \to \Grad^\Lambda (\mathcal{X})$ and
    $\Grad^F (\mathcal{Y}) \to \Grad^F (\mathcal{X})$
    are also quasi-compact and representable,
    for any $\Lambda \in \mathsf{Lat} (\mathbb{Z})$
    and $F \in \mathsf{Lat} (\mathbb{Q})$.
\end{lemma}

\begin{proof}
    For the statement on $\Grad^\Lambda$,
    it is enough to prove that for any quasi-compact scheme~$Z$
    and any morphism $Z \to \Grad^\Lambda (\mathcal{X})$,
    the pullback $Z \times_{\Grad^\Lambda (\mathcal{X})} \Grad^\Lambda (\mathcal{Y})$
    is a quasi-compact algebraic space.
    Let $Z \times \mathrm{B} T_\Lambda \to \mathcal{X}$
    be the corresponding morphism.
    Since~$f$ is quasi-compact and representable, we have
    $(Z \times \mathrm{B} T_\Lambda) \times_{\mathcal{X}} \mathcal{Y}
    \simeq U / T_\Lambda$ for a quasi-compact algebraic space~$U$.
    Then we have
    $Z \times_{\Grad^\Lambda (\mathcal{X})} \Grad^\Lambda (\mathcal{Y}) \simeq
    U^{\smash{T_\Lambda}}$.
    By \textcite[Proposition~1.4.1]{halpern-leistner-instability},
    $U^{\smash{T_\Lambda}}$ is a closed subspace of~$U$,
    and hence quasi-compact.

    For the statement on $\Grad^F$,
    the morphism in question
    is isomorphic to the one on $\Grad^\Lambda$ on each component,
    so it is also quasi-compact.
\end{proof}

\subsection{The component lattice}

\begin{para}
    For an algebraic stack~$\mathcal{X}$, we show that the set
    $\uppi_0 (\Grad (\mathcal{X})) \simeq \uppi_0 (\Filt (\mathcal{X}))$
    carries a natural structure of a formal $\mathbb{Z}$-lattice,
    called the \emph{component lattice} of~$\mathcal{X}$,
    denoted by~$\mathrm{CL} (\mathcal{X})$.

    This formal lattice, together with its rationalization
    $\mathrm{CL}_\mathbb{Q} (\mathcal{X})$,
    is a central object of study in this paper,
    and will be used to construct enumerative invariants of~$\mathcal{X}$
    in subsequent works.
\end{para}

\begin{para}[The component lattice]
    \label{para-component-lattice}
    Let $\mathcal{X}$ be a stack as in \cref{assumption-stack-basic}.

    The \emph{component lattice} of~$\mathcal{X}$
    is the formal $\mathbb{Z}$-lattice
    $\mathrm{CL} (\mathcal{X})$ defined by setting
    \[
        \mathrm{CL} (\mathcal{X}) (\Lambda) =
        \uppi_0 (\Grad^\Lambda(\mathcal{X}))
    \]
    for all $\Lambda \in \mathsf{Lat} (\mathbb{Z})$,
    where $\Grad^\Lambda (\mathcal{X})$ is defined in \cref{para-grad-lambda}.
    This formal lattice structure can be thought of as
    an extra combinatorial structure
    on the set
    $|\mathrm{CL} (\mathcal{X})| \simeq \uppi_0 (\Grad (\mathcal{X}))$.

    Let $\mathrm{CL}_\mathbb{Q} (\mathcal{X}) =
    \mathrm{CL} (\mathcal{X}) \otimes_\mathbb{Z} \mathbb{Q}$
    be the rationalization of~$\mathrm{CL} (\mathcal{X})$,
    as in \cref{para-restriction}.
    Equivalently, for $F \in \mathsf{Lat} (\mathbb{Q})$, we have
    \begin{equation*}
        \mathrm{CL}_\mathbb{Q} (\mathcal{X}) (F) \simeq
        \uppi_0 (\Grad^F (\mathcal{X})) \ ,
    \end{equation*}
    where $\Grad^F (\mathcal{X})$ is defined in \cref{para-rational-graded-points}.
    In particular, we have
    $|\mathrm{CL}_\mathbb{Q} (\mathcal{X})| \simeq \uppi_0 (\Grad_\mathbb{Q} (\mathcal{X}))$.

    Similarly, write
    $\mathrm{CL}_\mathbb{R} (\mathcal{X}) =
    \mathrm{CL} (\mathcal{X}) \otimes_\mathbb{Z} \mathbb{R}$
    for the extension of scalars defined in \cref{para-restriction}.

    We introduce shorthand notations for the categories
    \begin{alignat*}{2}
        \mathsf{Face} (\mathcal{X}) & =
        \mathsf{Face} (\mathrm{CL}_\mathbb{Q} (\mathcal{X})) \ ,
        \qquad &
        \mathsf{Face}^\mathrm{nd} (\mathcal{X}) & =
        \mathsf{Face}^\mathrm{nd} (\mathrm{CL}_\mathbb{Q} (\mathcal{X})) \ ,
        \\
        \mathsf{Cone} (\mathcal{X}) & =
        \mathsf{Cone} (\mathrm{CL}_\mathbb{Q} (\mathcal{X})) \ ,
        \qquad &
        \mathsf{Cone}^\mathrm{nd} (\mathcal{X}) & =
        \mathsf{Cone}^\mathrm{nd} (\mathrm{CL}_\mathbb{Q} (\mathcal{X})) \ .
    \end{alignat*}
    where the right-hand sides are defined in
    \crefrange{para-faces}{para-cones}.
    
    See
    \crefrange{eg-component-lattice-linear-quotient}{eg-component-lattice-linear-moduli}
    below for examples of component lattices.
    Also, in \cref{lemma-cl-unramified} below,
    we show that representable morphisms of stacks
    induce unramified maps on component lattices.
\end{para}

\begin{para}[The stack~\texorpdfstring{$\mathcal{X}_\alpha$}{X\_α}]
    \label{para-notation-x-alpha}
    For a face $\alpha \colon F \to \mathrm{CL}_\Q (\mathcal{X})$,
    define
    \begin{equation*}
        \mathcal{X}_\alpha \subset \Grad^F(\mathcal{X})
    \end{equation*}
    as the connected component corresponding to~$\alpha$.
    This is contravariant in
    $\alpha \in \mathsf{Face} (\mathcal{X})$.

    For an element $\lambda \in |\mathrm{CL}_\Q (\mathcal{X})|$,
    we also write
    $\mathcal{X}_\lambda \subset \Grad_\mathbb{Q} (\mathcal{X})$
    for the connected component corresponding to~$\lambda$.
    This can be seen as the special case of~$\mathcal{X}_\alpha$
    for the face $\alpha = \mathbb{Q} \cdot \lambda$.

    Additionally, for a morphism $\mathcal{Y} \to \mathcal{X}$ of algebraic stacks,
    we also write $\mathcal{Y}_\alpha \subset \Grad^F (\mathcal{Y})$
    for the preimage of~$\mathcal{X}_\alpha$ along the induced morphism
    $\Grad^F(\mathcal{Y})\to \Grad^F(\mathcal{X})$,
    which is an open and closed substack of~$\Grad^F (\mathcal{Y})$.

    We have the following rational version of \cref{lemma-grad-surj}.
\end{para}

\begin{lemma}
    \label{lemma-grad-surj-rational}
    Let~$\mathcal{X}$ be a stack as in \cref{assumption-stack-basic},
    and let $\alpha \colon F \to \mathrm{CL}_\Q (\mathcal{X})$ be a face.
    Then for any surjection $\pi \colon F' \to F$,
    writing $\beta = \alpha \circ \pi$,
    we have a canonical isomorphism
    \begin{equation*}
        \mathcal{X}_\alpha \longsimto \mathcal{X}_\beta \ .
        \eqno\qed
    \end{equation*}
\end{lemma}

In particular, the stack~$\mathcal{X}_\alpha$
only depends on the \emph{non-degenerate quotient face} of~$\alpha$,
which we define in \cref{def:tame-formal-lattice} below.

\begin{example}[Linear quotient stacks]
    \label{eg-component-lattice-linear-quotient}
    Suppose $S = \Spec K$, where $K$ is a field.
    Consider the quotient stack $\mathcal{X} = V / G$,
    where $V$ is a finite-dimensional $K$-vector space,
    $G$~is a smooth affine algebraic group over~$K$
    with a split maximal torus~$T \subset G$,
    and~$G$ acts linearly on~$V$.

    From the description in \cref{eg-grad-quotient-stack},
    we have isomorphisms of formal lattices
    \begin{align*}
        \mathrm{CL} (V / G) 
        & \simeq \Lambda_T / W \ ,
        \\
        \mathrm{CL}_\mathbb{Q} (V / G)
        & \simeq (\Lambda_T \otimes \mathbb{Q}) / W \ ,
    \end{align*}
    where
    $\Lambda_T = \mathrm{Hom} (\mathbb{G}_\mathrm{m}, T) \simeq \mathbb{Z}^{\dim T}$
    is the cocharacter lattice of~$T$,
    and $W = \mathrm{N}_G (T) / \mathrm{Z}_G (T)$
    is the Weyl group.
    Note that we have $\mathrm{Z}_G (T) = T$ when~$G$ is connected.
\end{example}

\begin{example}[Linear moduli stacks]
    \label{eg-component-lattice-linear-moduli}
    Let~$\mathcal{X}$ be a \emph{linear moduli stack},
    that is, a stack that behaves like
    the moduli stack of objects
    in an abelian category.
    We explain this in more detail in \cref{subsec-lms} below.
    For now, we think of $\mathcal{X}$ as
    the moduli stack of representations of a quiver,
    or the moduli stack of coherent sheaves on a projective variety,
    or similar examples.

    In this case, the set $\uppi_0 (\mathcal{X})$
    carries the structure of a commutative monoid,
    with monoid structure given by the direct sum of objects.

    For each $\mathbb{Z}$-lattice~$\Lambda$, the set
    $\mathrm{CL} (\mathcal{X}) (\Lambda) \simeq
    \uppi_0 (\Grad^\Lambda (\mathcal{X}))$
    can be identified with
    the set of maps $\gamma \colon \Lambda^\vee \to \uppi_0 (\mathcal{X})$
    such that the set
    $\mathrm{supp} (\gamma) = \{ v \in \Lambda^\vee \mid \gamma (v) \neq 0 \}$
    is finite.
    Here, $\Lambda^\vee$ is understood as the grading lattice.
    For each such $\gamma$, we have an identification
    \begin{equation*}
        \mathcal{X}_\alpha \simeq
        \prod_{v \in \mathrm{supp} (\gamma)} \mathcal{X}_{\gamma (v)} \ ,
    \end{equation*}
    where $\alpha$ is the face corresponding to~$\gamma$,
    and $\mathcal{X}_{\gamma (v)} \subset \mathcal{X}$
    is the corresponding connected component.

    Similarly, for each finite-dimensional $\mathbb{Q}$-vector space~$F$,
    the set $\mathrm{CL}_\mathbb{Q} (\mathcal{X}) (F)$
    can be identified with the set of maps
    $\gamma \colon F^\vee \to \uppi_0 (\mathcal{X})$
    such that the set
    $\mathrm{supp} (\gamma) = \{ v \in F^\vee \mid \gamma (v) \neq 0 \}$
    is finite, and we have
    $\mathcal{X}_\alpha \simeq
    \prod_{v \in \mathrm{supp} (\gamma)} \mathcal{X}_{\gamma (v)}$,
    where $\alpha$ is the corresponding face.

    This implies that we have the following explicit description:
    We have isomorphisms of formal $\mathbb{Z}$-lattices
    \begin{equation*}
        \mathrm{CL} (\mathcal{X}) \simeq
        \coprod_{\gamma \in \uppi_0 (\mathcal{X})}
        \mathrm{CL} (\mathcal{X}_\gamma) \simeq
        \coprod_{\gamma \in \uppi_0 (\mathcal{X})} {}
        \biggl[
            \biggl( \,
                \coprod_{\gamma = \gamma_1 + \cdots + \gamma_n}
                \mathbb{Z}^n
            \, \biggr) \bigg/ {\sim}
        \biggr] \ ,
    \end{equation*}
    where we run through all classes $\gamma \in \uppi_0 (\mathcal{X})$
    and all decompositions of $\gamma$
    into non-zero elements
    $\gamma_1, \dotsc, \gamma_n \in \uppi_0 (\mathcal{X}) \setminus \{ 0 \}$,
    and the relation~$\sim$ is generated by the following:

    \begin{itemize}
        \item
            The copy of $\mathbb{Z}^n$
            corresponding to~$(\gamma_1, \dotsc, \gamma_n)$
            is identified with the copy of~$\mathbb{Z}^n$
            corresponding to $(\gamma_{\sigma (1)}, \dotsc, \gamma_{\sigma (n)})$
            for a permutation~$\sigma$ of~$\{ 1, \dotsc, n \}$,
            where the basis vectors of~$\mathbb{Z}^n$
            are permuted accordingly.

        \item
            If $\gamma_i = \gamma_{i, 1} + \cdots + \gamma_{i, m_i}$
            for $i = 1, \dotsc, n$, where each $\gamma_{i, j}$ is non-zero,
            then the copy of~$\mathbb{Z}^n$
            corresponding to~$(\gamma_1, \dotsc, \gamma_n)$
            is identified with the sublattice
            $\mathbb{Z}^n \subset \mathbb{Z}^{m_1 + \cdots + m_n}$
            spanned by the sum of the first~$m_1$ basis vectors,
            the sum of the next~$m_2$ basis vectors, and so on,
            in the copy of~$\mathbb{Z}^{m_1 + \cdots + m_n}$
            corresponding to $(\gamma_{1, 1}, \dotsc, \gamma_{n, m_n})$.
    \end{itemize}
    A similar description holds for the rational component lattice,
    where we replace~$\mathbb{Z}^n$ with~$\mathbb{Q}^n$.
\end{example}

\begin{para}[The component lattice of~\texorpdfstring{$\mathcal{X}_\alpha$}{X\_α}]
    \label{para-component-lattice-x-alpha}
    Let $\mathcal{X}$ be a stack as in \cref{assumption-stack-basic},
    and let $\alpha \colon F \to \mathrm{CL}_\mathbb{Q} (\mathcal{X})$
    be a face. We describe the component lattice of
    the stack~$\mathcal{X}_\alpha$ defined in \cref{para-notation-x-alpha}.

    The isomorphism $\Grad^{F'} (\Grad^F (\mathcal{X}))
    \simeq \Grad^{F \times F'} (\mathcal{X})$
    for all $F' \in \mathsf{Lat} (\mathbb{Q})$ identifies
    $\mathrm{CL}_\mathbb{Q} (\mathrm{Grad}^F (\mathcal{X}))
    \simeq \mathrm{Map}(F,\mathrm{CL}_\mathbb{Q} (\mathcal{X}))$,
    where the right-hand side is the
    \emph{mapping formal lattice} defined by
    $\mathrm{Map}(F,\mathrm{CL}_\mathbb{Q} (\mathcal{X}))(F')=\mathrm{CL}_\mathbb{Q} (\mathcal{X}) (F \times F')$.
    Thus, for any~$F'$,
    the set $\mathrm{CL}_\mathbb{Q} (\mathcal{X}_\alpha) (F')$
    corresponds to faces
    $\beta \colon F \times F' \to \mathrm{CL}_\mathbb{Q} (\mathcal{X})$
    with $\beta |_{F \times \{ 0 \}} = \alpha$.

    In particular, a face
    $\beta \colon F' \to \mathrm{CL}_\mathbb{Q} (\mathcal{X})$
    together with a map of faces $\alpha\to \beta$
    lifts canonically to a face
    $\tilde{\beta} \colon F' \to \mathrm{CL}_\mathbb{Q} (\mathcal{X}_\alpha)$,
    corresponding to the face
    $\beta \circ (f, \mathrm{id}) \colon
    F \times F' \to F' \to \mathrm{CL}_\mathbb{Q} (\mathcal{X})$. 
    The lift~$\tilde{\beta}$ maps to~$\beta$
    under the induced map
    $\mathrm{CL}_\mathbb{Q} (\mathcal{X}_\alpha) \to \mathrm{CL}_\mathbb{Q} (\mathcal{X})$.
    Taking $\beta = \alpha$ and the identity morphism $\alpha \to \alpha$
    gives a tautological face
    $\tilde{\alpha} \colon F \to \mathrm{CL}_\mathbb{Q} (\mathcal{X}_\alpha)$.

    The assignment $(\alpha\to \beta)\mapsto (\tilde{\alpha}\to \tilde{\beta})$ induces equivalences of comma categories
    \[
        (\alpha\downarrow \mathsf{Face}(\cX))\simeq (\tilde{\alpha}\downarrow \mathsf{Face}(\cX_\alpha)) \quad \text{and} \quad (\alpha\downarrow \mathsf{Face}^\mathrm{nd}(\cX))\simeq (\tilde{\alpha}\downarrow \mathsf{Face}^\mathrm{nd}(\cX_\alpha))\ .
    \]
    These follow from \cref{para:uniqueness-factorizations} and \cref{prop:characterization-tame} \cref{item:condition-tameness} below. 
    Moreover, \cref{lemma-grad-surj-rational} yields natural isomorphisms
    $(\mathcal{X}_\alpha)_{\tilde{\beta}} \simeq \mathcal{X}_\beta$
    for all morphisms $\alpha \to \beta$ in $\mathsf{Face} (\mathcal{X})$.
\end{para}

\subsection{Tameness}

\begin{para}
    We introduce a nicely behaved class of formal $\bQ$-lattices, called \emph{tame}, that include the rational component lattice $\mathrm{CL}_\bQ(\cX)$ of a stack $\cX$ (\cref{theorem-non-degenerate-quotient-functor}). 
    Tame formal lattices interact nicely with maps that are unramified in the sense of\/~\cref{para-faces} (\cref{prop:characterization-tame,prop:unramified-base-change}). 
    This is useful for our proof of the finiteness theorem (\cref{thm-finiteness}).
\end{para}

\begin{definition}[Tame formal lattice]\label{def:tame-formal-lattice}
    A formal $\bQ$-lattice $X$ is \emph{tame} if its non-degenerate faces form a face arrangement, that is, if the inclusion $\mathsf{Face}^{\mathrm{nd}}(X) \hookrightarrow \mathsf{Face}(X)$
    has a left adjoint
    \[
        (-)^\mathrm{nd} \colon
        \mathsf{Face}(X) \longrightarrow \mathsf{Face}^{\mathrm{nd}}(X) \ ,
    \]
    called the \emph{non-degenerate quotient} functor.
\end{definition}

\begin{para}[Uniqueness of factorizations]
    \label{para:uniqueness-factorizations}
    If a face $\alpha\colon F_1\to X$ of a rational formal lattice $X$ factors through a surjection $f\colon F_1\twoheadrightarrow F_2$, then it does so uniquely. 
    This fact follows from the existence of a section of $f$, and it will be used repeatedly in this section.
\end{para}

\begin{para}[Null faces]
    \label{para:null}
    We introduce some terminology that will be useful to characterize tame formal lattices.
    A face $\alpha\colon F\to X$ is said to be \emph{null} if it factors through the zero map $F\to 0$. 
    Equivalently, $\alpha$ is null if and only if $\alpha\circ 0_F=\alpha$, where $0_F\colon F\to F$ is the zero map. 
    
    If $\alpha$ is null and $\alpha\twoheadrightarrow \beta$ is a surjective map of faces, then $\beta$ is also null. 
    Indeed, calling $u\colon F\twoheadrightarrow F_\beta$ the underlying map of vector spaces, we have
    \[\beta\circ u=\alpha=\alpha\circ 0_F=\beta\circ (F\xrightarrow{0} F_\beta)=\beta\circ 0_{F_\beta}\circ u\]
    and hence $\beta=\beta\circ 0_{F_\beta}$.
\end{para}

Tame formal lattices can be characterized in several useful ways as follows.

\begin{proposition}
    \label{prop:characterization-tame}
    Let $X$ be a formal $\bQ$-lattice. 
    Then the following properties are equivalent: 
    \begin{enumerate}
        \item \label{item:tame} The formal lattice $X$ is tame.
        \item \label{item:unramified-faces} Every non-degenerate face $\alpha\colon F\to X$ is an unramified map of formal lattices, in the sense of\/~\cref{para-faces}.
        \item \label{item:condition-tameness} For every face $\alpha\colon F\to X$ and every subspace $K\subset F$, if the restriction $\alpha|_K$ is null, then $\alpha$ factors through $F\twoheadrightarrow F/K$.
        \item \label{item:maximal-null-subface}For every face $\alpha\colon F\to X$, there is a (necessarily unique) subspace $K_\alpha \subset F$ such that $\alpha$ factors through $F\twoheadrightarrow F/K_\alpha$ and $K_\alpha$ contains every null subface of $\alpha$.
    \end{enumerate}
\end{proposition}

\begin{proof}
    First, assume \cref{item:tame}. 
    Let $\alpha$ be a non-degenerate face of $X$ and let $\beta\hookrightarrow \alpha$ be a subface. 
    We want to show that $\beta$ is non-degenerate. 
    By the adjunction property, we have injective arrows $\beta\hookrightarrow \beta^{\mathrm{nd}} \hookrightarrow \alpha$. 
    If $\beta$ is degenerate, then there is a non-trivial quotient $q\colon \beta \twoheadrightarrow \beta_1$. 
    Choosing $q$ so that $\dim \beta_1$ is minimal, we may assume that the face $\beta_1$ is non-degenerate. 
    Now, $q$ must factor through $\beta \hookrightarrow \beta^{\mathrm{nd}}$, which is injective, contradicting the non-triviality of $q$. 
    This shows~\cref{item:unramified-faces}.

    Assume \cref{item:unramified-faces}. 
    Take a face $\alpha\colon F\to X$ and a subspace $K\subset F$ such that $\alpha|_K$ is null. 
    Take a quotient $f\colon \alpha\twoheadrightarrow \beta$ of faces with $\dim \beta$ minimal, so that $\beta$ is non-degenerate. 
    Then $f(K)$ is a null subface of $\beta$, so $f(K)=0$ and thus $f$ factors through $F\twoheadrightarrow F/K$. 
    Thus \cref{item:unramified-faces} implies \cref{item:condition-tameness}.

    Now assume \cref{item:condition-tameness} and let $\alpha\colon F\to X$ be a face. 
    To show \cref{item:maximal-null-subface}, it is enough to prove that, for two null subfaces $K_1,K_2\subset F$ of $\alpha$, their sum $K_1+K_2$ is also a null subface. 
    The face $\alpha$ factors through $q\colon F\twoheadrightarrow F/K_1$, giving a face $\beta\colon F/K_1\to X$. 
    The quotient $q$ induces a surjection $\alpha\vert_{K_2}\twoheadrightarrow \beta\vert_{q(K_2)}$, so $\beta\vert_{q(K_2)}$ is null. 
    We also have a surjection $\alpha\vert_{K_1+K_2}\twoheadrightarrow \beta\vert_{q(K_1+K_2)}=\beta\vert_{q(K_2)}$, so $\alpha\vert_{K_1+K_2}$ is null, as desired.

    Finally, assume \cref{item:maximal-null-subface}. 
    For every face $\alpha\colon F\to X$, let $\alpha\to \alpha^\mathrm{nd}$ be the unique factorization of $\alpha$ through $F\twoheadrightarrow F/K_\alpha$. 
    To establish \cref{item:tame} it is enough to show that every map $f\colon \alpha \to \beta$ with $\beta$ non-degenerate factors uniquely through $\alpha\to \alpha^\mathrm{nd}$. Since we have a surjection $\alpha\vert_{K_\alpha}\twoheadrightarrow \beta\vert_{f(K_\alpha)}$, $f(K_\alpha)$ is a null subface of $\beta$, so $f(K_\alpha)\subset K_\beta=0$. 
    Thus, $f$ factors through $F\twoheadrightarrow F/K_\alpha$, and by \cref{para:uniqueness-factorizations} this factorization is unique.
\end{proof}

\begin{para}
    Let $X$ be a tame formal $\bQ$-lattice. 
    From the proof of \cref{prop:characterization-tame} we see that the non-degenerate quotient $\alpha^{\mathrm{nd}}$ of a face $\alpha\colon F\to X$ is the unique factorization of $\alpha$ through $F\twoheadrightarrow F/K_\alpha$. 
    In particular, the unit map $\alpha\to \alpha^{\mathrm{nd}}$ is surjective. 
    Moreover, it is clear from \cref{item:maximal-null-subface} that, for a map $f\colon \alpha\to \beta$ of faces, we have $K_\alpha = f^{-1}(K_\beta)$. 
    In particular, if $f$ is surjective, then $\alpha^{\mathrm{nd}} \to \beta^{\mathrm{nd}}$ is an isomorphism.
\end{para}

\begin{para}
    We say that a map $f\colon X\to Y$ of formal $\bQ$-lattices \emph{reflects null faces} if, for all $\alpha\in \mathsf{Face}(X)$ such that $f\circ \alpha$ is null, the face $\alpha$ is null too. 
\end{para}

\begin{lemma}
    \label{lemma:unramified-reflects-null}
    Let $f\colon X\to Y$ be a map of formal $\bQ$-lattices. If $f$ is unramified, as in\/~\cref{para-faces}, then it reflects null faces. 
    Conversely, if $f$ reflects null faces and $X$ is tame, then $f$ is unramified.
\end{lemma}
\begin{proof}
    Suppose that $f$ is unramified, and let $\alpha$ be a face of $X$ with $\dim \alpha>0$ such that $f\circ \alpha$ is null. 
    Then $\alpha$ must be degenerate, so there is a surjection of faces $\alpha \twoheadrightarrow \beta$ with $\dim \beta < \dim \alpha$. 
    We have a surjection $f\circ \alpha \twoheadrightarrow f\circ \beta$, so $f\circ \beta$ is also null, by \cref{para:null}. 
    By induction on dimension, $\beta$ is null, and hence so is $\alpha$.

    Now suppose that $f$ reflects null faces and that $X$ is tame, and let $\alpha\colon F\to X$ be a non-degenerate face of $X$. 
    If $f\circ \alpha$ is degenerate, then there is a non-zero subspace $K\subset F$ such that $(f\circ \alpha)|_K$ is null. Therefore, $\alpha|_K$ is also null, and by tameness of $X$ and \cref{prop:characterization-tame}, $\alpha$ is degenerate, a contradiction. 
\end{proof}

\begin{proposition}\label{prop:unramified-base-change}
    Consider a cartesian square
    \[
        \begin{tikzcd}
            X' \arrow{r}{f'} \ar[d, "g'"'] \arrow[rd , phantom, "\ulcorner", near start] & Y' \arrow{d}{g} \\
            X \arrow{r}{f} & Y
        \end{tikzcd}
    \]
    of formal $\bQ$-lattices, with $X$ and $Y'$ tame. Then $X'$ is tame. 
    If, moreover, $f$ is unramified, then so is $f'$.
\end{proposition}

\begin{proof}
    Observe that a face $\alpha\colon F\to X'$ is null if and only if both $f'\circ\alpha$ and $g'\circ \alpha$ are null. 
    Also, using \cref{para:uniqueness-factorizations} we see that $\alpha$ factors through a surjection $F\twoheadrightarrow F'$ if and only if $f'\circ\alpha$ and $g'\circ \alpha$ do. 
    By \cref{prop:characterization-tame} \cref{item:condition-tameness}, we conclude that $X'$ is tame. 

    Suppose in addition that $f$ is unramified. By \cref{lemma:unramified-reflects-null}, $f'$ is unramified if it reflects null faces. 
    Let $\alpha\colon F\to X'$ be a face such that $f'\circ \alpha$ is null. 
    Then $g\circ f'\circ \alpha$ is null so, since $f$ reflects null faces, $g'\circ \alpha$ is null too. 
    It follows that $\alpha$ is null.
\end{proof}

\begin{para}[Some properties of integral component lattices]\label{para-properties-integral-component-lattice}
    Let $\mathcal X$ be a stack as in \cref{assumption-stack-basic}. 
    \cref{lemma-grad-surj} yields an important property of $\mathrm{CL}(\cX)$, namely, that for a $\bZ$-lattice epimorphism $f\colon \Lambda_2\to \Lambda_1$, the induced map $\mathrm{CL}(\cX)(\Lambda_1)\to \mathrm{CL}(\cX)(\Lambda_2)$ is injective.
    
    Now, let $\alpha\colon F\to \mathrm{CL}_\bQ(\cX)$ be a face, given by a lattice $\Lambda\subset F$ of full rank and a face $\gamma\colon \Lambda\to \mathrm{CL}(\cX)$ of the integral component lattice. 
    It is proven in \cite[Proposition~1.3.9]{halpern-leistner-instability} that there is a subgroup $H_\gamma\subset T_\Lambda$ such that, for any field $K$ and any map $x\colon \mathrm{B}T_{\Lambda,K}\to \cX$ that gives a point of $\mathcal X_\gamma$, we have $\ker(T_{\Lambda,K}\to \mathrm{Aut}(x\vert_{\Spec K}))=(H_{\gamma})_K$.

    We note that the subgroup $H_\gamma\subset T_\Lambda$ is determined solely by $\mathrm{CL}(\cX)$ and $\gamma$. 
    Indeed, let $q\colon \Lambda\to \Hom(\Gm,T_\Lambda/H_\gamma)$ be induced by the quotient $T_\Lambda\to T_\Lambda/H_\gamma$. 
    Then $q$ is the maximum in the poset of $\bZ$-lattice epimorphisms $g\colon \Lambda\to \Lambda'$ through which $\gamma$ factors.
    Now, $q$ determines $H_\gamma$ by Cartier duality.
\end{para}

\begin{theorem}\label{theorem-non-degenerate-quotient-functor}
    Let $\cX$ be a stack as in \cref{assumption-stack-basic}. 
    Then the rational component lattice $\mathrm{CL}_\bQ(\cX)$ is tame.
\end{theorem}

\begin{proof}
    With notations as in \cref{para-properties-integral-component-lattice}, we let $K_\alpha=\Hom(\Gm,H_\gamma)_\bQ\subset F$. 
    Clearly $\alpha$ factors through $F\twoheadrightarrow F/K_\alpha$. 
    Now let $F_1\subset F$ be a subspace such that $\alpha\vert_{F_1}$ is null. 
    Then $\alpha\vert_{F_1\cap \Lambda}$ factors through $F_1\cap \Lambda\to 0$, so $F_1\cap \Lambda\subset \Hom(\Gm,H_\gamma)$ and $F_1\subset K_\alpha$. 
    Therefore, $\mathrm{CL}_\bQ(\cX)$ satisfies property \cref{item:maximal-null-subface} in \cref{prop:characterization-tame}, and thus it is tame.
\end{proof}

\begin{lemma}
    \label{lemma-cl-unramified}
    Let $\mathcal{X}$ and~$\mathcal{Y}$ be stacks as in \cref{assumption-stack-basic},
    and let $f \colon \mathcal{X} \to \mathcal{Y}$ be a representable morphism.
    Then the induced map
    $f_* \colon \mathrm{CL}_\mathbb{Q} (\mathcal{X}) \to \mathrm{CL}_\mathbb{Q} (\mathcal{Y})$
    is unramified in the sense of\/~\cref{para-faces}.
\end{lemma}

\begin{proof}
    Let $\alpha$ and $\gamma$ be as in \cref{para-properties-integral-component-lattice}. 
    By representability, we have that $H_\gamma=H_{f_*\gamma}$. 
    In particular, the equality $K_\alpha=K_{f_*\alpha}$ holds by the proof of \cref{theorem-non-degenerate-quotient-functor}. 
    The result follows.
\end{proof}

\subsection{Comparison with the component fan}

\begin{para}
    The component lattice $\mathrm{CL}(\cX)$ is a replacement for the
    \emph{component fan} $\mathrm{CF}_\bullet(\cX)$ introduced by
    \textcite[Definition~3.6.1]{halpern-leistner-instability}.
    In the following, we explain the relationship between the two objects.
\end{para}

\begin{para}[The component fan]
    As in \cite[Definition~3.1.1]{halpern-leistner-instability},
    define a category $\mathfrak{Cone}$ of
    integral simplicial cones as the
    subcategory of $\mathsf{Lat}(\bZ)$ whose objects are $\bZ^n$ for $n \geq 0$,
    and morphisms are injective homomorphisms $\bZ^m \to \bZ^n$ with
    non-negative matrix coefficients, that is, mapping the standard basis in
    $\bZ^m$ to the cone spanned by the standard basis in $\bZ^n$.
    A \emph{formal fan} is a functor
    $Y_\bullet \colon \mathfrak{Cone}^{\rm{op}} \to \mathsf{Set}$,
    and we write $Y_n = Y_\bullet(\bZ^n)$.

    For a stack~$\mathcal{X}$ as in \cref{assumption-stack-basic},
    its \emph{component fan} $\mathrm{CF}_\bullet(\cX)$ is the formal fan
    defined in the notations of this paper by
    \[
        \mathrm{CF}_\bullet (\cX) (\mathbb{Z}^n) =
        \mathrm{CL} (\cX) (\mathbb{Z}^n)^\mathrm{nd}
    \]
    for objects $\mathbb{Z}^n$ in $\mathfrak{Cone}^\mathrm{op}$,
    where
    $\mathrm{CL} (\mathcal{X}) (\mathbb{Z}^n)^\mathrm{nd} \subset
    \mathrm{CL} (\mathcal{X}) (\mathbb{Z}^n)$
    is the subset of non-degenerate faces.
\end{para}

\begin{para}[Geometric realization]
    Any formal fan $Y_\bullet$ has a \emph{geometric realization},
    which is the topological space
    \begin{equation*}
        |Y_\bullet| = \colim_{(\bZ^n,\alpha) \in (\mathfrak{Cone} \downarrow Y_\bullet)} (\bR_{\geq 0})^n \ ,
    \end{equation*}
    equipped with the colimit topology of the standard topologies on the copies of $(\bR_{\geq 0})^n$ \cite[\S 3.1]{halpern-leistner-instability}. Here $(\mathfrak{Cone} \downarrow Y_\bullet)$ denotes the overcategory $\{\alpha \colon \bZ^n \to Y_\bullet\}$ with morphisms consisting of injective homomorphisms over $Y_\bullet$ with non-negative matrix coefficients. The geometric realization $|Y_\bullet|$ has a projectivized variant 
    \[
        \mathbb{P}(Y_\bullet) = (|Y_\bullet| \setminus \{\text{rank 0 points}\}) / \bR^\times_{> 0} \ ,
    \]
    where rank $0$ points are the image of $Y_0$, and the group $\bR_{> 0}^\times$ acts on each copy of $(\bR_{\geq 0})^n$ in the colimit by scaling.

    On the other hand, recall from \cref{para-geometric-realization}
    that any formal $\mathbb{Z}$-lattice~$X$ has a geometric realization~$|X_\mathbb{R}|$,
    which can be written as a similar colimit~\cref{eq-real-geometric-realization}
    indexed by the category~$\mathsf{Face}(X)$,
    which is also an overcategory $(\mathsf{Lat}(\bZ) \downarrow X)$.

    The space $\bP(\mathrm{CF}_\bullet(\cX))$ is of principal interest in \cite{halpern-leistner-instability}. A numerical invariant on $\cX$ with values in a totally ordered set $\Gamma$ is a function $\mu \colon \bP(Y_\bullet) \to \Gamma$ for some sub-fan $Y_\bullet \subset \mathrm{CF}_\bullet(\cX)$ (see \cite[Definition~4.1.1]{halpern-leistner-instability}). The main results of \cite{halpern-leistner-instability} establish conditions under which a numerical invariant defines a $\Theta$-stratification of $\cX$. The following lemma shows that $\mathrm{CL}(\cX)$ can be regarded as a different combinatorial structure on the same underlying space $|\mathrm{CF}_\bullet(\cX)|$, and is therefore useful for applying results based on \cite{halpern-leistner-instability} to $\mathrm{CL}(\cX)$.
\end{para}

\begin{lemma}
    Let $\cX$ be a stack as in \cref{assumption-stack-basic}. The overcategory $(\mathfrak{Cone} \downarrow \mathrm{CF}_\bullet(\cX))$ is a subcategory of $\mathsf{Face}(\mathrm{CL}(\cX))$, and this inclusion of indexing categories along with the standard inclusions $(\bR_{\geq 0})^n \hookrightarrow \bR^n$ determines a homeomorphism of colimits
    \[
        |\mathrm{CF}_\bullet(\cX)| \simeq |\mathrm{CL}_\mathbb{R} (\cX)| \ .
    \]
\end{lemma}
\begin{proof} From the definition of $\mathrm{CF}_\bullet(\cX)$, we can re-express $(\mathfrak{Cone} \downarrow \mathrm{CF}_\bullet(\cX)) = (\mathfrak{Cone} \downarrow \mathrm{CL}(\cX))^{\rm{nd}})$, where on the right-hand-side we regard $\mathfrak{Cone}$ as a subcategory of $\mathsf{fLat}(\bZ)$, and the superscript denotes the full subcategory of non-degenerate morphisms $\bZ^n \to \mathrm{CL}(\cX)$ in $\mathsf{fLat}(\bZ)$. \cref{theorem-non-degenerate-quotient-functor} implies that $\mathsf{Face}^{\rm nd}(\mathrm{CL}(\cX)) \subset \mathsf{Face}(\mathrm{CL}(\cX))$ is cofinal, so it suffices to show that the map 
\[
\colim_{(\bZ^n,\alpha) \in (\mathfrak{Cone}\downarrow \mathrm{CL}(\cX))^{\rm{nd}}} (\bR_{\geq 0})^n \longrightarrow \colim_{(F,\alpha) \in (\mathsf{Lat}(\bZ)\downarrow \mathrm{CL}(\cX))^{\rm{nd}}} F_\bR
\]
coming from the faithful functor $(\mathfrak{Cone}\downarrow \mathrm{CL}(\cX)) \hookrightarrow (\mathsf{Lat}(\bZ) \downarrow \mathrm{CL}(\cX))$ and the inclusion $(\bR_{\geq 0})^n \hookrightarrow \bR^n$ is a homeomorphism.

To ease notation, we will denote $X = \mathrm{CL}(\cX)$. Now consider the category
\[
    \cD = \bigl\{
        \bZ^n \overset{\beta}{\longrightarrow} F \overset{\alpha}{\longrightarrow} X
        \bigm|
        \text{$\alpha$, $\beta$  non-degenerate} 
    \bigr\} \ ,
\]
where morphisms are pairs $(f \colon \bZ^{n_1} \to \bZ^{n_2}, \ g \colon F_1 \to F_2)$ making the relevant diagram commute, and such that $f$ has non-negative matrix coefficients. We have functors
\begin{equation*}
    (\mathfrak{Cone}\downarrow X)^{\rm{nd}}
    \overset{p}{\longleftarrow}
    \cD
    \overset{q}{\longrightarrow}
    (\mathsf{Lat}(\bZ)\downarrow X)^{\rm{nd}} \ ,
\end{equation*}
where $p(\bZ^n \to F \to X)$ is the composition $(\bZ^n \to X)$, and $q(\bZ^n \to F \to X) = (F \to X)$. The functor $p$ is cofinal because for any $y = (\bZ^n \to X) \in (\mathfrak{Cone}\downarrow X)^{\rm{nd}}$, the comma category $(y \downarrow p)$ has an initial object $(x=(\bZ^n \overset{\mathrm{id}}{\to} \bZ^n \to X), \ y \simeq p(x))$. The functor $q$ is a cocartesian fibration, where the pushforward of $\beta_1 \colon \bZ^{n} \to F_1$ along a morphism $(F_1,\alpha_1) \to (F_2,\alpha_2)$ in $(\mathsf{Lat}(\bZ)\downarrow X)^{\rm{nd}}$ is the composition $\beta_2 \colon \bZ^n \to F_1 \to F_2$. The details of the arguments that $p$ is cofinal and that $q$ is a cocartesian fibration both use the fact that, for any morphism of non-degenerate faces $(F_1,\alpha_1) \to (F_2,\alpha_2)$ of $X$, the underlying homomorphism of lattices $F_1 \to F_2$ is injective.

The fibre $q^{-1}(F \to X)$ is the category $(\mathfrak{Cone} \downarrow F)^{\rm{nd}}$, so the facts that $p$ is cofinal and $q$ is a cocartesian fibration imply that
\begin{align*}
\colim_{(\mathfrak{Cone}\downarrow X)^{\rm{nd}}} (\bR_{\geq 0})^n &\simeq \colim_{\cD} (\bR_{\geq 0})^n
\simeq \colim_{(F,\alpha) \in (\mathsf{Lat}(\bZ)\downarrow X)^{\rm{nd}}} \
\biggl(
    \colim_{(\mathfrak{Cone} \downarrow F)^{\rm{nd}}} (\bR_{\geq 0})^n
\biggr) \ .
\end{align*}
Therefore, to prove the claim, it suffices to show that for any $N\geq 0$, the canonical map $\colim_{(\mathfrak{Cone} \downarrow \bZ^N)^{\rm{nd}}} (\bR_{\geq 0})^n \to \bR^N$ is a homeomorphism. The left hand side is the geometric realization of the formal fan denoted $R_\bullet(\{\bR^N\})$ in \cite[Construction 3.1.5]{halpern-leistner-instability}, so the desired homeomorphism follows from \cite[Lemma~3.1.7]{halpern-leistner-instability}. 
\end{proof}

\section{Special faces}

\label{sec-faces}

\subsection{Special faces}
\label{subsec-faces}

\begin{para}
    We introduce the notion of \emph{special faces} for an algebraic stack~$\mathcal{X}$,
    which are faces~$\alpha$ in the component lattice~$\mathrm{CL}_\mathbb{Q} (\mathcal{X})$,
    such that one cannot enlarge~$\alpha$
    without changing the stack~$\mathcal{X}_\alpha$
    defined in \cref{para-notation-x-alpha}.
    The special faces are the ones that cut out
    the wall-and-chamber structure on the component lattice
    described in \cref{para-intro-component-lattice}.

    In \cref{thm-special-face-closure},
    we show that there is a \emph{special face closure} functor
    \[
        (-)^{\mathrm{sp}} \colon \mathsf{Face}(\mathcal{X})
        \longrightarrow \mathsf{Face}^\mathrm{sp}(\mathcal{X}) \ ,
    \]
    where $\mathsf{Face}^\mathrm{sp}(\mathcal{X})$
    is the category of special faces.
    For a face~$\alpha$,
    its special face closure~$\alpha^{\mathrm{sp}}$
    is the minimal special face containing~$\alpha$,
    in the sense that the morphism $\alpha \to \alpha^{\mathrm{sp}}$
    is universal among morphisms from~$\alpha$ to special faces.
    We have a canonical isomorphism
    $\mathcal{X}_{\smash{\alpha^{\mathrm{sp}}}} \simeq \mathcal{X}_{\alpha}$.

    In particular, this implies that the special faces of~$\mathcal{X}$
    form a face arrangement on~$\mathrm{CL}_{\mathbb{Q}} (\mathcal{X})$
    in the sense of \cref{para-face-arrangement}.
\end{para}

\begin{para}[Special faces]
    \label{para-special-faces}
    Let $\mathcal{X}$ be an algebraic stack as in \cref{assumption-stack-basic}.
    Recall from \cref{para-component-lattice}
    the categories
    $\mathsf{Face}^\mathrm{nd} (\mathcal{X}) \subset
    \mathsf{Face} (\mathcal{X})$
    of faces of~$\mathcal{X}$,
    and from \cref{para-notation-x-alpha}
    the stack~$\mathcal{X}_\alpha$ associated to a face~$\alpha$.

    A \emph{special face} of~$\mathcal{X}$
    is a non-degenerate face
    $\alpha \in \mathsf{Face}^\mathrm{nd} (\mathcal{X})$,
    such that for any morphism $f \colon \alpha \to \beta$
    in $\mathsf{Face}^\mathrm{nd} (\mathcal{X})$,
    if the induced morphism
    $\mathcal{X}_\beta \to \mathcal{X}_\alpha$
    is an isomorphism, then~$f$ is an isomorphism.
    We denote by
    \begin{equation*}
        \mathsf{Face}^\mathrm{sp} (\mathcal{X}) \subset
        \mathsf{Face}^\mathrm{nd} (\mathcal{X})
    \end{equation*}
    the full subcategory of special faces.
\end{para}

\begin{para}[The central rank]
    \label{para-central-rank}
    Let $\mathcal{X}$ be a stack as in \cref{assumption-stack-basic},
    and assume that~$\mathcal{X}$ is connected.

    The \emph{central rank} of~$\mathcal{X}$,
    denoted by $\crk (\mathcal{X})$,
    is the largest integer~$n \geq 0$
    such that $\mathrm{CL}_\mathbb{Q} (\mathcal{X})$
    admits a non-degenerate face~$\alpha$ of dimension~$n$,
    such that the morphism
    $\mathrm{tot}_\alpha \colon \mathcal{X}_\alpha \to \mathcal{X}$
    is an isomorphism,
    as in \textcite[\S3.4]{ibanez-nunez-filtrations}.

    For example, if $\mathcal{X} = \mathrm{B} G$,
    where $G$ is an affine algebraic group over a field~$K$ with a split maximal torus,
    then the central rank of~$\mathcal{X}$
    is the rank of the centre of~$G$.
\end{para}

    We now give a useful characterization of special faces that
    will be proven in \cref{proof-thm-special-face-crk} below.

\begin{theorem}
    \label{thm-special-face-crk}
    Let $\mathcal{X}$ be a stack as in \cref{assumption-stack-basic}.
    Then for any non-degenerate face
    $\alpha \in \mathsf{Face}^\mathrm{nd} (\mathcal{X})$,
    we have the inequality
    \begin{equation*}
        \crk (\mathcal{X}_\alpha) \geq \dim \alpha \ ,
    \end{equation*}
    with equality if and only if~$\alpha$ is special.
\end{theorem}

The following is the main result in this section, whose proof will be given in \cref{subsec-special-face-closure}.

\begin{theorem}
    \label{thm-special-face-closure}
    Let $\mathcal{X}$ be a stack as in \cref{assumption-stack-basic}.
    Then there is a functor
    \[
        (-)^{\mathrm{sp}} \colon \mathsf{Face}(\mathcal{X})
        \longrightarrow \mathsf{Face}^\mathrm{sp}(\mathcal{X}) \ ,
    \]
    called the special face closure,
    which is left adjoint to the inclusion
    $\mathsf{Face}^\mathrm{sp} (\mathcal{X}) \hookrightarrow
    \mathsf{Face} (\mathcal{X})$.
    In particular, $\mathsf{Face}^\mathrm{sp} (\mathcal{X})$
    is a face arrangement in~$\mathrm{CL}_{\mathbb{Q}}(\mathcal{X})$
    in the sense of\/~\cref{para-face-arrangement}.

    Moreover, for each face $\alpha \in \mathsf{Face} (\mathcal{X})$,
    the adjunction unit
    $\alpha \to \alpha^{\mathrm{sp}}$
    induces an isomorphism
    \begin{equation*}
        \mathcal{X}_{\smash{\alpha^{\mathrm{sp}}}}
        \longsimto \mathcal{X}_{\alpha} \ .
    \end{equation*}

    When~$\mathcal{X}$ is connected,
    the category~$\mathsf{Face}^\mathrm{sp} (\mathcal{X})$
    has an initial object~$\alpha_{\mathrm{ce}}$
    called the maximal central face of\/~$\mathcal{X}$,
    with~$\dim \alpha_{\mathrm{ce}} = \mathrm{crk}(\mathcal{X})$.
    It is the maximal non-degenerate face such that the morphism
    $\mathrm{tot}_{\alpha_{\mathrm{ce}}} \colon
    \mathcal{X}_{\alpha_{\mathrm{ce}}} \to \mathcal{X}$
    is an isomorphism.
\end{theorem}

\begin{example}[Linear quotient stacks]
    \label{eg-special-face-linear-quotient}
    Let $\mathcal{X} = V / G$ be a linear quotient stack,
    as in \cref{eg-component-lattice-linear-quotient},
    and assume that~$G$ is connected.
    Recall that $\mathrm{CL}_\mathbb{Q} (V / G) \simeq (\Lambda_T \otimes \mathbb{Q}) / W$.
    We describe the special faces of~$\mathcal{X}$.

    Let $\Lambda^T = \Lambda_T^\vee$
    be the weight lattice of~$T$.
    Consider the hyperplane arrangement~$\Phi_{V / G}$
    in $\Lambda_T \otimes \mathbb{Q}$
    given by the following two types of hyperplanes:

    \begin{itemize}
        \item
            Hyperplanes~$P_w$
            dual to non-zero weights~$w \in \Lambda^T$
            present in the $T$-representation~$V$.

        \item
            Hyperplanes~$Q_r$
            dual to the roots~$r \in \Lambda^T$ of~$G$.
    \end{itemize}
    The flats of~$\Phi_{V / G}$ (see \cref{para-hyperplane-arrangements})
    define a face arrangement in
    $(\Lambda_T \otimes \mathbb{Q}) / W$.

    In fact, this face arrangement is precisely
    $\mathsf{Face}^\mathrm{sp} (\mathcal{X})$,
    since for a face $\alpha \colon F \to \mathrm{CL}_\mathbb{Q} (V / G)$,
    we have $\mathcal{X}_\alpha \simeq V^\alpha / L_\alpha$,
    with notations analogous to that in \cref{eg-grad-quotient-stack},
    and the flats of~$\Phi_{V / G}$
    precisely give faces~$\alpha$ that are maximal among faces
    with given~$V^\alpha$ and~$L_\alpha$.
\end{example}

\begin{example}[Linear moduli stacks]
    \label{eg-special-face-lms}
    Let~$\mathcal{X}$ be a linear moduli stack,
    as in \cref{eg-component-lattice-linear-moduli}.
    Recall that the faces $F \to \mathrm{CL}_\mathbb{Q} (\mathcal{X})$
    correspond to maps $\alpha \colon F^\vee \to \uppi_0 (\mathcal{X})$
    such that the set
    $\mathrm{supp} (\alpha) = \{ v \in F^\vee \mid \alpha (v) \neq 0 \}$
    is finite, and we have
    $\mathcal{X}_\alpha \simeq \prod_{v \in F^\vee} \mathcal{X}_{\alpha (v)}$.
    Such a face is non-degenerate if and only if
    $\mathrm{supp} (\alpha)$ spans~$F^\vee$.

    For classes
    $\gamma_1, \ldots, \gamma_n \in \uppi_0 (\mathcal{X})$,
    consider the face
    \begin{equation*}
        \alpha (\gamma_1, \dotsc, \gamma_n) \in \mathsf{Face} (\mathcal{X})
    \end{equation*}
    corresponding to the map
    $\mathbb{Q}^n \to \uppi_0 (\mathcal{X})$
    sending the $i$-th standard basis vector to~$\gamma_i$,
    and all other elements to~$0$,
    as explained in \cref{eg-component-lattice-linear-moduli}.

    We argue that special faces of~$\mathcal{X}$
    are necessarily of this form,
    with $\gamma_1, \ldots, \gamma_n$ non-zero.
    Indeed, for a face
    corresponding to a map $\alpha \colon F^\vee \to \uppi_0 (\mathcal{X})$,
    such a face is non-degenerate if and only if
    $\mathrm{supp} (\alpha)$ spans~$F^\vee$,
    and if it is special,
    then it is maximal among faces with
    the same multiset of non-zero values of~$\alpha$,
    so that $\mathrm{supp} (\alpha)$ has to be a basis of~$F^\vee$.

    This will be stated precisely
    in \cref{para-lms-special-faces} below.
\end{example}

\subsection{The special face closure}
\label{subsec-special-face-closure}

\begin{para}
    This section is devoted to the proof of \cref{thm-special-face-closure}.
\end{para}

\begin{para}[Central faces]
    Let~$\mathcal{X}$ be a stack as in \cref{assumption-stack-basic},
    and suppose that~$\mathcal{X}$ is connected.

    A \emph{central face} of~$\mathcal{X}$
    is a face $\alpha \colon F \to \mathrm{CL}_\mathbb{Q} (\mathcal{X})$
    such that the morphism
    $\mathrm{tot}_\alpha \colon \mathcal{X}_{\alpha} \to \mathcal{X}$
    is an isomorphism.
    Denote by
    $\mathsf{Face}^\mathrm{ce}(\mathcal{X}) \subset \mathsf{Face}(\mathcal{X})$
    the full subcategory of central faces.

    Similarly, a \emph{central integral face} of~$\mathcal{X}$
    is a face $\alpha \colon \Lambda \to \mathrm{CL} (\mathcal{X})$
    such that the morphism
    $\mathrm{tot}_\alpha \colon \mathcal{X}_{\alpha} \to \mathcal{X}$
    is an isomorphism.
    Equivalently, the rationalization of~$\alpha$
    is a central face.
\end{para}

\begin{para}[Central actions]
    \label{para-central-actions}
    For a central integral face
    $\alpha \colon \Lambda \to \mathrm{CL} (\mathcal{X})$,
    the isomorphism $\mathcal{X}_{\alpha} \simeq \mathcal{X}$ induces a
    $\mathrm{B} T_\Lambda$-action
    $\mathrm{B} T_\Lambda \times \mathcal{X} \to \mathcal{X}$,
    which in turn induces an action
    $\Lambda \times \mathrm{CL}(\mathcal{X}) \to \mathrm{CL}(\mathcal{X})$.

    Taking rationalizations, for a central face
    $(F, \alpha) \in \mathsf{Face}^\mathrm{ce}(\mathcal{X})$,
    there is an action
    \[
        m_{\alpha} \colon F \times \mathrm{CL}_{\mathbb{Q}}(\mathcal{X})
        \longrightarrow \mathrm{CL}_{\mathbb{Q}}(\mathcal{X}) \ ,
    \]
    which we call the \emph{central action} of~$\alpha$.

    For a central face
    $(F, \alpha) \in \mathsf{Face}^\mathrm{ce}(\mathcal{X})$
    and a face $(F', \alpha') \in \mathsf{Face}(\mathcal{X})$,
    define a face $(F \times F', \alpha \boxplus \alpha') \in \mathsf{Face}(\mathcal{X})$
    by the composition
    \[
        \alpha \boxplus \alpha' \colon F \times F'
        \overset{\id_F \times \alpha'}{\longrightarrow}
        F \times \mathrm{CL}_{\mathbb{Q}}(\mathcal{X})
        \overset{m_{\alpha}}{\longrightarrow}
        \mathrm{CL}_{\mathbb{Q}}(\mathcal{X}) \ .
    \]
    Applying $\Grad^{F'}$ to the isomorphism
    $\mathcal{X}_{\alpha} \simeq \mathcal{X}$,
    we obtain an isomorphism
    $\mathcal{X}_{\alpha \, \boxplus \, \alpha'} \simeq \mathcal{X}_{\alpha'}$.
\end{para}

\begin{para}[The maximal central face]
    Let~$\mathcal{X}$ be as in \cref{assumption-stack-basic},
    and suppose that~$\mathcal{X}$ is connected.
    There is a non-degenerate central face
    $(F_{\mathrm{ce}}, \alpha_{\mathrm{ce}}) \in \mathsf{Face}^\mathrm{ce}(\mathcal{X})$,
    such that $\dim \alpha_{\mathrm{ce}} = \mathrm{crk} (\mathcal{X})$. By the definition of central rank, $\alpha_{\mathrm{ce}}$ is special.

    We show below that~$\alpha_\mathrm{ce}$
    is the \emph{maximal central face} of~$\mathcal{X}$,
    in that it is the final object of~$\mathsf{Face}^\mathrm{ce}(\mathcal{X})$.
    In particular, $\alpha_{\mathrm{ce}}$ does not depend on choices and
    $\mathrm{Aut}(\alpha_{\mathrm{ce}})$ is trivial.  
\end{para}

\begin{para}[Unique map from a central face to a special face]
    \label{para-central-face-unique-map}
    Let $\alpha \colon F \to \mathrm{CL}_{\mathbb{Q}}(\mathcal{X})$ be a central face and let $\beta \colon F' \to \mathrm{CL}_{\mathbb{Q}}(\mathcal{X})$ be a special face. 
    Then we have
    $(\alpha \boxplus \beta)^{\mathrm{nd}} \simeq \beta$,
    which follows from the isomorphism
    $\mathcal{X}_{\alpha \, \boxplus \, \beta} \simeq \mathcal{X}_{\beta}$
    from \cref{para-central-actions}.
    In particular, there is a canonical morphism 
    $ \iota_{\alpha,\beta} \colon \alpha
        \to \beta $.

    We claim that this is the unique morphism
    $\alpha \to \beta$.
    Indeed, for any map $j \colon \alpha \to \beta$,
    consider the induced morphism
    $\alpha \boxplus j \colon
    \alpha \boxplus \alpha \to \alpha \boxplus \beta$.
    Precomposing with the morphisms
    $(0, \mathrm{id}), (\mathrm{id}, 0) \colon
    \alpha \to \alpha \boxplus \alpha$,
    we see that the composition of the non-degenerate quotient $\alpha\to \alpha^{\mathrm{nd}}$ with the morphism
    $(\alpha \boxplus j)^{\mathrm{nd}} \colon \alpha^{\mathrm{nd}} \to \beta$
    is equal to both~$j$ and $\iota_{\alpha,\beta}$.

    Taking $\beta = \alpha_{\mathrm{ce}}$, we obtain the claim above that $\alpha_{\mathrm{ce}}$ is the final object of $\mathsf{Face}^\mathrm{ce}(\mathcal{X})$. Also, since $\alpha_{\mathrm{ce}}$ is both special and central, it follows that the comma category $(\alpha\downarrow \mathsf{Face}^\mathrm{sp}(\mathcal{X}))$ has an initial object, given by the canonical map $\alpha\to \alpha_{\mathrm{ce}}$.

    Taking $\alpha=\alpha_{\mathrm{ce}}$, we obtain that $\alpha_{\mathrm{ce}}$ is the initial object of $\mathsf{Face}^\mathrm{sp}(\mathcal{X})$, as claimed in \cref{thm-special-face-closure}. 
\end{para}

\begin{para}[Special faces of \texorpdfstring{$\mathcal{X}_\alpha$}{X\_α}]
    \label{para-special-faces-x-alpha}
    Let $\alpha \colon F \to \mathrm{CL}_{\mathbb{Q}}(\mathcal{X})$ be a face. 
    The equivalence 
    $(\alpha \downarrow \mathsf{Face}(\mathcal{X}))\simeq (\tilde{\alpha}\downarrow \mathsf{Face}(\mathcal{X}_{\alpha}))$
    in \cref{para-component-lattice-x-alpha} restricts to an equivalence $(\alpha \downarrow \mathsf{Face}^\mathrm{sp}(\mathcal{X}))\simeq (\tilde{\alpha}\downarrow \mathsf{Face}^\mathrm{sp}(\mathcal{X}_{\alpha}))$.
    Since $\tilde{\alpha}$ is a central face of $\mathcal{X}_{\alpha}$, we have $\mathsf{Face}^\mathrm{sp}(\mathcal{X}_\alpha)
        \simeq (\tilde{\alpha} \downarrow \mathsf{Face}^\mathrm{sp}(\mathcal{X}_{\alpha}))$ by \cref{para-central-face-unique-map}.
    Thus, we obtain an equivalence
    \[
        \mathsf{Face}^\mathrm{sp}(\mathcal{X}_\alpha)
        \simeq (\alpha\downarrow \mathsf{Face}^\mathrm{sp}(\mathcal{X})) \ .
    \]
\end{para}

\begin{para}[Proof of \texorpdfstring{\cref{thm-special-face-closure}}{Theorem \ref{thm-special-face-closure}}]
    \label{proof-thm-special-face-closure}
    It remains to establish the adjunction, for which we need to show that, for every face $\alpha \in \mathsf{Face}(\mathcal{X})$, the comma category $(\alpha \downarrow \mathsf{Face}^\mathrm{sp}(\mathcal{X}))$ has an initial object. 
    This follows from \cref{para-special-faces-x-alpha} and \cref{para-central-face-unique-map}.
    \qed
\end{para}

\begin{para}[Proof of \texorpdfstring{\cref{thm-special-face-crk}}{Theorem \ref{thm-special-face-crk}}]
    \label{proof-thm-special-face-crk}
    By \cref{para-component-lattice-x-alpha}, after replacing $\mathcal{X}$ with $\mathcal{X}_\alpha$, we may assume that $\alpha$ is a central face of $\mathcal{X}$. There is a unique map $\iota\colon \alpha\hookrightarrow \alpha_{\mathrm{ce}}$, so $\dim \alpha \leq \dim \alpha_{\mathrm{ce}}=\mathrm{crk}(\mathcal X)$, with equality if and only if $\iota$ is an isomorphism, that is, if an only if $\alpha$ is special.
\end{para}

\section{Special cones}

\label{sec-filt-cones}

\subsection{Cone filtrations}
\label{subsec-filt-cones}

\begin{para}[Idea]
    \label{para-idea-filt-cones}
    We introduce a generalization of the stack
    $\Filt^n (\mathcal{X})$ of $\mathbb{Z}^n$-filtered points,
    based on the following idea.

    As we have seen in \cref{eg-special-face-linear-quotient},
    when~$\mathcal{X} = V / G$ is a linear quotient stack,
    the component lattice $\mathrm{CL}_\mathbb{Q} (\mathcal{X})$
    admits a hyperplane arrangement~$\Phi_{V / G}$
    given by the non-zero weights in~$V$ and the roots of~$G$.
    For a point~$\lambda \in |\mathrm{CL} (\mathcal{X})|$,
    as in \cref{eg-grad-quotient-stack},
    the components
    $\mathcal{X}_\lambda \subset \mathrm{Grad} (\mathcal{X})$
    and
    $\mathcal{X}_\lambda^+ \subset \mathrm{Filt} (\mathcal{X})$
    corresponding to~$\lambda$ are given by
    $\mathcal{X}_\lambda \simeq V^\lambda / L_\lambda$
    and
    $\mathcal{X}_\lambda^+ \simeq V^{\lambda, +} / P_\lambda$,
    and they do not change when we vary~$\lambda$
    inside a chamber in a flat in~$\Phi_{V / G}$.
    This is because~$\mathcal{X}_\lambda$ only depends on
    the set of weights whose pairing with~$\lambda$ is zero,
    and $\mathcal{X}_\lambda^+$ only depends on
    the set of weights whose pairing with~$\lambda$ is non-negative.

    On the other hand, for any cone
    $\sigma \in \mathsf{Cone} (\mathcal{X})$,
    we are going to construct a stack~$\mathcal{X}_\sigma^+$,
    which is, roughly speaking,
    an analogue of the stack~$\mathcal{X}_\lambda^+$ above,
    but with respect to weights with non-negative pairing with all points in~$\sigma$.
    When~$\sigma$ is of the standard form~$(\mathbb{Q}_{\geq 0})^n \subset \mathbb{Q}^n$,
    this construction coincides with a component in $\Filt^n (\mathcal{X})$.

    The construction of the stack~$\mathcal{X}_\sigma^+$
    is similar to the construction of~$\Filt^n (\mathcal{X})$,
    but we replace the stack~$\Theta^n = \mathbb{A}^n / \mathbb{G}_\mathrm{m}^n$
    by a stack of the form~$R / \mathbb{G}_\mathrm{m}^n$,
    where~$R$ is an affine toric variety corresponding to the cone~$\sigma$.
    We then define~$\mathcal{X}_\sigma^+$
    as a component in the mapping stack~$\Map (R / \mathbb{G}_\mathrm{m}^n, \mathcal{X})$.

    As we will see in \cref{lemma-filt-sigma-image-rational},
    this construction does not produce new stacks,
    since each stack~$\mathcal{X}_\sigma^+$
    can always be identified with a component in~$\Filt^n (\mathcal{X})$ for some~$n$.
    More precisely, if~$\sigma$ is generated by the vectors~$v_1, \dotsc, v_n$,
    then~$\mathcal{X}_\sigma^+$ is isomorphic to the component of~$\Filt^n_\mathbb{Q} (\mathcal{X})$
    corresponding to the possibly degenerate face
    $\mathbb{Q}^n \to \mathrm{CL}_\mathbb{Q} (\mathcal{X})$
    sending the standard basis vectors to the vectors~$v_1, \dotsc, v_n$.

    We remark that the idea to consider general cones
    also appears in \textcite[Section 3]{odaka2024stability}.
\end{para}

\begin{para}[The stack~\texorpdfstring{$\Theta_\Sigma$}{Θ\_Σ}]
    \label{para-theta-sigma}
    Let~$\Sigma$ be an integral cone in the sense of \cref{para-integral-cones}.
    We now construct a stack~$\Theta_\Sigma$,
    which simultaneously generalizes the stacks $\mathrm{B} \mathbb{G}_\mathrm{m}^n$
    and~$\Theta^n \simeq \mathbb{A}^n / \mathbb{G}_\mathrm{m}^n$
    used to define the stacks of graded and filtered points
    in \cref{para-grad-filt}.
    These two cases correspond to the integral cones~$\mathbb{Z}^n$ and~$\mathbb{N}^n$, respectively.

    For an integral cone~$\Sigma$,
    define its \emph{dual cone} as the integral cone
    $\Sigma^\vee = \mathrm{Hom} (\Sigma, \mathbb{N})$,
    which is the set of monoid homomorphisms $\Sigma \to \mathbb{N}$.
    For example, if~$\Sigma = \mathbb{Z}^m \oplus \mathbb{N}^n$,
    then~$\Sigma^\vee \simeq \mathbb{N}^n$.

    Let~$\Lambda_\Sigma$ be the groupification of~$\Sigma$,
    and let $\Lambda_\Sigma^\vee = \mathrm{Hom} (\Lambda_\Sigma, \mathbb{Z})$.

    Define the quotient stack
    \begin{equation*}
        \Theta_\Sigma = R_\Sigma / T_\Sigma \ ,
    \end{equation*}
    where $R_\Sigma = \Spec \mathbb{Z} [\Sigma^\vee]$,
    and $T_\Sigma = \Spec \mathbb{Z} [\Lambda_\Sigma^\vee]$
    is the split torus with cocharacter lattice~$\Lambda_\Sigma$.
    The $T_\Sigma$-action on~$R_\Sigma$ is given by
    the opposite of the natural
    $\Lambda_\Sigma^\vee$-grading on~$\mathbb{Z} [\Sigma^\vee]$,
    also given by the algebra homomorphism
    $\Z[\Sigma^\vee] \to \Z[\Lambda_\Sigma^\vee] \otimes \Z[\Sigma^{\vee}]$
    defined by
    $s \mapsto s \otimes s$ for $s \in \Sigma^\vee$.
    For example, we have
    $\Theta_{\mathbb{Z}^m \oplus \mathbb{N}^n} \simeq
    \mathrm{B} \mathbb{G}_\mathrm{m}^m \times \Theta^n$
    for integers $m, n \geq 0$.

    This construction is covariant in~$\Sigma$.
    There are two special points~$0, 1 \in \Theta_\Sigma (\mathbb{Z})$,
    given by the homomorphisms $\mathbb{Z} [\Sigma^\vee] \to \mathbb{Z}$
    sending all non-zero elements of~$\Sigma^\vee$ to~$0$ and~$1$, respectively.
    The first homomorphism is well-defined because
    no non-zero element in $\Sigma^\vee$ has an opposite.
\end{para}

\begin{para}[The stack of \texorpdfstring{$\Sigma$}{Σ}-filtered points]
    \label{para-filt-sigma}
    Let~$\Sigma$ be a monoid as above,
    and let~$\mathcal{X}$ be a stack as in \cref{assumption-stack-basic}.
    Define the \emph{stack of\/ $\Sigma$-filtered points}
    of~$\mathcal{X}$ as the mapping stack
    \begin{equation*}
        \Filt^\Sigma (\mathcal{X}) = \Map (\Theta_\Sigma, \mathcal{X}) \ .
    \end{equation*}
    For example, we have
    $\Filt^{\mathbb{Z}^n} (\mathcal{X}) \simeq \Grad^n (\mathcal{X})$ and
    $\Filt^{\mathbb{N}^n} (\mathcal{X}) \simeq \Filt^n (\mathcal{X})$.

    Note the inconsistency with the previous terminology of
    \emph{$\mathbb{Z}^n$-filtered points},
    which now correspond to \emph{$\mathbb{N}^n$-filtered points}.

    The stack~$\Filt^\Sigma (\mathcal{X})$
    is an algebraic stack over~$S$ satisfying conditions in
    \cref{assumption-stack-basic}.
    Indeed, since the projection $\Theta_\Sigma \to \Spec \mathbb{Z}$
    is a good moduli space,
    by \textcite[Lemma~2.2.2]{ibanez-nunez-filtrations},
    the stack~$\Filt^\Sigma (\mathcal{X})$
    is identified with the mapping stack
    $\Map_S (\Theta_\Sigma \times S, \mathcal{X})$ relative to the base~$S$,
    and is naturally defined over~$S$.
    It then follows from \textcite[Theorem~10.17]{alper-hall-rydh-etale-local}
    that $\Filt^\Sigma (\mathcal{X})$ is algebraic and locally of finite type over~$S$,
    where \cite[Setup~10.12]{alper-hall-rydh-etale-local} is satisfied
    since $\Theta_\Sigma \times S$ satisfies the condition~(N) there.
    Finally, it follows from \cref{thm-filt-sigma-image} below
    that the conditions in \cref{assumption-stack-basic}
    are satisfied.

    The assignment
    $(\Sigma, \mathcal{X}) \mapsto \Filt^\Sigma (\mathcal{X})$
    is contravariant in~$\Sigma$ and covariant in~$\mathcal{X}$.
    The morphisms
    \begin{equation*}
        \begin{tikzcd}
            \mathrm{B} T_{\Sigma}
            \ar[shift left=0.5ex, r, "0"]
            &
            \Theta_\Sigma
            \ar[shift left=0.5ex, l, "\mathrm{pr}"]
            &
            \Spec \mathbb{Z}
            \ar[shift left=0.5ex, l, "1"]
            \ar[shift right=0.5ex, l, "0"']
            \ar[ll, bend right, start anchor=north west, end anchor=north east, looseness=.8]
        \end{tikzcd}
    \end{equation*}
    induce morphisms
    \begin{equation*}
        \begin{tikzcd}
            \Grad^{\Lambda_\Sigma} (\mathcal{X})
            \ar[rr, bend left, start anchor=north east, end anchor=north west, looseness=.6, "\smash{\mathrm{tot}}"]
            \ar[r, shift right=0.5ex, "\mathrm{sf}"']
            &
            \Filt^\Sigma (\mathcal{X})
            \ar[l, shift right=0.5ex, "\mathrm{gr}"']
            \ar[r, shift left=0.5ex, "\mathrm{ev}_0"]
            \ar[r, shift right=0.5ex, "\mathrm{ev}_1"']
            &
            \mathcal{X} \vphantom{^{\Lambda_\Sigma}}
        \end{tikzcd}
    \end{equation*}
    on the mapping stacks,
    where $\Grad^{\Lambda_\Sigma} (\mathcal{X})$
    is defined in \cref{para-grad-lambda}.
\end{para}

\begin{theorem}
    \label{thm-filt-sigma-image}
    Let~$\mathcal{X}$ be a stack as in \cref{assumption-stack-basic}.
    Let $\pi \colon \Sigma' \to \Sigma$ be a map of integral cones,
    such that the induced map $\Sigma'_{\smash{\Q}} \to \Sigma_\Q$ is surjective,
    with notations as in \cref{para-integral-cones}.

    Then the induced morphism
    \begin{equation*}
        \pi^* \colon \Filt^\Sigma (\mathcal{X}) \longrightarrow \Filt^{\Sigma'} (\mathcal{X})
    \end{equation*}
    is an open and closed immersion.
\end{theorem}

This theorem is an extension of \cref{lemma-grad-surj} to filtered points.
It also implies that for any integral cone~$\Sigma$,
the stack $\Filt^\Sigma (\mathcal{X})$
is isomorphic to an open and closed substack
of~$\Filt^n (\mathcal{X})$ for some~$n$,
since we may apply the theorem to a surjection $\mathbb{N}^n \to \Sigma$.
In particular,
$\Filt^\Sigma (\mathcal{X})$ satisfies the conditions in
\cref{assumption-stack-basic},
as we claimed in \cref{para-filt-sigma}.

\begin{proof}
    First we claim that, under the assumptions of the theorem,
    the morphism $\Theta_{\Sigma'}\to \Theta_\Sigma$
    induced by~$\pi$
    is a \emph{good moduli space morphism}
    in the sense of \textcite[Remark~4.4]{alper-2013-good-moduli}.
    Denote $\Lambda=\Lambda_\Sigma$ and $\Lambda'=\Lambda_{\Sigma'}$, and let $H$ be the kernel of the surjective homomorphism $T_\Lambda\to T_{\Lambda'}$. We have a cartesian square
    \begin{equation*}
        \begin{tikzcd}
            R_{\Sigma'}/H
            \ar[r] \ar[d]
            \ar[dr, phantom, "\ulcorner", very near start]
            & R_\Sigma \ar[d] \\
            \Theta_{\Sigma'} \ar[r]
            & \Theta_\Sigma \rlap{ ,}
        \end{tikzcd}
    \end{equation*}
    and the right vertical arrow is a smooth cover.
    Thus, it is enough to show that
    $R_{\Sigma'}/H\to R_\Sigma$
    is a good moduli space,
    or, in more concrete terms,
    that we have $\Z[(\Sigma')^\vee]^H=\Z[\Sigma^\vee]$.
    Note that the equality
    $\Sigma^\vee=(\Sigma')^\vee\cap \Lambda^\vee$
    inside $(\Lambda')^\vee$
    holds due to the surjectivity of $\Sigma'_{\Q}\to \Sigma_\Q$.
    Since $\Z[(\Sigma')^\vee]^H=\Z[(\Sigma')^\vee\cap \Lambda^\vee]$,
    the claim follows.

    Let $R$ be an affine scheme and consider the map
    $\Theta_{\Sigma'}\times R\to \Theta_\Sigma\times R$ induced by $\pi$.
    By the universal property of good moduli space morphisms
    in \textcite[Theorem~7.22]{alper-hall-rydh-etale-local},
    the groupoid $\mathcal X(\Theta_\Sigma\times R)$
    is the full subgroupoid of $\mathcal X(\Theta_{\Sigma'}\times R)$
    consisting of those morphisms $\phi\colon \Theta_{\Sigma'}\times R\to \mathcal X$
    satisfying that, for every field-valued point
    $x \colon {\Spec K} \to R$,
    the homomorphism
    $T_{\Lambda',K}=\mathrm{Aut}(0,x)\to \mathrm{Aut}(\phi(0,x))$
    induced by~$\phi$ factors through $T_{\Lambda',K}\to T_{\Lambda,K}$.
    By the proof of
    \textcite[Proposition~1.3.9]{halpern-leistner-instability},
    this is an open and closed condition on $\vert R\vert$.
\end{proof}

\begin{para}[Rational filtrations]
    \label{para-rational-filt}
    We now define a rationalized version
    of the stack of $\Sigma$-filtered points,
    for rational cones instead of integral cones.

    Let~$\mathcal{X}$ be a stack as in \cref{assumption-stack-basic},
    and let $C$ be a rational cone as in \cref{para-cones}.
    Similarly to \cref{para-rational-graded-points},
    we define the \emph{stack of\/ $C$-filtered points} of~$\mathcal{X}$ by
    \begin{equation*}
        \Filt^C (\mathcal{X}) =
        \operatornamewithlimits{colim}_{\Sigma \subset C}
        \Filt^\Sigma (\mathcal{X}) \ ,
    \end{equation*}
    where the colimit is indexed by integral cones
    $\Sigma \subset C$ with $C = \mathbb{Q}_{\geq 0} \cdot \Sigma$.

    We denote
    $\Filt^n_\mathbb{Q} (\mathcal{X}) = \Filt^{(\mathbb{Q}_{\geq 0})^n} (\mathcal{X})$
    and
    $\Filt_\mathbb{Q} (\mathcal{X}) = \Filt^{\mathbb{Q}_{\geq 0}} (\mathcal{X})$.

    By \cref{thm-filt-sigma-image},
    all morphisms in the colimit diagram are open and closed immersions,
    so $\Filt^C (\mathcal{X})$ is an algebraic stack satisfying conditions in
    \cref{assumption-stack-basic},
    and each component of $\Filt^C (\mathcal{X})$
    is isomorphic to a component of $\Filt^\Sigma (\mathcal{X})$
    for an integral cone $\Sigma \subset C$
    with $C = \mathbb{Q}_{\geq 0} \cdot \Sigma$.
    In particular,
    each component of $\Filt^n_\mathbb{Q} (\mathcal{X})$
    is isomorphic to a component of $\Filt^n (\mathcal{X})$.

    The assignment
    $(C, \mathcal{X}) \mapsto \Filt^C (\mathcal{X})$
    is contravariant in~$C$ and covariant in~$\mathcal{X}$.
    We have the induced morphisms
    \begin{equation*}
        \begin{tikzcd}
            \Grad^{F_C} (\mathcal{X})
            \ar[rr, bend left, start anchor=north east, end anchor=north west, looseness=.6, "\smash{\mathrm{tot}}"]
            \ar[r, shift right=0.5ex, "\mathrm{sf}"']
            &
            \Filt^C (\mathcal{X})
            \ar[l, shift right=0.5ex, "\mathrm{gr}"']
            \ar[r, shift left=0.5ex, "\mathrm{ev}_0"]
            \ar[r, shift right=0.5ex, "\mathrm{ev}_1"']
            &
            \mathcal{X} \ , \vphantom{^{F_C}}
        \end{tikzcd}
    \end{equation*}
    defined by the colimits of the morphisms in \cref{para-filt-sigma},
    where~$F_C$ is the groupification of~$C$,
    and $\Grad^{\smash{F_C}} (\mathcal{X})$
    is defined in \cref{para-rational-graded-points}.

    By \cref{thm-filt-sigma-image}, if $C' \to C$ is
    a surjective morphism of rational cones, then the induced morphism
    $\Filt^C (\mathcal{X}) \to \Filt^{C'} (\mathcal{X})$
    is an open and closed immersion.
    In particular, for any rational cone~$C$, we can choose a
    surjection $(\Q_{\geq 0})^n \to C$ for some $n$, so
    $\Filt^C (\mathcal{X})$ can be regarded as an open and closed substack
    of~$\Filt^n_\mathbb{Q} (\mathcal{X})$ for some~$n$.
\end{para}

\begin{para}[The stack~\texorpdfstring{$\mathcal{X}_\sigma^+$}{X\_σ\string^+}]
    \label{para-x-sigma-plus}
    Let~$(C, \sigma) \in \mathsf{Cone} (\mathcal{X})$ be a cone,
    and let $(F, \alpha) = \mathrm{span} (C, \sigma) \in \mathsf{Face} (\mathcal{X})$.
    Similarly to \cref{para-notation-x-alpha},
    we define
    \begin{equation*}
        \mathcal{X}_\sigma^+ \subset \Filt^C (\mathcal{X})
    \end{equation*}
    as the preimage of $\mathcal{X}_\alpha \subset \Grad^F (\mathcal{X})$
    under the morphism
    $\mathrm{gr} \colon \Filt^C (\mathcal{X}) \to \Grad^{F} (\mathcal{X})$.
    The assignment $\sigma \mapsto \mathcal{X}_\sigma^+$
    is contravariant in~$\sigma$.

    We will see in \cref{lemma-theta-retract} below that~$\mathcal{X}_\sigma^+$
    is a connected component of~$\Filt^C (\mathcal{X})$.
    There are induced morphisms
    \begin{equation*}
        \begin{tikzcd}
            \mathcal{X}_\alpha
            \ar[rr, bend left, start anchor=north east, end anchor=north west, looseness=.8, "{\smash[t]{\mathrm{tot}_\alpha}}"]
            \ar[r, shift right=0.5ex, "\mathrm{sf}_\sigma"']
            &
            \mathcal{X}_\sigma^+
            \ar[l, shift right=0.5ex, "\mathrm{gr}_\sigma"']
            \ar[r, shift left=0.5ex, "\mathrm{ev}_{0, \sigma}"]
            \ar[r, shift right=0.5ex, "\mathrm{ev}_{1, \sigma}"']
            &
            \mathcal{X} \ .
        \end{tikzcd}
    \end{equation*}

    Note that the~`$+$' in the notation~$\mathcal{X}_\sigma^+$ is redundant,
    since if~$(C, \sigma) = (F, \alpha)$ is a full face,
    then~$\mathcal{X}_\sigma^+ \simeq \mathcal{X}_\alpha$,
    so we could have written~$\mathcal{X}_\sigma$ without any ambiguity.
    However, we keep the plus sign for more clarity.

    For an element $\lambda \in |\mathrm{CL}_\mathbb{Q} (\mathcal{X})|$,
    we also write
    \begin{equation*}
        \mathcal{X}_\lambda^+ \subset \Filt_\mathbb{Q} (\mathcal{X})
    \end{equation*}
    for the connected component corresponding to~$\lambda$.
    This can be seen as the special case of~$\mathcal{X}_\sigma^+$
    for the cone $\sigma = \mathbb{Q}_{\geq 0} \cdot \lambda$.

    Similarly to \cref{para-notation-x-alpha},
    for a morphism~$\mathcal{Y} \to \mathcal{X}$ of algebraic stacks,
    we also write $\mathcal{Y}_\sigma^+ \subset \Filt^C (\mathcal{Y})$
    for the preimage of~$\mathcal{X}_\sigma^+$ under the induced morphism
    $\Filt^C (\mathcal{Y}) \to \Filt^C (\mathcal{X})$.

    We have the following direct consequence of
    \cref{thm-filt-sigma-image}.
\end{para}

\begin{lemma}
    \label{lemma-filt-sigma-image-rational}
    Let~$\mathcal{X}$ be a stack as in \cref{assumption-stack-basic}.
    Let $\pi \colon (C', \sigma') \to (C, \sigma)$
    be a morphism in $\mathsf{Cone} (\mathcal{X})$
    which is surjective as a map $C' \to C$.
    Then $\pi$ induces an isomorphism
    \begin{equation*}
        \mathcal{X}_\sigma^+ \longsimto \mathcal{X}_{\smash{\sigma'}}^+ \ .
        \eqno\qed
    \end{equation*}
\end{lemma}

This can be seen as a generalization and an analogue of
\cref{lemma-grad-surj-rational} for filtered points.

In particular, this implies that the stack~$\mathcal{X}_\sigma^+$
only depends on the \emph{underlying non-degenerate cone} of~$\sigma$,
defined as the image of~$\sigma$
in the non-degenerate quotient of the face
$\mathrm{span} (\sigma)$,
as in \cref{theorem-non-degenerate-quotient-functor}.

\begin{example}[Quotient stacks]
    \label{eg-filt-cones-quotient-stack}
    We continue the discussion in
    \cref{eg-grad-quotient-stack},
    and describe the stacks~$\mathcal{X}_\alpha$ and~$\mathcal{X}_\sigma^+$
    when~$\mathcal{X}$ is a quotient stack.

    Let~$\mathcal{X} = U / G$ be a quotient stack
    as in \cref{eg-grad-quotient-stack},
    under one of the two situations there.
    Assume for convenience that~$S$ is connected.
    Recall that we have
    $\mathrm{CL}_\mathbb{Q} (\mathrm{B} G) \simeq (\Lambda_T \otimes \mathbb{Q}) / W$,
    where~$T \subset G$ is a split maximal torus,
    $\Lambda_T$ is its cocharacter lattice,
    and~$W$ is the Weyl group.

    For a face $\alpha \in \mathsf{Face} (\mathrm{B} G)$
    and a cone $\sigma \in \mathsf{Cone} (\mathrm{B} G)$,
    we use the notations~$\mathcal{X}_\alpha$ and~$\mathcal{X}_\sigma^+$
    as defined in \cref{para-notation-x-alpha,para-x-sigma-plus}
    for the projection $\mathcal{X} \to \mathrm{B} G$.
    These stacks are given by
    \begin{align*}
        \mathcal{X}_\alpha
        & \simeq U^{\alpha} / L_\alpha \ ,
        \\
        \mathcal{X}_\sigma^+
        & \simeq U^{\sigma, +} / P_\sigma \ ,
    \end{align*}
    where $L_\alpha$ and $P_\sigma$ are
    the Levi and parabolic subgroups,
    and $U^{\alpha}$ and $U^{\sigma, +}$ are
    the fixed and attractor loci,
    generalizing the corresponding notions in \cref{eg-grad-quotient-stack}.
    The groups~$L_\alpha$ and~$P_\sigma$
    can be defined as the stabilizer groups in
    $(\mathrm{B} G)_\alpha$ and~$(\mathrm{B} G)_\sigma^+$
    at their unique $S$-points,
    and the spaces~$U^{\alpha}$ and~$U^{\sigma, +}$
    as fibres of the representable morphisms
    $\mathcal{X}_\alpha \to (\mathrm{B} G)_\alpha$
    and $\mathcal{X}_\sigma^+ \to (\mathrm{B} G)_\sigma^+$
    at these $S$-points.

    Consequently, for a face $\tilde{\alpha} \in \mathsf{Face} (\mathcal{X})$,
    writing~$\alpha$ for its image in $\mathsf{Face} (\mathrm{B} G)$,
    the stack~$\mathcal{X}_{\tilde{\alpha}}$
    is isomorphic to a connected component of $U^{\alpha} / L_\alpha$;
    a similar statement holds for cones.
\end{example}

\begin{para}[Cone filtrations in cone filtrations]
    \label{para-cone-filt-cone-filt}
    Let~$\mathcal{X}$ be a stack as in \cref{assumption-stack-basic},
    and let $C$ be a rational cone.
    We describe a construction analogous to \cref{para-component-lattice-x-alpha}
    for cone filtrations.

    By \cref{lemma-theta-retract} below, the set $\mathrm{Hom}(C,\mathrm{CL}_\bQ(\cX))$ of cones of $\cX$ of the form $(C,\sigma)$ is in bijection with $\pi_0(\mathrm{Filt}^C(\cX))$. Given such a cone $(C,\sigma)\in \mathsf{Cone}(\cX)$ and another rational cone $D$, the isomorphism
    $\Filt^{D} (\Filt^C (\mathcal{X}))
    \simeq \Filt^{C \times D} (\mathcal{X})$
    gives an identification of $\mathrm{Hom}(D,\mathrm{CL}_\bQ(\cX_\sigma^+))$ with the set of cones $f\colon D\times C \to \mathrm{CL}_\bQ(\cX)$ with $f|_{\{0\}\times C}=\sigma$. In particular, we have a natural identification $\mathrm{CL}_\bQ(\cX_\sigma^+)\simeq \mathrm{CL}_\bQ(\cX_\alpha)$, where $\alpha=\mathrm{span}(\sigma)$.

    The equivalence $(\alpha\downarrow \mathsf{Face}(\cX))\simeq (\tilde{\alpha}\downarrow \mathsf{Face}(\cX_\alpha))$ from \cref{para-component-lattice-x-alpha}
    then yields a further equivalence of comma categories
    \begin{align*}
        (\sigma\downarrow \mathsf{Cone}(\cX))
        & \longsimto
        (\tilde{\sigma}\downarrow \mathsf{Cone}(\cX_\sigma^+)) \ ,
        \\
        (\sigma \to \tau)
        & \longmapsto
        (\tilde{\sigma} \to \tilde{\tau}) \ ,
    \end{align*}
    and, by \cref{lemma-filt-sigma-image-rational}, natural isomorphisms
    $(\mathcal{X}_\sigma^+)_{\tilde{\tau}}^+ \simeq \mathcal{X}_\tau^+$
    for all morphisms $\sigma \to \tau$ in $\mathsf{Cone} (\mathcal{X})$.

    In particular, if $\sigma = \alpha$ is a full face,
    then $\mathcal{X}_\sigma^+ \simeq \mathcal{X}_\alpha$,
    and the above can be rewritten as
    \[
    (\alpha\downarrow \mathsf{Cone}(\cX))\simeq (\tilde{\alpha}\downarrow \mathsf{Cone}(\cX_\alpha))\quad\text{and}\quad (\mathcal{X}_\alpha)_{\tilde{\tau}}^+ \simeq \mathcal{X}_\tau^+\ ,
    \]
    for morphisms $\alpha\to \tau$ in $\mathsf{Cone}(\cX)$.
\end{para}

\begin{para}[\texorpdfstring{$\Theta$}{Θ}-action retracts]
    The next property we would like to prove is that the morphism
    $\mathrm{gr} \colon \Filt^\Sigma (\mathcal{X}) \to
    \Grad^{\Lambda_\Sigma} (\mathcal{X})$
    is a \emph{$\Theta$-action retract},
    where~$\Theta = \mathbb{A}^1 / \mathbb{G}_\mathrm{m}$,
    which is a stronger notion than
    an \emph{$\mathbb{A}^1$-deformation retract}
    or a \emph{$\Theta$-deformation retract}
    in the sense of \textcite[above Lemma~1.3.5]{halpern-leistner-instability}.
    This is proved when $\Sigma = \mathbb{N}$ in
    \textcite[Propositions~1.2.4 and~1.4.1]{halpern-leistner-derived}.
    In particular, it will follow that the morphism~$\mathrm{gr}$
    induces a bijection on the sets of connected components.

    Let $\cY$ and $\cZ$ be algebraic stacks, not necessarily satisfying the basic assumptions~\cref{assumption-stack-basic}, and let $\pi \colon \mathcal{Y} \to \mathcal{Z}$ be a morphism.
    We define the following notions:

    \begin{itemize}
        \item
            A structure of an \emph{$\mathbb{A}^1$-action retract} on~$\pi$
            consists of an action
            $r \colon \mathbb{A}^1 \times \mathcal{Y} \to \mathcal{Y}$ of the multiplicative monoid $\mathbb{A}^1$ and
            a map $s \colon \mathcal{Z} \to \mathcal{Y}$,
            called the \emph{zero section}, such that there exist equivalences $\pi\circ s\simeq \mathrm{id}_{\cZ}$ and $r (0, -) \simeq s \circ \pi$.

        \item
            A structure of a \emph{$\Theta$-action retract} on~$\pi$
            consists of a monoid action
            $r \colon \Theta \times \mathcal{Y} \to \mathcal{Y}$ and a map $s \colon \mathcal{Z} \to \mathcal{Y}$,
            called the \emph{zero section},
            such that there exist equivalences $\pi\circ s\simeq \mathrm{id}_{\cZ}$ and $r (0, -) \simeq s \circ \pi$.
    \end{itemize}
    A structure of an $\Theta$-action retract
    naturally gives a structure of an $\mathbb{A}^1$-action retract
    and a structure of a $\Theta$-deformation retract;
    a structure of an $\mathbb{A}^1$-action retract
    gives a structure of an $\mathbb{A}^1$-deformation retract.

    When~$\pi$ has a structure of a $\Theta$-action retract,
    by \cite[Lemmas~1.3.5 and~1.3.7]{halpern-leistner-instability},
    it induces a bijection on the sets of connected components,
    and any open set in~$\mathcal{Y}$ containing
    the image of the zero section is the whole of~$\mathcal{Y}$.
\end{para}

\begin{para}[Action retracts from actions]
An algebraic stack $\mathcal{Y}$ endowed with an $\mathbb{A}^1$-action or a $\Theta$-action admits canonically a structure of an $\mathbb{A}^1$-action retract or a $\Theta$-action retract, respectively, as follows. 
The monoid homomorphisms $\{0,1\}\to \mathbb{A}^1$ or $\{0,1\}\to \Theta$ induce an action of $\{0,1\}$ on $\mathcal{Y}$, which amounts to an idempotent $e\colon \mathcal Y\to \mathcal Y$ with the relevant coherence data.
Since the category of algebraic stacks has finite limits \cite[Lemma~4.35]{antieau-gepner-brauer-groups} and is a 2-category, the idempotent $e$ is effective \cite[Corollary~4.4.5.14]{lurie-2009-htt}, meaning that it comes from a (uniquely determined) retraction $\pi\colon \cY\to \cZ$ with section $s\colon \cZ\to \cY$ and the relevant coherence data. 
This gives a structure of an $\mathbb{A}^1$-action retract or a $\Theta$-action retract on the morphism $\pi\colon \cY\to \cZ$. 

Conversely, any structure of an $\mathbb{A}^1$- or $\Theta$-action retract on a morphism $\pi\colon \cY\to \cZ$ is determined by the $\bA^1$-action or $\Theta$-action $r$ on $\cY$, since the splitting of the idempotent $r(0,-)$ is unique in the homotopy category of the category of algebraic stacks.
\end{para}

\begin{lemma}
    \label{lemma-theta-retract}
    Let~$\mathcal{X}$ be a stack as in \cref{assumption-stack-basic},
    let~$K$ denote an integral cone $\Sigma$ or a rational cone $C$, and write $V=\Lambda_\Sigma$ or $F_C$,
    respectively. 

    Then the associated graded map
    $\mathrm{gr} \colon \Filt^K(\mathcal{X}) \to \Grad^{V} (\mathcal{X})$
    admits a non-canonical structure of a $\Theta$-action retract,
    whose zero section is given by the split filtration map
    $\mathrm{sf} \colon \Grad^{V} (\mathcal{X}) \to \Filt^K(\mathcal{X})$.

    In particular, $\mathrm{gr}$ induces a bijection
    $\uppi_0 (\Filt^K(\mathcal{X})) \simto \uppi_0 (\Grad^{V} (\mathcal{X}))$,
    and any open set in a component of\/~$\Filt^K(\mathcal{X})$
    containing the image of\/~$\mathrm{sf}$ in that component
    is the whole component.
\end{lemma}

\begin{proof}
    Let~$\lambda \in K$ be an element
    that does not lie on the boundary of~$K$,
    that is, we have $x (\lambda) > 0$
    for any $x \in K^\vee \setminus \{ 0 \}$.
    The morphism $\bN\to K\colon 1\mapsto \lambda$ gives a monoid homomorphism $\Theta\to \Theta_K$ which, by the hypothesis on $\lambda$, sends $0$ to $0$. In the case of a rational cone $C$, the monoid stack $\Theta_C$ is defined in \cref{para-rational-mapping} below, where it is shown that $\Filt^C (\mathcal{X})=\mathrm{Map}(\Theta_C,\cX)$. 
    Thus, in both cases, the natural $\Theta_K$-action on $\Filt^K (\mathcal{X})$ restricts to a $\Theta$-action $r\colon \Theta\times \Filt^K (\mathcal{X})\to \Filt^K (\mathcal{X})$ with $r(0,-)\simeq \mathrm{sf} \circ \mathrm{gr}$, as desired.
\end{proof}

\begin{para}[Rational filtrations as mapping stacks]
    \label{para-rational-mapping}
    Finally, we explain how to view
    the stacks of rational graded and filtered points
    as mapping stacks from non-algebraic sources.

    Let $\cX$ and $C$ be as in
    \cref{para-rational-filt}, let
    $F_C^\vee=\Hom(F_C,\bQ)$ be the dual vector space, and let
    $C^\vee=\Hom(C,\bQ_{\geq 0})$ be the dual cone. As in
    \cref{para-theta-sigma}, we define
    \[
        R_C=\Spec \bZ[C^\vee] \ , \qquad
        T_C=\Spec \bZ[F_C^\vee] \ , \qquad
        \Theta_C=R_C/T_C \ .
    \]
    In this case, the group scheme $T_C$ is flat,
    but not finitely presented over $\Spec \bZ$.
    The stack $\Theta_C$ is thus fpqc-algebraic,
    but not algebraic in the usual sense.
    Still, we have natural identifications
    \[
        \Filt^C(\cX)
        \simeq \Map_S(\Theta_C\times S,\cX)
        \simeq \Map(\Theta_C,\cX) \ .
    \]
    To prove this, we may assume that $S$ is an affine scheme.
    The general case then follows by descent. 
    For $\Sigma$ an integral cone, let
    $(R_{\Sigma,i})_{i\in \Delta_{\leq 2}}$ be the truncated \v{C}ech nerve
    for the presentation $R_\Sigma\to \Theta_\Sigma$,
    where $\Delta_{\leq 2}$ is the full subcategory of the simplex category
    $\Delta$ spanned by $[0],[1],[2]$. Explicitly,
    $R_{\Sigma,0}=R_\Sigma$, $R_{\Sigma,1}=T_\Sigma\times R_\Sigma$, and $R_{\Sigma,2}=T_\Sigma\times T_\Sigma\times R_\Sigma$.
    We have $\Theta_\Sigma \simeq \colim_{i\in \Delta_{\leq 2}} R_{\Sigma,i}$. For a
    rational cone $C$, we define $(R_{C,i})_{i\in \Delta_{\leq 2}}$ in the same
    way, so that $\Theta_C \simeq \colim_{i\in \Delta_{\leq 2}} R_{C,i}$. Given such a
    cone $C$, let $\cD$ denote the opposite of the poset category of integral
    cones $\Sigma\subset C$ with $C=\bQ_{\geq 0} \cdot \Sigma$, which is filtered. The
    assignment $(\Sigma,i)\mapsto R_{\Sigma,i}$ defines a functor out of the
    product category $\cD\times \Delta_{\leq 2}^\mathrm{op}$. Note that
    $\lim_{\Sigma\in \cD} R_{\Sigma,i} \simeq R_{C,i}$ for $i=0,1,2$. The result
    follows from the identifications
    \begin{align*}
        \Filt^C (\mathcal{X})
        & \simeq \colim_{\Sigma\in\cD} \Map(\Theta_\Sigma,\cX)
        \simeq \colim_{\Sigma\in\cD} \Map_S(\Theta_\Sigma\times S,\cX)
        \simeq \colim_{\Sigma\in\cD} \Map_S\Bigl( \colim_{i\in \Delta_{\leq 2}} (R_{\Sigma,i}\times S),\cX\Bigr) \\
        & \simeq \colim_{\Sigma\in\cD} \lim_{i\in \Delta_{\leq 2}} \Map_S(R_{\Sigma,i}\times S,\cX)
        \simeq \lim_{i\in \Delta_{\leq 2}} \colim_{\Sigma\in\cD} \Map_S(R_{\Sigma,i}\times S,\cX) \\
        & \simeq \lim_{i\in \Delta_{\leq 2}} \Map_S\Bigl(\lim_{\Sigma\in\cD} (R_{\Sigma,i}\times S),\cX\Bigr)
        \simeq \lim_{i\in \Delta_{\leq 2}} \Map_S(R_{C,i}\times S,\cX) \\
        & \simeq \Map_S\Bigl(\colim_{i\in \Delta_{\leq 2}} (R_{C,i}\times S),\cX\Bigr)
        \simeq \Map_S(\Theta_C\times S,\cX)
        \simeq \Map(\Theta_C,\cX) \ .
    \end{align*}
    The second equality is an application of
    \cite[Lemma~2.2.2]{ibanez-nunez-filtrations}. The fifth equality follows by
    commutativity of finite limits with filtered colimits. The sixth one holds
    by \cite[Proposition~4.18]{laumon-champs}
    because $\cX\to S$ is locally finitely presented.
    The ninth follows because $\Delta_{\leq 2}^{\mathrm{op}}$ is a sifted
    category. The last equality follows by the argument in the proof of
    \cite[Lemma~2.2.2]{ibanez-nunez-filtrations} noting that,
    for every affine scheme $T$, the map $\Theta_C\times T\to T$ is universal
    for maps from $\Theta_C\times T$ to affine schemes.
    The other steps are straightforward.
\end{para}

\subsection{Deformation theory}
\label{subsec-filt-deformation-theory}

\begin{para}
    \label{para-deformation-theory-intro}
    In this section,
    we study the deformation theory of the stacks of cone filtrations
    $\Filt^\Sigma (\mathcal{X})$.
    This will be a crucial tool in studying the structure of these stacks,
    such as proving that certain connected components
    in stacks of cone filtrations are isomorphic to each other,
    even outside of the situation of \cref{thm-filt-sigma-image},
    which we do in \cref{subsec-special-cones} below.
    In particular, the constancy theorem,
    \cref{thm-constancy},
    will be one of the results of this type.
\end{para}

\begin{para}[Derived algebraic geometry]
    \label{para-derived-algebraic-geometry}
    The deformation theory of mapping stacks, such as $\Filt^\Sigma (\mathcal{X})$,
    is perhaps more straightforward to understand
    in the context of \emph{derived algebraic geometry},
    where the cotangent complex of a mapping stack is easily described,
    and behaves better than in classical algebraic geometry.

    We list here a minimum set of concepts and properties that we will need.
    See \textcite{lurie-2009-htt,lurie-2004-dag,toen-vezzosi-2008-hag-ii}
    for relevant background.

    \begin{itemize}
        \item
            A \emph{higher stack} is a stack of $\infty$-groupoids
            over the category $\mathsf{Aff} = \mathsf{CRing}^\mathrm{op}$
            of affine schemes,
            equipped with the étale topology.
            Higher stacks form an $\infty$-category~$\mathsf{hSt}$.

        \item
            A \emph{derived stack} is a stack of $\infty$-groupoids
            over the $\infty$-category $\mathsf{dAff} = \mathsf{sCRing}^\mathrm{op}$,
            the opposite category of simplicial commutative rings,
            equipped with the étale topology.
            Derived stacks form an $\infty$-category~$\mathsf{dSt}$.
            A \emph{derived $1$-stack} is a derived stack
            whose values on classical affine schemes are $1$-groupoids.

        \item
            There is a fully faithful embedding
            $\mathsf{hSt} \hookrightarrow \mathsf{dSt}$,
            with a right adjoint $(-)_\mathrm{cl} \colon \mathsf{dSt} \to \mathsf{hSt}$,
            called the \emph{classical truncation}.
            We naturally regard higher stacks as derived stacks.

        \item
            There is the notion of an \emph{algebraic} derived stack,
            called a \emph{geometric $D^-$-stack}
            in \textcite[\S 2.2]{toen-vezzosi-2008-hag-ii}.

        \item
            A derived stack~$\mathcal{X}$ that is algebraic and
            \emph{locally almost of finite presentation}
            over a base~$\mathcal{S}$ admits a \emph{cotangent complex}
            $\mathbb{L}_{\mathcal{X} / \mathcal{S}}$,
            which is an \emph{almost perfect complex} on~$\mathcal{X}$,
            that is, a complex that is locally quasi-isomorphic to
            a bounded-above complex of vector bundles.
            This is sometimes abbreviated as $\mathbb{L}_\mathcal{X}$ when the base is clear.

        \item
            When $\cX$ and $\mathcal S$ are algebraic derived 1-stacks, there is a natural map
            $\mathbb{L}_{\mathcal{X} / \mathcal{S}} |_{\mathcal{X}_\mathrm{cl}}
            \to \mathbb{L}_{\mathcal{X}_\mathrm{cl} / \mathcal{S}_\mathrm{cl}}$
            whose induced map on cohomology sheaves
            $\mathrm{H}^i (\mathbb{L}_{\mathcal{X} / \mathcal{S}}) |_{\mathcal{X}_\mathrm{cl}}
            \to \mathrm{H}^i (\mathbb{L}_{\mathcal{X}_\mathrm{cl} / \mathcal{S}_\mathrm{cl}})$
            is an isomorphism for $i = 0, 1$, and a surjection for $i = -1$,
            by \textcite[Lemma~1.2.4]{halpern-leistner-instability} and \textcite[Proposition~25.3.5.1]{lurie-sag}.

        \item
            An algebraic stack~$\mathcal{X}$ in the classical sense,
            locally of finite presentation over a base algebraic stack~$\mathcal{S}$,
            is algebraic and locally almost of finite presentation in the derived sense,
            and has an almost perfect cotangent complex
            $\mathbb{L}_{\mathcal{X} / \mathcal{S}}$.

        \item
            Denote mapping stacks in~$\mathsf{dSt}$ by $\dMap (-, -)$,
            and call them \emph{derived mapping stacks},
            to distinguish them from classical mapping stacks.
            For a classical $1$-stack~$\mathcal{Y}$ and a derived $1$-stack~$\mathcal{X}$,
            we have $\dMap (\mathcal{Y}, \mathcal{X})_\mathrm{cl}
            \simeq \Map (\mathcal{Y}, \mathcal{X}_\mathrm{cl})$,
            as in \textcite[Lemma~1.2.1]{halpern-leistner-instability}.
            Stacks of graded and filtered points have derived versions,
            denoted by
            \begin{equation*}
                \dGrad (\mathcal{X}) \ , \qquad
                \dFilt (\mathcal{X}) \ , \qquad
                \dFilt^\Sigma (\mathcal{X}) \ , \quad
            \end{equation*}
            etc.
            Their classical truncations coincide with the classical versions.
            These are algebraic derived stacks if $\cX$ is algebraic and locally almost of finite presentation over a base algebraic space $S$, and if $\cX_\mathrm{cl}$ satisfies the basic assumptions \cref{assumption-stack-basic}.
            This follows from the algebraicity of the classical versions
            and \textcite[Theorem~2.2.6.11]{toen-vezzosi-2008-hag-ii},
            where the existence of a cotangent complex follows from
            \textcite[Proposition~5.1.10]{halpern-leistner-preygel-2023-mapping-stacks}
            and \cref{lemma-sheaves-on-theta} below.
    \end{itemize}
\end{para}

\begin{lemma}
    \label{lemma-sheaves-on-theta}
    Let~$\mathcal{X}$ be an algebraic derived stack
    as in \cref{para-derived-algebraic-geometry},
    over a base algebraic space~$S$ as in \cref{assumption-stack-basic},
    and let~$\Sigma$ be an integral cone. Consider the morphisms
    \begin{equation*}
        \mathrm{B} T_\Sigma \times \mathcal{X}
        \overset{p}{\longleftarrow}
        \Theta_\Sigma \times \mathcal{X}
        \overset{\pi}{\longrightarrow}
        \mathcal{X} \ ,
    \end{equation*}
    where~$T_\Sigma$ and~$\Theta_\Sigma$ are as in \cref{para-theta-sigma}.
    Identify quasi-coherent sheaves on
    $\mathrm{B} T_\Sigma \times \mathcal{X}$
    with $\Lambda_\Sigma^\vee$-graded quasi-coherent sheaves on~$\mathcal{X}$.

    Then the pullback functor
    $\pi^* \colon \mathsf{QCoh} (\mathcal{X}) \to
    \mathsf{QCoh} (\Theta_\Sigma \times \mathcal{X})$
    admits a left adjoint~$\pi_+$,
    and the functors $\pi_* \circ p^*$
    and $\pi_+ \circ p^*$
    send a $\Lambda_\Sigma^\vee$-graded quasi-coherent sheaf\/~$E$
    on~$\mathcal{X}$ to the parts
    \begin{equation*}
        E_{\Sigma, +} =
        \bigoplus_{\lambda \smash{\geq_\Sigma^{}} 0} E_\lambda \ ,
        \qquad
        E_{\Sigma, -} =
        \bigoplus_{\lambda \smash{\leq_\Sigma^{}} 0} E_\lambda \ ,
    \end{equation*}
    respectively,
    where we define the partial order~$\leq_\Sigma^{}$
    on~$\Lambda_\Sigma^\vee$ by
    $\lambda \leq_\Sigma^{} \lambda'$ if and only if
    $\lambda' - \lambda \in \Sigma^\vee$.
\end{lemma}

See also \textcite[Lemma~1.3.1]{halpern-leistner-derived}
for the special case when $\Sigma = \mathbb{N}$.

\begin{proof}
    It is enough to prove this when
    $\mathcal{X} = \Spec A$ is an affine derived scheme,
    so that quasi-coherent sheaves on $\Theta_\Sigma \times \Spec A$
    can be identified with $\Lambda_\Sigma^\vee$-graded
    $A [\Sigma^\vee]$-modules,
    where $A [\Sigma^\vee]$ is equipped with the $\Lambda_\Sigma^\vee$-grading
    opposite to the natural one.
    The functor~$\pi^*$ sends an $A$-module~$F$
    to the $A [\Sigma^\vee]$-module~$F [\Sigma^\vee]$,
    which gives copies of~$F$ in some of the degrees
    with respect to the $\Lambda_\Sigma^\vee$-grading, and~$0$ in the others.
    It thus preserves all colimits and limits,
    and has a left adjoint~$\pi_+$
    by the adjoint functor theorem.

    Let $E = (E_\lambda)_{\lambda \in \Lambda_\Sigma^\vee}$
    be a $\Lambda_\Sigma^\vee$-graded $A$-module.
    Then $p^* (E)$ is the $\Lambda_\Sigma^\vee$-graded
    $A [\Sigma^\vee]$-module
    $\bigoplus_{\lambda \in \Lambda_\Sigma^\vee} E_\lambda [\Sigma^\vee]$,
    and taking~$\pi_*$ picks out its $0$-th graded piece,
    which can be identified with~$E_{\Sigma, +}$.
    On the other hand, for any $A$-module~$F$, we have
    $\mathrm{Hom}_{A [\Sigma^\vee]} (p^* (E), F [\Sigma^\vee]) \simeq
    \bigoplus_{\lambda \leq_\Sigma^{} 0} \mathrm{Hom}_A (E_\lambda, F)$,
    which shows that~$\pi_+$ picks out~$E_{\Sigma, -}$.
\end{proof}

\begin{theorem}
    \label{thm-tangent-filt-sigma}
    Let $\mathcal{X}$ be a stack as in \cref{assumption-stack-basic},
    and let~$\Sigma$ be an integral cone.
    Then we have
    \begin{equation}
        \mathbb{L}_{\smash{\dFilt^\Sigma (\mathcal{X})}}
        |_{\smash{\dGrad^{\Lambda_\Sigma} (\mathcal{X})}}
        \simeq \bigl(
            \mathbb{L}_\mathcal{X}
            |_{\smash{\dGrad^{\Lambda_\Sigma} (\mathcal{X})}}
        \bigr)_{\Sigma, -} \ ,
    \end{equation}
    where we use the derived versions of the morphisms
    $\mathrm{sf}$ and\/~$\mathrm{tot}$ in \cref{para-filt-sigma},
    and~$(-)_{\Sigma, -}$ is as in \cref{lemma-sheaves-on-theta}
    for the canonical $\Lambda_\Sigma^\vee$-grading on
    $\mathbb{L}_\mathcal{X}
    |_{\smash{\dGrad^{\Lambda_\Sigma} (\mathcal{X})}}$.

    In particular, if\/ $\Sigma = \Lambda$ is a $\mathbb{Z}$-lattice, then
    \begin{equation}
        \label{eq-tangent-grad-lambda}
        \mathbb{L}_{\smash{\dGrad^\Lambda (\mathcal{X})}}
        \simeq \bigl(
            \mathbb{L}_\mathcal{X}
            |_{\smash{\dGrad^{\Lambda} (\mathcal{X})}}
        \bigr)_0 \ ,
    \end{equation}
    where $(-)_0$ denotes the $0$-th graded piece
    for the canonical $\Lambda^\vee$-grading.
\end{theorem}

This generalizes
\textcite[Lemma~1.2.3]{halpern-leistner-instability}
from the case of $\Sigma = \mathbb{N}$ to general cones,
and the special case \cref{eq-tangent-grad-lambda} was also proved in
\cite[Lemma~1.2.2]{halpern-leistner-instability}.

\begin{proof}
    It is enough to prove the first statement.
    Let $\mathrm{ev} \colon \Theta_\Sigma \times \dFilt^\Sigma (\mathcal{X}) \to \mathcal{X}$
    be the evaluation morphism.
    By \textcite[Proposition~5.1.10]{halpern-leistner-preygel-2023-mapping-stacks},
    we have
    \begin{equation*}
        \mathbb{L}_{\smash{\dFilt^\Sigma (\mathcal{X})}} \simeq
        \pi_+ \circ \mathrm{ev}^* (\mathbb{L}_\mathcal{X}) \ ,
    \end{equation*}
    where $\pi \colon
    \Theta_\Sigma \times \dFilt^\Sigma (\mathcal{X})
    \to \dFilt^\Sigma (\mathcal{X})$
    is the projection. Consider the following commutative diagram:
    \begin{equation*}
        \begin{tikzcd}
            \dGrad^\Lambda(\cX)
            \ar[d, "\mathrm{sf}"']
            &
            \Theta_\Sigma \times \dGrad^\Lambda(\cX)
            \ar[r, "\mathrm{\mathrm{pr}\times \mathrm{id}}"]
            \ar[l, "\pi'"']
            \ar[d, "\mathrm{id} \times \mathrm{sf}"]
            &
            \mathrm{B}T_\Lambda \times \dGrad^\Lambda(\cX)
            \ar[d, "\mathrm{ev}'"]
            \\
            \dFilt^\Sigma (\mathcal{X})
            &
            \Theta_\Sigma \times \dFilt^\Sigma (\mathcal{X})
            \ar[r, "\mathrm{ev}"]
            \ar[l, "\pi"']
            &
            \cX \rlap{\ .}
        \end{tikzcd}
    \end{equation*}
    Here, $\pi'$ is the projection, and $\mathrm{ev}'$ is the evaluation morphism for $\dGrad^\Lambda(\cX)$. We have
    \begin{align*}
        \mathrm{sf}^*(\mathbb{L}_{\dFilt^\Sigma(\cX)})
        & \simeq \mathrm{sf}^*\circ \pi_+\circ \mathrm{ev}^*(\mathbb{L}_\cX)
        \\
        & \simeq \pi'_+\circ (\mathrm{id}\times \mathrm{sf})^*\circ \mathrm{ev}^*(\mathbb{L}_\cX)
        \\
        & \simeq \pi'_+\circ (\mathrm{pr}\times \mathrm{id})^*\circ (\mathrm{ev}')^*(\mathbb{L}_\cX) \ ,
    \end{align*}
    where the second equality holds by \cite[Lemma~5.1.8]{halpern-leistner-preygel-2023-mapping-stacks} and cartesianity of the left square.
    The theorem then follows from \cref{lemma-sheaves-on-theta}.
\end{proof}

\begin{para}[Remark]
    It is perhaps more intuitive to think of
    \cref{thm-tangent-filt-sigma}
    as a statement about tangent complexes, that is, we have
    \begin{align*}
        \mathbb{T}_{\smash{\dFilt^\Sigma (\mathcal{X})}}
        & \simeq
        \pi_* \circ \mathrm{ev}^* (\mathbb{T}_\mathcal{X}) \ ,
        \\
        \mathbb{T}_{\smash{\dFilt^\Sigma (\mathcal{X})}}
        |_{\smash{\dGrad^{\Lambda_\Sigma} (\mathcal{X})}}
        & \simeq
        \bigl(
            \mathbb{T}_\mathcal{X}
            |_{\smash{\dGrad^{\Lambda_\Sigma} (\mathcal{X})}}
        \bigr)_{\Sigma, +} \ .
    \end{align*}
    We chose to state the theorem using cotangent complexes,
    because of the technical issue that
    taking the dual complex may not be a faithful functor,
    as the cotangent complex may not be perfect,
    so the statement about tangent complexes is weaker
    than \cref{thm-tangent-filt-sigma}.
\end{para}

The following lemma will be useful in producing
isomorphisms between connected components of stacks of cone filtrations,
as mentioned in \cref{para-deformation-theory-intro}.

\begin{lemma}
    \label{lemma-retract-etale-implies-iso}
    Let $\mathcal{X}_1, \mathcal{X}_2, \mathcal{Y}$
    be algebraic stacks, and let
    \begin{equation*}
        \begin{tikzcd}[column sep={3em,between origins}, row sep=1.5em]
            \mathcal{X}_1 \ar[rr, "f"] \ar[dr, "p_1"' pos=.4] &&
            \mathcal{X}_2 \ar[dl, "p_2" pos=.4] \\
            & \mathcal{Y}
        \end{tikzcd}
    \end{equation*}
    be a commutative diagram,
    such that $p_1$ and $p_2$ are $\mathbb{A}^1$-action retracts,
    and~$f$ is $\mathbb{A}^1$-equivariant and étale.
    Then $f$ is an isomorphism.
\end{lemma}

\begin{proof}
    We first prove the lemma assuming that $f$ is a monomorphism. In this case it is enough to show that $f$ is surjective. Let~$K$ be an algebraically closed field,
    and let $x \in \mathcal{X}_2 (K)$ be a $K$-point.
    Form the pullback diagram
    \begin{equation*}
        \begin{tikzcd}
            U \ar[r, "{\smash[t]{g'}}"] \ar[d, "f'"']
            \ar[dr, phantom, pos=.2, "\ulcorner"] &
            \mathcal{X}_1 \ar[d, "f"] \\
            \mathbb{A}^1_K \ar[r, "g"] &
            \mathcal{X}_2 \rlap{ ,}
        \end{tikzcd}
    \end{equation*}
    where~$g$ is given by $g (t) = t \cdot x$.
    Then~$U$ admits an $\mathbb{A}^1$-action retract to $\Spec K$,
    and~$f'$ is $\mathbb{A}^1$-equivariant and étale.
    It is enough to show that~$f'$ is surjective. But the image of $f'$ is a $\Gm$-equivariant
    open subset of~$\mathbb{A}^1_K$ containing $0$, and hence all of~$\mathbb{A}^1_K$.

    Now, we show the lemma assuming that $f$ is representable. Then the diagonal $\Delta_f$ is an étale $\mathbb A^1$-equivariant monomorphism with source and target admitting compatible $\mathbb A^1$-action retracts to $\mathcal Y$. Therefore, $\Delta_f$ is an isomorphism, and thus $f$ is a monomorphism, which implies in turn that $f$ is an isomorphism.

    Finally, for general $f$, the representable case applied to the diagonal $\Delta_f$ implies that $\Delta_f$ is an isomorphism, and we conclude in the same way that $f$ is an isomorphism.
\end{proof}

This lemma has the following useful corollary.

\begin{theorem}
    \label{thm-gr-etale-pb}
    Let $\mathcal{X}, \mathcal{Y}$ be stacks as in \cref{assumption-stack-basic},
    and let $f \colon \mathcal{Y} \to \mathcal{X}$ be an étale morphism.
    Then for any integral cone~$\Sigma$ and any rational cone~$C$,
    we have pullback squares
    \begin{equation*}
        \begin{tikzcd}[column sep=1.5em]
            \Filt^\Sigma (\mathcal{Y}) \ar[r] \ar[d, "\mathrm{gr}"']
            \ar[dr, phantom, pos=.2, "\ulcorner"]
            & \Filt^\Sigma (\mathcal{X}) \ar[d, "\mathrm{gr}"]
            \\ \Grad^{\Lambda_\Sigma} (\mathcal{Y}) \ar[r]
            & \Grad^{\Lambda_\Sigma} (\mathcal{X}) \rlap{\textnormal{ ,}}
        \end{tikzcd}
        \qquad
        \begin{tikzcd}[column sep=1.5em]
            \Filt^C (\mathcal{Y}) \ar[r] \ar[d, "\mathrm{gr}"']
            \ar[dr, phantom, pos=.2, "\ulcorner"]
            & \Filt^C (\mathcal{X}) \ar[d, "\mathrm{gr}"]
            \\ \Grad^{F_C} (\mathcal{Y}) \ar[r]
            & \Grad^{F_C} (\mathcal{X}) \rlap{\textnormal{ ,}}
        \end{tikzcd}
    \end{equation*}
    where the horizontal morphisms are induced by~$f$.
\end{theorem}

\begin{proof}
    It is enough to prove the first pullback square.
    Computing the cotangent complexes using \cref{thm-tangent-filt-sigma},
    we see that the horizontal morphisms are étale, so the morphism
    $\Filt^\Sigma (\mathcal{Y}) \to \Filt^\Sigma (\mathcal{X})
    \times_{\Grad^{\smash{\Lambda_\Sigma}} (\mathcal{X})} \Grad^{\Lambda_\Sigma} (\mathcal{Y})$
    is also étale.
    Moreover, both sides of this morphism are $\Theta$-action retracts
    over~$\Grad^{\Lambda_\Sigma} (\mathcal{Y})$,
    and this morphism is $\Theta$-equivariant,
    so the theorem follows from \cref{lemma-retract-etale-implies-iso}.
\end{proof}

\subsection{Cotangent weights}
\label{subsec-cotangent-weights}

\begin{para}
    We discuss \emph{cotangent weights} and
    \emph{cotangent arrangements} of a stack.
    These are another key tool to generalize our observations
    for a linear quotient stack $V / G$
    discussed in \cref{para-idea-filt-cones}
    to general stacks.
    In particular, the cotangent arrangement
    can be seen as a generalization
    of the hyperplane arrangement $\Phi_{V / G}$
    defined in \cref{eg-special-face-linear-quotient}
    to general stacks.

    These notions will be used to study \emph{special cones}
    in \cref{subsec-special-cones} below.
\end{para}

\begin{para}[Cotangent arrangements]
    \label{para-cotangent-arrangement}
    Let~$\mathcal{X}$ be a stack as in \cref{assumption-stack-basic}.
    For a face~$(F, \alpha) \in \mathsf{Face} (\mathcal{X})$,
    consider the complex $\mathbb{L}_\alpha = \mathrm{tot}_\alpha^* (\mathbb{L}_{\mathcal{X} / S})$
    on~$\mathcal{X}_\alpha$, which admits a canonical $F^\vee$-grading,
    and consider the truncation
    $\mathbb{L}_\alpha^{\smash{[0, 1]}} = \tau^{\geq 0} (\mathbb{L}_\alpha)$.
    The set of \emph{cotangent weights} of~$\mathcal{X}$ in the face~$\alpha$
    is defined as
    \begin{equation*}
        W^- (\mathcal{X}, \alpha) =
        \{ \lambda \in F^\vee \mid
        (\mathbb{L}_\alpha^{\smash{[0, 1]}})_\lambda \nsimeq 0 \} \ ,
    \end{equation*}
    and the set of \emph{tangent weights} is
    $W^+ (\mathcal{X}, \alpha) =
    -W^- (\mathcal{X}, \alpha) \subset F^\vee$.

    For a face~$(F, \alpha) \in \mathsf{Face} (\mathcal{X})$,
    define the following hyperplane arrangement
    and cone arrangement (see \cref{para-cone-arrangement}) in~$F$:
    \begin{alignat*}{2}
        \Phi^\mathrm{cot} (\mathcal{X}, \alpha)
        & = \bigl\{
            & \lambda^\perp
            & =
            \{ v \in F \mid \lambda (v) = 0 \}
            \bigm|
            \lambda \in W^+ (\mathcal{X}, \alpha) \setminus \{ 0 \}
        \bigr\} \ ,
        \\
        \Psi^\mathrm{cot} (\mathcal{X}, \alpha)
        & = \bigl\{
            & C_I
            & =
            \{ v \in F \mid \lambda (v) \geq 0 \text{ for all } \lambda \in I \}
            \bigm|
            I \subset W^+ (\mathcal{X}, \alpha)
        \bigr\} \ ,
    \end{alignat*}
    where we allow $\Phi^\mathrm{cot} (\mathcal{X}, \alpha)$
    to be infinite,
    but we only define $\Psi^\mathrm{cot} (\mathcal{X}, \alpha)$
    when~$\mathcal{X}$ has finite cotangent weights as in \cref{para-finite-cotangent-weights} below,
    to ensure that its elements are polyhedral cones.

    We call these the \emph{cotangent arrangement}
    and the \emph{cotangent cone arrangement} of~$\mathcal{X}$ on~$\alpha$.
    We also naturally regard elements of
    $\Psi^\mathrm{cot} (\mathcal{X}, \alpha)$
    as objects of $\mathsf{Cone} (\mathcal{X})$.

    Cotangent arrangements in a subspace of a given face
    refine those on the original face restricted to the subspace.
    More precisely, for a subspace $E \subset F$,
    we have $\Phi^\mathrm{cot} (\mathcal{X}, \alpha) |_E
    \subset \Phi^\mathrm{cot} (\mathcal{X}, \alpha |_E)$
    and $\Psi^\mathrm{cot} (\mathcal{X}, \alpha) |_E
    \subset \Psi^\mathrm{cot} (\mathcal{X}, \alpha |_E)$
    whenever they are defined,
    where the restriction of a hyperplane or cone arrangement to~$E$
    consists of the intersections of each element with~$E$
    and, in the case of hyperplanes, we exclude the face~$E$ itself.
\end{para}

\begin{para}[Stacks with finite cotangent weights]
    \label{para-finite-cotangent-weights}
    We say that a stack~$\mathcal{X}$ has \emph{finite cotangent weights},
    if for all $(F, \alpha) \in \mathsf{Face} (\mathcal{X})$,
    the set~$W^- (\mathcal{X}, \alpha) / {\sim}$ is finite,
    where~$\sim$ denotes positive scaling.

    For example, every quasi-compact stack with \emph{quasi-compact graded points}
    in the sense of \cref{para-quasi-compact-graded-points} below
    has finite cotangent weights.
    In particular, this is the case for any quasi-compact quotient stack.
    Also, every linear moduli stack in the sense of \cref{para-linear-moduli-stacks}
    has finite cotangent weights, as explained in \cref{para-lms-special-faces}.

    Additionally,
    we show below that the mapping stack into a stack with finite cotangent weights
    also has finite cotangent weights.
\end{para}

\begin{example}[Mapping stacks]
    Suppose that $T=\Spec(A)$ for a G-ring $A$ that admits a dualizing complex and that $\pi \colon \cY \to T$ is a flat and formally proper morphism in the sense of~\cite{halpern-leistner-preygel-2023-mapping-stacks} with geometrically connected fibres. Let $\cX$ be a stack locally of finite presentation over~$T$ and with finite cotangent weights. We claim that the stack $\mathrm{Map}_T(\cY,\cX)$, which is algebraic by \cite[Theorem~5.1.1]{halpern-leistner-preygel-2023-mapping-stacks}, also has finite cotangent weights. Note though that it will not automatically satisfy the conditions in \cref{assumption-stack-basic} unless $\cX$ has quasi-affine diagonal.
    
    By a formal manipulation of the functors of points, we have
    \[\mathrm{Grad}^n(\mathrm{Map}(\cY,\cX)) \simeq \mathrm{Map}(\cY, \mathrm{Grad}^n(\cX)) \ .\]
    If $\cM$ denotes a connected component of this stack, then we claim $\cM \times_T \cY$ is connected. To see this, first note that because $\pi_\cM : \cM \times_T \cY \to \cM$ is flat and universally closed by \cite[Proposition~2.4.5 ]{halpern-leistner-preygel-2023-mapping-stacks}, $\cM = \mathrm{im}(\pi_\cM) \cup \cM \setminus \mathrm{im}(\pi_\cM)$ is a decomposition of $\cM$ into disjoint open substacks, and therefore $\pi_\cM$ must be surjective because $\cM$ is connected. Next assume that $\cM \times_T \cY = \cU \cup \cV$ is a decomposition into disjoint open substacks. The fact that $\pi_\cM$ has connected geometric fibers implies that $\cU = \pi_\cM^{-1}(\pi_\cM(\cU))$ and likewise for $\cV$. It follows that $\cM = \pi_\cM(\cU) \cup \pi_\cM(\cV)$ is a decomposition into disjoint open substacks, which implies that either $\cU$ or $\cV$ is empty. Therefore $\cM \times_{T} \cY$ is connected.
    
    Because $\cM \times_T \cY$ is connected, the universal morphism factors as $\mathrm{ev} \colon \cM \times \cY \to \cX_\alpha$ for some $\alpha \colon \bZ^n \to \mathrm{CL}(\cX)$. The pullback of the cotangent complex of the mapping stack to $\cM$ is $\pi_+(\mathrm{ev}^\ast(\mathbb{L}_\alpha))$, where $\pi_+$ is the left adjoint of the pullback $\pi^\ast$ on derived categories of quasi-coherent sheaves. The pullback functor $\pi^\ast$ preserves coconnective complexes because $\pi$ is flat, so its left adjoint $\pi_+$ preserves connective complexes. It follows that
    \[
        \bigl(
            \pi_+ \circ \mathrm{ev}^\ast(\mathbb{L}_\alpha)
        \bigr)^{[0,1]}
        \simeq
        \bigl(
            \pi_+ \circ \mathrm{ev}^\ast(\mathbb{L}_\alpha^{[0,1]})
        \bigr)^{[0,1]} \ ,
    \]
    so up to positive scaling only finitely many weights appear in this complex.
\end{example}

\begin{para}[Pathological examples]
    \label{para-pathological-examples}
    We provide two pathological examples,
    hopefully to illustrate some of the subtleties
    involved in the results in \cref{subsec-special-cones} and below.

    Firstly, a quasi-compact stack need not have finite cotangent weights.
    Consider the following example:
    Let $\bG_\mathrm{m}^2$ act linearly on $\bA^1$,
    corresponding to some character $\chi \in \bZ^2$,
    and let $\phi \colon \bG_\mathrm{m}^2 \simto \bG_\mathrm{m}^2$
    be an automorphism such that the sequence
    $(\chi, \phi^* (\chi), \phi^{*2} (\chi),\ldots)$
    is not periodic.
    Let $\cX$ be the stack obtained by gluing $\bA^1 / \bG_m^2$
    to $\bP^1 \times \mathrm{B} \bG_\mathrm{m}^2$
    along the points $0 \in \bA^1$ and $1 \in \bP^1$ respectively,
    and then gluing $0$ to $\infty$ in $\bP^1$,
    identifying the fibres
    $\{0\} \times \mathrm{B} \bG_\mathrm{m}^2 \simeq \{\infty\} \times \mathrm{B} \bG_\mathrm{m}^2$
    via the automorphism $\phi$.
    Then~$\mathcal{X}$ does not have finite cotangent weights.

    Secondly, there is a stack~$\mathcal{X}$
    with $\mathrm{CL}_\mathbb{Q} (\mathcal{X}) \simeq \mathbb{Q}^2$,
    such that it has trivial cotangent arrangements on all faces,
    and it has quasi-compact graded points as in
    \cref{para-quasi-compact-graded-points},
    while every linear subspace of~$\mathbb{Q}^2$ is a special face.
    Indeed, consider the affine group scheme $G \to \mathbb{A}^1$
    with underlying scheme $\{ (x, t) \in \mathbb{A}^2 \mid x t \neq 1 \}$
    and multiplication $(x, t) \cdot (y, t) = (x + y - t x y, t)$.
    Consider the stack
    $\mathcal{G} = \mathrm{B} (G \times \mathbb{G}_\mathrm{m})$
    defined over~$\mathbb{A}^1$.
    Let~$\mathcal{X}$ be the stack obtained by gluing
    infinitely many copies of~$\mathcal{G}$ into a chain,
    identifying the fibre
    $\{ 1 \} \times \mathrm{B} \mathbb{G}_\mathrm{m}^2$
    in the $n$-th copy with the fibre
    $\{ -1 \} \times \mathrm{B} \mathbb{G}_\mathrm{m}^2$
    in the $(n+1)$-th copy using an automorphism
    $\phi_n \colon \mathbb{G}_\mathrm{m}^2 \simto \mathbb{G}_\mathrm{m}^2$,
    such that the sequence $(L, \phi_1 (L), \phi_2 \phi_1 (L), \ldots)$
    traverses all lines in~$\mathbb{Q}^2$,
    where $L = \{ 0 \} \times \mathbb{Q} \subset \mathbb{Q}^2$,
    and each $\phi_n$ induces an automorphism of
    $\mathbb{Q}^2 \simeq \mathrm{CL}_\mathbb{Q} (\mathrm{B} \mathbb{G}_\mathrm{m}^2)$.
\end{para}

\begin{lemma}
    \label{lemma-cot-arr-special-faces}
    Let~$\mathcal{X}$ be a stack as in \cref{assumption-stack-basic},
    and let~$(F, \alpha) \in \mathsf{Face}^\mathrm{sp} (\mathcal{X})$.
    Then for any linear subspace $E \subset F$
    which is an intersection of hyperplanes in $\Phi^\mathrm{cot} (\mathcal{X}, \alpha)$,
    the face $(E, \alpha |_E)$ is special.
\end{lemma}

\begin{proof}
    By \cref{thm-special-face-closure},
    the collection of subspaces of~$F$ that are special faces
    form a face arrangement in~$F$,
    and hence is closed under taking intersections.
    Therefore, it is enough to show that
    any hyperplane $E \subset F$
    in~$\Phi^\mathrm{cot} (\mathcal{X}, \alpha)$ is a special face.

    Suppose the contrary, and write $\beta = \alpha |_E$.
    Then the special closure of the face $(E, \beta)$
    has to be $(F, \alpha)$,
    and we have an induced isomorphism
    $\mathcal{X}_\alpha \simto \mathcal{X}_\beta$.
    Consider the complexes~$\mathbb{L}_\alpha$ and~$\mathbb{L}_\beta$
    defined in \cref{para-cotangent-arrangement},
    with a canonical $F^\vee$- and $E^\vee$-grading, respectively.
    By \cref{thm-tangent-filt-sigma},
    or by \cite[Lemma~1.2.2]{halpern-leistner-instability},
    we have isomorphisms
    $\mathrm{H}^i ((\mathbb{L}_\alpha)_0) \simeq
    \mathrm{H}^i (\mathbb{L}_{\mathcal{X}_\alpha}) \simeq
    \mathrm{H}^i ((\mathbb{L}_\beta)_0)$
    for $i = 0, 1$ of cohomology sheaves on~$\mathcal{X}_\alpha$,
    where $(-)_0$ denotes the $0$-th graded piece.
    But $(\mathbb{L}_\beta)_0 \simeq \bigoplus_{\lambda \in E^\perp} (\mathbb{L}_\alpha)_\lambda$,
    where $E^\perp \subset F^\vee$ is the annihilator of~$E$,
    which implies that $(\mathbb{L}_\alpha)_\lambda \simeq 0$
    for all $\lambda \in E^\perp \setminus \{ 0 \}$,
    so that $E \notin \Phi^\mathrm{cot} (\mathcal{X}, \alpha)$,
    a contradiction.
\end{proof}

\subsection{Special cones}
\label{subsec-special-cones}

\begin{para}
    We introduce the notion of \emph{special cones} of a stack,
    similarly to special faces introduced in \cref{sec-faces}.
    They are cones~$\sigma$ in the component lattice
    such that one cannot enlarge~$\sigma$
    without changing either the stack~$\mathcal{X}_\sigma^+$
    or the stack~$\mathcal{X}_{\mathrm{span} (\sigma)}$
    defined in \cref{para-notation-x-alpha,para-x-sigma-plus}.
    Roughly, they are the cones that arise in the wall-and-chamber structure
    on the component lattice described in \cref{para-intro-component-lattice}.

    Under a mild condition, special cones can be described explicitly
    using special faces and the cotangent arrangement
    introduced in \cref{para-cotangent-arrangement},
    by \cref{thm-cotangent-arrangement}.
    This result will then be used to show that
    special cones form a cone arrangement
    in the sense of \cref{para-cone-arrangement},
    in \cref{thm-special-cone-closure}.
\end{para}

\begin{para}[Special cones]
    Let~$\mathcal{X}$ be a stack as in \cref{assumption-stack-basic}.

    A \emph{special cone} of~$\mathcal{X}$ is a cone
    $\sigma \in \mathsf{Cone}^\mathrm{nd} (\mathcal{X})$,
    such that its span~$\alpha$ is a special face,
    and for any cone $\sigma' \subset \alpha$ containing~$\sigma$,
    if the induced morphism
    $\mathcal{X}_{\smash{\sigma'}}^+ \to \mathcal{X}_\sigma^+$
    is an isomorphism, then $\sigma = \sigma'$.
    We denote by
    \begin{equation*}
        \mathsf{Cone}^\mathrm{sp} (\mathcal{X}) \subset
        \mathsf{Cone}^\mathrm{nd} (\mathcal{X})
    \end{equation*}
    the full subcategory of special cones.

    Equivalently, a cone $\sigma \in \mathsf{Cone}^\mathrm{nd} (\mathcal{X})$
    is special if and only if for any morphism $f \colon \sigma \to \sigma'$
    in $\mathsf{Cone}^\mathrm{nd} (\mathcal{X})$,
    if the induced morphisms
    $\mathcal{X}_{\alpha'} \to \mathcal{X}_\alpha$
    and $\mathcal{X}_{\smash{\sigma'}}^+ \to \mathcal{X}_\sigma^+$
    are isomorphisms,
    where $\alpha = \mathrm{span} (\sigma)$ and $\alpha' = \mathrm{span} (\sigma')$,
    then $f$ is an isomorphism.

    For example, all special faces of~$\mathcal{X}$ are special cones, that is,
    $\mathsf{Face}^\mathrm{sp} (\mathcal{X}) \subset
    \mathsf{Cone}^\mathrm{sp} (\mathcal{X})$.
\end{para}

\begin{example}[Linear quotient stacks]
    Let $\mathcal{X} = V / G$ be a linear quotient stack,
    as in \cref{eg-component-lattice-linear-quotient,eg-special-face-linear-quotient},
    so that $\mathrm{CL}_\mathbb{Q} (\mathcal{X}) \simeq
    (\Lambda_T \otimes \mathbb{Q}) / W$.
    Assume that~$G$ is connected.

    Let $W^+ \subset \Lambda^T$ be the set of weights in~$V$ and roots of~$G$.
    By the explicit description of the stacks
    $\mathcal{X}_\alpha$ and $\mathcal{X}_\sigma^+$
    in \cref{eg-filt-cones-quotient-stack},
    the special cones of~$\mathcal{X}$
    are $W$-orbits of cones
    $C \subset \Lambda_T \otimes \mathbb{Q}$
    which are intersections of the following two types of cones:
    \begin{itemize}
        \item
            hyperplanes of the form
            $\lambda^\perp = \{ v \in \Lambda_T \otimes \mathbb{Q} \mid \lambda (v) = 0 \}$
            for $\lambda \in W^+ \setminus \{ 0 \}$, and
        \item
            half spaces of the form
            $\lambda^+ = \{ v \in \Lambda_T \otimes \mathbb{Q} \mid \lambda (v) \geq 0 \}$
            for $\lambda \in W^+ \setminus \{ 0 \}$.
    \end{itemize}
    We also allow the empty intersection, which is the whole space.

    For example, if $\mathcal{X} = \mathbb{A}^1 / \mathbb{G}_m$
    with the scaling action,
    then $\mathrm{CL}_\mathbb{Q} (\mathcal{X}) \simeq \mathbb{Q}$
    and $W^+ = \{ 0, 1 \} \subset \mathbb{Z}$,
    so the special cones of~$\mathcal{X}$ are
    $\{ 0 \}$, $\mathbb{Q}_{\geq 0}$, and~$\mathbb{Q}$.
\end{example}

\begin{lemma}
    \label{lemma-cotangent-arrangement-cones}
    Let~$\mathcal{X}$ be a stack as in \cref{assumption-stack-basic},
    with finite cotangent weights,
    and let $(F, \alpha) \in \mathsf{Face} (\mathcal{X})$ be a face.
    Then for any cone $(C, \sigma) \subset (F, \alpha)$ of full dimension,
    writing~$(C^\mathrm{cot}, \sigma^\mathrm{cot})$ for the minimal cone
    in $\Psi^\mathrm{cot} (\mathcal{X}, \alpha)$
    containing~$C$,
    we have an induced isomorphism
    \begin{equation*}
        \mathcal{X}_{\smash{\sigma^\mathrm{cot}}}^+ \longsimto
        \mathcal{X}_\sigma^+ \ .
    \end{equation*}
    Conversely, if there are cones
    $(C, \sigma) \subset (C', \sigma') \subset (F, \alpha)$
    of full dimension, such that the induced morphism
    $\mathcal{X}_{\smash{\sigma'}}^+ \to \mathcal{X}_\sigma^+$
    is an isomorphism, then
    $C^\mathrm{cot} = (C')^\mathrm{cot}$.
\end{lemma}

\begin{proof}
    By \cref{lemma-theta-retract},
    the stacks $\mathcal{X}_\sigma^+$ and~$\mathcal{X}_{\smash{\sigma^\mathrm{cot}}}^+$
    are $\Theta$-action retracts onto~$\mathcal{X}_\alpha$,
    so that by \cref{lemma-retract-etale-implies-iso},
    it is enough to show that the induced morphism
    $f \colon \mathcal{X}_{\smash{\sigma^\mathrm{cot}}}^+ \to \mathcal{X}_\sigma^+$
    is étale.
    By the definition of $C^\mathrm{cot}$,
    each tangent weight in $W^+ (\mathcal{X}, \alpha)$
    pairs non-negatively with~$C$
    if and only if it pairs non-negatively with~$C^\mathrm{cot}$.
    Therefore, writing~$\dmathcal{X}_\sigma^+$
    and~$\dmathcal{X}_{\smash{\sigma^\mathrm{cot}}}^+$
    for the derived versions,
    by \cref{thm-tangent-filt-sigma},
    the induced map on cotangent complexes
    $\mathrm{H}^i (f^* \mathbb{L}_{\dmathcal{X}_\sigma^+}\vert_{\cX_\alpha}) \to
    \mathrm{H}^i (\mathbb{L}_{\dmathcal{X}_{\smash{\sigma^\mathrm{cot}}}^+}\vert_{\cX_\alpha})$
    is an isomorphism for $i = 0, 1$, where we restrict along the split filtration morphism
    $\mathcal{X}_\alpha \to \mathcal{X}_{\smash{\sigma^\mathrm{cot}}}^+$. Since, by \cref{thm-tangent-filt-sigma} again, the morphism $f^* \mathbb{L}_{\dmathcal{X}_\sigma^+}\vert_{\cX_\alpha} \to
    \mathbb{L}_{\dmathcal{X}_{\smash{\sigma^\mathrm{cot}}}^+}\vert_{\cX_\alpha}$ is a split surjection, the map
    $\mathrm{H}^{-1} (f^* \mathbb{L}_{\dmathcal{X}_\sigma^+}\vert_{\cX_\alpha}) \to
    \mathrm{H}^{-1} (\mathbb{L}_{\dmathcal{X}_{\smash{\sigma^\mathrm{cot}}}^+}\vert_{\cX_\alpha})$ is surjective.
    We thus have
    $\mathrm{H}^i (\mathbb{L}_{\dmathcal{X}_{\smash{\sigma^\mathrm{cot}}}^+
    / \dmathcal{X}_\sigma^+}|_{\mathcal{X}_\alpha}) \simeq 0$
    for $i = -1, 0, 1$ and, by right exactness of pullback $(-)|_{\mathcal{X}_\alpha}$, also $\mathrm{H}^i (\mathbb{L}_{\dmathcal{X}_{\smash{\sigma^\mathrm{cot}}}^+
    / \dmathcal{X}_\sigma^+})|_{\mathcal{X}_\alpha} \simeq 0$
    for $i = -1, 0, 1$.
    By \cref{lemma-theta-retract},
    any open substack of~$\mathcal{X}_{\smash{\sigma^\mathrm{cot}}}^+$
    containing the split filtrations is the whole stack,
    so we have
    $\mathrm{H}^i (\mathbb{L}_{\dmathcal{X}_{\smash{\sigma^\mathrm{cot}}}^+
    / \dmathcal{X}_\sigma^+}) \simeq 0$
    for $i = -1, 0, 1$, and hence
    $\mathrm{H}^i (\mathbb{L}_{\mathcal{X}_{\smash{\sigma^\mathrm{cot}}}^+
    / \mathcal{X}_\sigma^+}) \simeq 0$
    for $i = -1, 0, 1$, so that~$f$ is étale.

    For the final statement,
    note that we have isomorphisms
    $\mathrm{H}^i (\mathbb{L}_{\dmathcal{X}_\sigma^+}|_{\mathcal{X}_\alpha}) \simto
    \mathrm{H}^i (\mathbb{L}_{\dmathcal{X}_{\smash{\sigma'}}^+} |_{\mathcal{X}_\alpha})$
    for $i = 0, 1$,
    so by \cref{thm-tangent-filt-sigma},
    each tangent weight in $W^+ (\mathcal{X}, \alpha)$
    pairs non-negatively with~$C$
    if and only if it pairs non-negatively with~$C'$,
    which implies that $C^\mathrm{cot} = (C')^\mathrm{cot}$.
\end{proof}

\begin{theorem}
    \label{thm-cotangent-arrangement}
    Let~$\mathcal{X}$ be a stack as in \cref{assumption-stack-basic},
    with finite cotangent weights,
    and let~$(F, \alpha) \in \mathsf{Face}^\mathrm{sp} (\mathcal{X})$
    be a special face.
    Then we have the following:

    \begin{enumerate}[resume]
        \item
            \label{item-cotangent-arrangement-1}
            Every cone $(C, \sigma) \in \Psi^\mathrm{cot} (\mathcal{X}, \alpha)$
            is special.

        \item
            \label{item-cotangent-arrangement-2}
            A cone $(C, \sigma) \subset (F, \alpha)$
            of full dimension is special
            if and only if\/ $(C, \sigma) \in \Psi^\mathrm{cot} (\mathcal{X}, \alpha)$.

        \item
            \label{item-cotangent-arrangement-3}
            In particular, all boundary cones (of all dimensions)
            of a special cone are special cones.
    \end{enumerate}
\end{theorem}

\begin{proof}
    For \cref{item-cotangent-arrangement-1},
    note that~$C$ belongs to the cotangent cone arrangement
    in its span, which is a special face by
    \cref{lemma-cot-arr-special-faces}.
    Therefore, the statement follows from~\cref{item-cotangent-arrangement-2}.

    The statement \cref{item-cotangent-arrangement-2}
    is a direct consequence of \cref{lemma-cotangent-arrangement-cones}.

    For \cref{item-cotangent-arrangement-3},
    note that for a special face~$(F, \alpha)$
    and a cone $(C, \sigma) \in \Psi^\mathrm{cot} (\mathcal{X}, \alpha)$,
    each boundary cone of~$C$
    belongs to the cotangent cone arrangement in its own span,
    and that the span is a special face by
    \cref{lemma-cot-arr-special-faces}.
    Therefore, all boundary cones of~$C$ are special cones.
\end{proof}

\begin{theorem}
    \label{thm-special-cone-closure}
    Let~$\mathcal{X}$ be a stack as in \cref{assumption-stack-basic},
    and suppose that~$\mathcal{X}$ has finite cotangent weights
    in the sense of\/~\cref{para-finite-cotangent-weights}.
    Then there is a special cone closure functor
    \begin{equation*}
        (-)^\mathrm{sp} \colon
        \mathsf{Cone} (\mathcal{X}) \longrightarrow
        \mathsf{Cone}^\mathrm{sp} (\mathcal{X}) \ ,
    \end{equation*}
    which is left adjoint to the inclusion
    $\mathsf{Cone}^\mathrm{sp} (\mathcal{X}) \hookrightarrow
    \mathsf{Cone} (\mathcal{X})$.
    Its restriction to $\mathsf{Face} (\mathcal{X})$
    agrees with the special face closure
    in \cref{thm-special-face-closure}.

    In particular, special cones form a cone arrangement
    in the sense of\/~\cref{para-cone-arrangement},
    which we call the special cone arrangement of~$\mathcal{X}$.

    Moreover, for any $\sigma \in \mathsf{Cone} (\mathcal{X})$, we have
    $\mathrm{span} (\sigma^\mathrm{sp}) \simeq
    \mathrm{span} (\sigma)^\mathrm{sp}$,
    and the adjunction unit
    $\sigma \to \sigma^\mathrm{sp}$
    induces isomorphisms
    \begin{align*}
        \mathcal{X}_{\smash{\sigma^\mathrm{sp}}}^+
        & \longsimto
        \mathcal{X}_\sigma^+ \ ,
        \\
        \mathcal{X}_{\mathrm{span} (\sigma^\mathrm{sp})}
        & \longsimto
        \mathcal{X}_{\mathrm{span} (\sigma)} \ .
    \end{align*}
\end{theorem}

\begin{proof}
    For a cone
    $(C, \sigma) \in \mathsf{Cone} (\mathcal{X})$,
    let $(F, \alpha) = \mathrm{span} (C, \sigma)$,
    and let $(F^\mathrm{sp}, \alpha^\mathrm{sp})$
    be the special face closure of~$(F, \alpha)$.
    Define the cone
    $(C^\mathrm{sp}, \sigma^\mathrm{sp}) \subset
    (F^\mathrm{sp}, \alpha^\mathrm{sp})$
    to be the minimal cone in
    $\Psi^\mathrm{cot} (\mathcal{X}, \alpha^\mathrm{sp})$
    containing the image of~$\sigma$, which is of full dimension by \cref{lemma-cot-arr-special-faces}. The cone
    $(C^\mathrm{sp}, \sigma^\mathrm{sp})$
    is special by \cref{thm-cotangent-arrangement}.

    We show that there is an induced isomorphism
    $\mathcal{X}_{\smash{\sigma^\mathrm{sp}}}^+ \simto \mathcal{X}_\sigma^+$.
    Indeed, by
    \cref{lemma-theta-retract,thm-special-face-closure},
    both sides are $\Theta$-action retracts onto
    $\mathcal{X}_{\alpha^\mathrm{sp}} \simeq \mathcal{X}_\alpha$,
    and by \cref{lemma-retract-etale-implies-iso},
    it is enough to show that the morphism
    $\mathcal{X}_{\smash{\sigma^\mathrm{sp}}}^+ \to \mathcal{X}_\sigma^+$
    is étale.
    This follows from computing the cotangent complexes at split filtrations,
    as in the proof of \cref{lemma-cotangent-arrangement-cones}.

    To prove the adjunction property, let
    $(C', \sigma') \in \mathsf{Cone}^\mathrm{sp} (\mathcal{X})$
    be a special cone, and let
    $f \colon \sigma \to \sigma'$
    be a morphism.
    We need to show that~$f$ factors uniquely
    through the canonical map
    $\sigma \to \sigma^\mathrm{sp}$.
    Write $(F, \alpha) = \mathrm{span} (C, \sigma)$
    and $(F', \alpha') = \mathrm{span} (C', \sigma')$.
    By \cref{thm-special-face-closure},
    the induced map~$f \colon \alpha \to \alpha'$ factors uniquely
    through the canonical map $\alpha \to \alpha^\mathrm{sp}$,
    inducing a map
    $(F^\mathrm{sp}, \alpha^\mathrm{sp}) \to (F', \alpha')$,
    which must be injective because $\alpha^{\rm{sp}}$ is non-degenerate.
    This identifies $F^\mathrm{sp}$ with a subspace of~$F'$,
    and it is enough to show that~$C^\mathrm{sp} \subset C'$
    under this identification.
    This is because by \cref{thm-cotangent-arrangement},
    we have $C' \in \Psi^\mathrm{cot} (\mathcal{X}, \alpha')$,
    so the cone
    $C'' = C' \cap C^\mathrm{sp} \subset F^\mathrm{sp}$
    belongs to
    $\Psi^\mathrm{cot} (\mathcal{X}, \alpha^\mathrm{sp})$
    and contains the image of~$C$,
    and hence $C'' = C^\mathrm{sp}$.
\end{proof}

We also have the following version of
\cref{thm-special-cone-closure}
relative to a face:

\begin{corollary}
    \label{cor-special-cone}
    Let~$\mathcal{X}$ be a stack as in \cref{assumption-stack-basic},
    with finite cotangent weights,
    and let~$(F, \alpha) \in \mathsf{Face} (\mathcal{X})$.
    Then for any cone $(C, \sigma) \subset (F, \alpha)$,
    there exists a unique maximal cone
    $(\bar{C}, \bar{\sigma}) \subset (F, \alpha)$
    containing~$C$,
    such that the induced morphisms
    $\mathcal{X}_{\smash{\bar{\sigma}}}^+ \to \mathcal{X}_\sigma^+$
    and $\mathcal{X}_{\mathrm{span} (\bar{\sigma})}
    \to \mathcal{X}_{\mathrm{span} (\sigma)}$
    are isomorphisms.
    Moreover, if~$C$ is of full dimension in~$F$,
    then~$\bar{C}$ is the minimal element of\/
    $\Psi^\mathrm{cot} (\mathcal{X}, \alpha)$
    \textnormal{(see \cref{para-cotangent-arrangement})}
    containing~$C$.
\end{corollary}

\begin{proof}
    Consider the special face closure
    $(F^\mathrm{sp}, \alpha^\mathrm{sp})$
    of $(F, \alpha)$, and let
    $i \colon F \to F^\mathrm{sp}$
    be the canonical map.
    The special cone closure $\sigma^\mathrm{sp}$
    has a canonical map to $\alpha^\mathrm{sp}$,
    identifying it with a cone of full dimension
    $C^\mathrm{sp} \subset F^\mathrm{sp}$.
    Setting $\bar{C} = i^{-1} (C^\mathrm{sp})$
    gives the desired cone.
    The final claim follows from applying
    \cref{thm-cotangent-arrangement}~\cref{item-cotangent-arrangement-2}
    to $\alpha^\mathrm{sp}$, then restricting to~$\alpha$.
\end{proof}

\begin{para}[Cotangent chambers and special closures of rays]
    \label{para-special-closure-rays}
    Let~$\mathcal{X}$ be a stack as in \cref{assumption-stack-basic},
    with finite cotangent weights.
    We now discuss three closely related types of cones:

    \begin{enumerate}
        \item
            \label{item-chambers}
            \emph{Cotangent chambers},
            i.e., chambers (as in~\cref{para-hyperplane-arrangements})
            in cotangent arrangements
            $\Phi^\mathrm{cot} (\mathcal{X}, \alpha)$
            on special faces~$\alpha$.
        \item
            \label{item-special-closure-chambers}
            Special cones
            that are special closures of cotangent chambers.
        \item
            \label{item-special-closure-rays}
            Special cones
            that are special closures of rays
            $\mathbb{Q}_{\geq 0} \cdot \lambda$
            for $\lambda \in |\mathrm{CL}_\mathbb{Q} (\mathcal{X})|$.
    \end{enumerate}

    Cones~$\sigma$ belonging to type~\cref{item-special-closure-rays}
    have the special property that
    $\mathcal{X}_\sigma^+ \simeq \mathcal{X}_\lambda^+$
    is isomorphic to a connected component of $\mathrm{Filt} (\mathcal{X})$,
    which is not true for all special cones.
    However, cones of types \cref{item-chambers,item-special-closure-chambers}
    often also have this property, as we discuss below.

    By \cref{thm-cotangent-arrangement},
    taking the special closure gives a bijection of cones
    of type~\cref{item-chambers} and type~\cref{item-special-closure-chambers}.
    Additionally, the two types coincide
    when~$\mathcal{X}$ has \emph{symmetric cotangent weights},
    that is, when for any face
    $(F, \alpha) \in \mathsf{Face} (\mathcal{X})$,
    we have $W^+ (\mathcal{X}, \alpha) / {\sim} = W^- (\mathcal{X}, \alpha) / {\sim}$
    as subsets of $(F^\vee \setminus \{ 0 \}) / {\sim}$,
    where~$\sim$ denotes positive scaling.

    We have an inclusion
    \cref{item-special-closure-rays}~$\subset$
    \cref{item-special-closure-chambers}
    by \cref{thm-cotangent-arrangement}.
    The reverse inclusion is also usually true,
    and only fails in the pathological case
    where a cotangent chamber~$\sigma$
    is covered by infinitely many lower-dimensional special faces.
    This can indeed happen,
    as shown in the second example in \cref{para-pathological-examples}.
    However, this pathology cannot occur if~$\mathcal{X}$ is quasi-compact
    and has quasi-compact graded points in the sense of
    \cref{para-quasi-compact-graded-points},
    which follows from the finiteness theorem, \cref{thm-finiteness},
    or for linear moduli stacks
    in the sense of \cref{para-linear-moduli-stacks},
    which follows from \cref{para-lms-special-faces}.

    In particular, if this pathology does not occur,
    then the cones~$\sigma$ of all three types above
    satisfy the property that their special closures
    belong to type~\cref{item-special-closure-rays},
    and hence $\mathcal{X}_\sigma^+ \simeq \mathcal{X}_\lambda^+$
    for some $\lambda \in |\mathrm{CL}_\mathbb{Q} (\mathcal{X})|$.
\end{para}

\begin{para}[Very special cones]
    Finally, we briefly discuss a notion of \emph{very special cones},
    which will not be used in the rest of the paper.

    Let~$\mathcal{X}$ be a stack as in \cref{assumption-stack-basic}.
    We say that a cone $\sigma \in \mathsf{Cone} (\mathcal{X})$
    is \emph{very special},
    if for any morphism $f \colon \sigma \to \sigma'$
    in $\mathsf{Cone} (\mathcal{X})$,
    if the induced morphism
    $\mathcal{X}_{\smash{\sigma'}}^+ \to \mathcal{X}_\sigma^+$
    is an isomorphism, then $f$ is an isomorphism.

    All very special cones are special cones,
    and using an argument
    similar to that of \cref{subsec-special-face-closure},
    one can show that
    for a stack~$\mathcal{X}$ with finite cotangent weights,
    its very special cones also form a cone arrangement.
\end{para}

\section{Main results}

\label{sec-main-results}

\subsection{The constancy theorem}
\label{subsec-constancy-theorem}

\begin{para}
    In this section, we prove the \emph{constancy theorem},
    \cref{thm-constancy},
    which makes precise the idea mentioned in
    \cref{para-intro-component-lattice}
    that the component lattice of a stack~$\mathcal{X}$
    usually has a wall-and-chamber structure,
    such that the isomorphism type of components of
    $\Grad (\mathcal{X})$ and $\Filt (\mathcal{X})$
    do not vary inside a chamber.

    The constancy theorem is a combination of the main results
    in the previous sections,
    \cref{thm-special-face-closure,thm-special-cone-closure},
    and does not introduce any essentially new ideas.
\end{para}

\begin{theorem}
    \label{thm-constancy}
    Let~$\mathcal{X}$ be a stack as in \cref{assumption-stack-basic},
    and suppose that~$\mathcal{X}$ has finite cotangent weights
    in the sense of\/~\cref{para-finite-cotangent-weights}.

    For any $\lambda \in |\mathrm{CL}_\mathbb{Q} (\mathcal{X})|$,
    let $\alpha \in \mathsf{Face}^\mathrm{sp} (\mathcal{X})$
    be the special face closure of
    the face $\mathbb{Q} \cdot \lambda$.
    Then~$\lambda$ is contained in the interior of a chamber
    $\sigma \subset \alpha$
    in the cotangent arrangement
    $\Phi^\mathrm{cot} (\mathcal{X}, \alpha)$
    defined in \cref{para-cotangent-arrangement},
    and we have induced isomorphisms
    \begin{alignat}{2}
        \label{eq-constancy-grad}
        \mathcal{X}_\alpha
        & \longsimto \mathcal{X}_{\lambda}
        && \ ,
        \\
        \label{eq-constancy-filt}
        \mathcal{X}_\sigma^+
        & \longsimto \mathcal{X}_{\smash{\lambda}}^+
        && \ ,
    \end{alignat}
    where $\mathcal{X}_\lambda \subset \mathrm{Grad}_\mathbb{Q} (\mathcal{X})$
    and $\mathcal{X}_{\smash{\lambda}}^+ \subset \mathrm{Filt}_\mathbb{Q} (\mathcal{X})$
    are the connected components corresponding to~$\lambda$.
\end{theorem}

In particular, the theorem implies that
the stacks~$\mathcal{X}_\lambda$ and~$\mathcal{X}_{\smash{\lambda}}^+$
do not change if we vary the point~$\lambda$
inside a chamber of the cotangent arrangement,
apart from a lower dimensional subset
given by special faces strictly contained in~$\alpha$,
hence the name \emph{constancy theorem}.
See, however,
\cref{para-special-closure-rays},
for a pathological situation
where these lower-dimensional faces can cover the whole face~$\alpha$.

\begin{proof}
    The isomorphism~\cref{eq-constancy-grad}
    follows from \cref{thm-special-face-closure}.
    The existence of such a chamber follows from
    \cref{lemma-cot-arr-special-faces},
    and in this case,
    the special cone closure of the cone
    $\mathbb{Q}_{\geq 0} \cdot \lambda$
    must agree with the special cone closure of~$\sigma$,
    by \cref{thm-cotangent-arrangement}.
    The isomorphism~\cref{eq-constancy-filt}
    then follows from \cref{thm-special-cone-closure}.
\end{proof}

\begin{para}[A derived version]
    \label{para-constancy-derived}
    We also discuss a version of \cref{thm-constancy}
    for derived stacks.

    Let~$\mathcal{X}$ be a derived algebraic stack,
    locally finitely presented over a base
    derived algebraic space~$S$,
    such that their classical truncations satisfy the assumptions
    in \cref{assumption-stack-basic}.

    We have a derived version of
    \cref{thm-filt-sigma-image}
    for the stacks $\dFilt^\Sigma (\mathcal{X})$,
    which follows from the classical version and
    \cref{thm-tangent-filt-sigma},
    together with the last part of \cref{lemma-theta-retract},
    which altogether imply that the morphism in question
    is an open and closed immersion on the classical truncations
    and is étale, so it is an open and closed immersion.

    Similarly, a derived version of
    \cref{lemma-retract-etale-implies-iso}
    can be deduced from the classical version.

    We may also define \emph{special faces}
    and \emph{special cones} of~$\mathcal{X}$,
    similarly to the classical case,
    but using the stacks $\dmathcal{X}_\alpha$ and $\dmathcal{X}_\sigma^+$
    instead of~$\mathcal{X}_\alpha$ and~$\mathcal{X}_\sigma^+$.
    We have derived versions of
    \cref{thm-special-face-closure,thm-cotangent-arrangement,thm-special-cone-closure},
    which follow from the same proofs with minor changes,
    using the above derived versions of preliminary results,
    and the definition of cotangent arrangements
    needs to be modified by using the full cotangent complex
    $\mathbb{L}_\alpha = \mathrm{tot}_\alpha^* (\mathbb{L}_{\mathcal{X} / S})$
    in place of $\mathbb{L}_\alpha^{\smash{[0, 1]}}$
    in \cref{para-cotangent-arrangement}.

    Consequently, a derived version of \cref{thm-constancy}
    involving the stacks~$\dmathcal{X}_\alpha$,
    $\dmathcal{X}_\lambda$,
    $\dmathcal{X}_\sigma^+$,
    and $\dmathcal{X}_\lambda^+$,
    also holds using special faces, special cones,
    and cotangent arrangements in the derived sense,
    where we assume~$\mathcal{X}$ to have finite cotangent weights
    in the derived sense using the full cotangent complex.
\end{para}

\subsection{The finiteness theorem}

\begin{para}
    In this section, we prove the \emph{finiteness theorem},
    \cref{thm-finiteness},
    which provides a criterion for a stack
    to have only finitely many special faces.
\end{para}

\begin{para}[Stacks with quasi-compact graded points]
    \label{para-quasi-compact-graded-points}
    Let~$\mathcal{X}$ be a stack as in \cref{assumption-stack-basic}.
    We say that~$\mathcal{X}$ has
    \emph{quasi-compact graded points} if, for any face
    $\alpha \colon F \to \mathrm{CL}_{\mathbb{Q}} (\mathcal{X})$,
    the morphism
    $\mathrm{tot}_\alpha \colon \mathcal{X}_\alpha \to \mathcal{X}$
    is quasi-compact.

    This is a mild condition on the stack~$\mathcal{X}$.
    For example, if~$\mathcal{X}$ is a quotient stack,
    then this is satisfied by the explicit description
    in \cref{eg-grad-quotient-stack}.
    Also, if~$\mathcal{X}$ has a
    \emph{norm on graded points}
    as in \textcite[Definition~4.1.12]{halpern-leistner-instability},
    then it has quasi-compact graded points
    by \cite[Proposition~3.8.2]{halpern-leistner-instability}.
    A norm on graded points exists in many examples.
\end{para}

\begin{theorem}
    \label{thm-finiteness}
    Let $\mathcal{X}$ be a quasi-compact stack
    with quasi-compact graded points.
    Then~$\mathcal{X}$ has only finitely many special faces
    and special cones,
    and there are only finitely many morphisms between them.
\end{theorem}

In particular,
by \cref{thm-special-face-closure,thm-special-cone-closure},
this implies that for such a stack~$\mathcal{X}$, the connected components of the stacks
\begin{alignat*}{4}
    \Grad^n (\mathcal{X}) \ , && \qquad
    \Grad^\Lambda (\mathcal{X}) \ , && \qquad
    \Filt^n (\mathcal{X}) \ , && \qquad
    \Filt^\Sigma (\mathcal{X}) \ , \\
    \Grad^n_\mathbb{Q} (\mathcal{X}) \ , && \qquad
    \Grad^F (\mathcal{X}) \ , && \qquad
    \Filt^n_\mathbb{Q} (\mathcal{X}) \ , && \qquad
    \Filt^C (\mathcal{X}) \ ,
\end{alignat*}
for all integers $n \geq 0$,
all $\mathbb{Z}$-lattices~$\Lambda$,
all integral cones~$\Sigma$,
all finite-dimensional $\mathbb{Q}$-vector spaces~$F$,
and all rational cones~$C$,
can only take finitely many isomorphism types.

\begin{para}[Properties of formal lattices]
    Let~$X, Y$ be formal $\mathbb{Q}$-lattices,
    and $f \colon Y \to X$ a map.
    Recall that $\mathsf{Lat} (\mathbb{Q})$ denotes
    the category of finite-dimensional $\mathbb{Q}$-vector spaces.
    We make the following definitions:

    \begin{enumerate}
        \item
            $f$ is \emph{injective} (or~\emph{surjective})
            if for any $F \in \mathsf{Lat} (\mathbb{Q})$,
            the induced map $Y (F) \to X (F)$ is injective (or~surjective).

        \item
            $f$ is a \emph{covering family} of~$X$
            if it is surjective and
            $Y \simeq \coprod_{i \in I} F_i$,
            with $F_i \in \mathsf{Lat} (\mathbb{Q})$.

        \item
            $X$ is \emph{quasi-compact}
            if it admits a finite covering family,
            meaning that the index set~$I$ is finite;
            $f$ is \emph{quasi-compact}
            if for any map $F \to X$ with $F \in \mathsf{Lat} (\mathbb{Q})$,
            the fibre product $F \times_X Y$ is quasi-compact.

        \item
            $X$ is \emph{quasi-separated}
            if the diagonal map $X \to X \times X$
            is quasi-compact.

        \item
            $X$ is \emph{coherent}
            if it is quasi-compact and quasi-separated.

        \item
            $X$ is \emph{locally finite}
            if, for any $\alpha, \beta \in \mathsf{Face} (X)$
            with~$\alpha$ non-degenerate,
            the set $\mathrm{Hom} (\beta, \alpha)$ is finite.
            In particular, $\mathrm{Aut} (\alpha)$ is finite
            for all $\alpha \in \mathsf{Face}^\mathrm{nd} (X)$.

        \item
            Finally, recall from \cref{para-faces}
            that~$f$ is \emph{unramified}
            if it preserves non-degenerate faces.
    \end{enumerate}
    Here, the term \emph{coherent} comes from the theory of \emph{coherent topoi},
    as in \textcite[D3.3]{johnstone-2002-elephant-ii}.
    Indeed, the category of formal $\mathbb{Q}$-lattices
    is a coherent topos,
    and there is a notion of coherent objects in any coherent topos,
    which coincides with the one given above for formal $\mathbb{Q}$-lattices.

    Note that being coherent is stronger than
    being a compact object in the category of formal $\mathbb{Q}$-lattices.
    For example, the quotient of~$\mathbb{Q}^2$
    by a $\mathbb{Z}$-action of infinite order
    is a compact object, as it is a coequalizer of two maps
    $\mathbb{Q}^2 \rightrightarrows \mathbb{Q}^2$,
    but it is not coherent, as it is not quasi-separated.
\end{para}

\begin{theorem}
    \label{thm-cl-loc-finite}
    Let~$\mathcal{X}$ be a stack with quasi-compact graded points.
    Then $\mathrm{CL}_{\mathbb{Q}} (\mathcal{X})$ is locally finite.
\end{theorem}

\begin{proof}
    Let $\alpha, \beta \in \mathsf{Face} (\mathcal{X})$
    with~$\alpha$ non-degenerate,
    and write $n = \dim \alpha$.
    Choose a field~$K$ and a representable morphism
    $x \colon \mathcal{Z} = \mathrm{B} \mathbb{G}_{\smash{\mathrm{m}, K}}^n \to \mathcal{X}$
    such that the map
    $\mathrm{CL}_\mathbb{Q} (\mathcal{Z}) \simeq \mathbb{Q}^n
    \to \mathrm{CL}_\mathbb{Q} (\mathcal{X})$
    is isomorphic as a face to~$\alpha$,
    and we identify this~$\mathbb{Q}^n$ with the source of~$\alpha$.
    Let~$F$ be the source of~$\beta$,
    and let $\mathcal{Z}_\beta \subset \Grad^F (\mathcal{Z})$
    be the preimage of~$\mathcal{X}_\beta$ under the morphism
    $\Grad^F (\mathcal{Z}) \to \Grad^F (\mathcal{X})$.
    By \cref{lemma-grad-preserves-qc-morphisms},
    the morphism $\mathcal{Z}_\beta \to \mathcal{X}_\beta$ is quasi-compact;
    by assumption, the morphism $\mathcal{X}_\beta \to \mathcal{X}$ is quasi-compact,
    so that $\mathcal{Z}_\beta \to \mathcal{X}$ is quasi-compact,
    and so is the morphism $\mathcal{Z}_\beta \to \mathcal{Z}$.
    It follows that $\mathrm{Hom} (\beta, \alpha) = \uppi_0 (\mathcal{Z}_\beta)$ is finite.
\end{proof}

\begin{lemma}
    \label{lemma-lattice-coherent-gluing}
    Let~$X$ be a tame formal $\mathbb{Q}$-lattice,
    $F$ a finite-dimensional $\mathbb{Q}$-vector space,
    $\alpha_1, \alpha_2 \colon F \to X$ two non-degenerate faces,
    and $Y = \mathrm{coeq} (\alpha_1, \alpha_2)$ the coequalizer.
    Assume that~$X$ is coherent,
    that the map $X \to Y$ is unramified,
    and that~$Y$ is locally finite.
    Then~$Y$ is tame and coherent.
\end{lemma}

\begin{proof}
    First we show that $Y$ is tame. Let $\alpha\colon E\to Y$ be a non-degenerate face. By \cref{prop:characterization-tame}, it is enough to show that the map $E \to Y$ is unramified. The face $\alpha$ lifts to a non-degenerate face $E\to X$ of $X$ and, since $X$ is tame, the map $E \to X$ is unramified. Since $X \to Y$ is unramified by assumption, the composition $E \to X \to Y$ is unramified.

    To prove coherence of $Y$ it is enough to show that the natural injection
    $X \times_Y X \hookrightarrow X \times X$
    is quasi-compact, since it is the base change
    of the diagonal $Y \hookrightarrow Y \times Y$
    along the surjection
    $X \times X \twoheadrightarrow Y \times Y$.
    Since $X \times X$ is coherent,
    it is enough to show that the lattice
    $R = X \times_Y X$ is quasi-compact.

    Write $S = F \sqcup F \sqcup X$,
    and let $q_1, q_2 \colon S \to X$
    be maps given by $(\alpha_1, \alpha_2, \mathrm{id}_X)$
    and $(\alpha_2, \alpha_1, \mathrm{id}_X)$, respectively,
    which are unramified.
    For each integer $n \geq 0$, define $R_n \subset X \times X$
    as the image of the map
    $(q_1, q_2) \colon S_n = S \times_{q_2, X, q_1} \cdots
    \times_{q_2, X, q_1} S \to X \times X$,
    where $S$ appears~$n$ times,
    and the map uses~$q_1$ of the first factor and~$q_2$ of the last factor.
    Then $R_n \subset R_{n+1}$ for all~$n$,
    and $R = \bigcup_{n \geq 0} R_n$,
    since it is the equivalence relation on~$X$ generated by~$(\alpha_1, \alpha_2)$.
    Note that~$S_n$ and~$X \times X$ are coherent, so~$R_n$ is coherent.
    It is thus enough to show that $R = R_n$ for some~$n$.

    Let $Z_n \subset S_n$ be the sublattice
    which is the preimage of $\Delta \subset X^{n+1}$
    under the natural map $S_n \to X^{n+1}$,
    where $\Delta$ is the locus where the $(n + 1)$ coordinates are not all distinct.
    Let $T_n \subset S_n$ be the sublattice spanned by faces not in~$Z_n$,
    which is coherent, since for a covering family for~$S_n$,
    removing those in~$Z_n$ gives a covering family for~$T_n$.
    Then~$R_n$ is the image of the projection
    $\coprod_{0 \leq m \leq n} T_m \to X \times X$.

    For integers $m \geq n \geq 0$, there is a projection $S_m \to S_n$ to the first~$n$ factors. It is unramified by a repeated application of \cref{prop:unramified-base-change}, since the maps $q_1$ and $q_2$ are unramified by \cref{prop:characterization-tame}. The induced projection
    $T_m \to T_n$ is thus also unramified.
    Let $T_{n, m} \subset T_n$ be its image,
    so $T_n = T_{n, n} \supset T_{n, n+1} \supset \cdots$.
    This is a descending sequence of coherent lattices,
    and must stabilize. Let $T_{n, \infty} = \bigcap_{m \geq n} T_{n, m} \subset T_n$.
    Then $T_{m, \infty}$ surjects onto~$T_{n, \infty}$ for all $m \geq n$,
    as $T_{m, \infty} = T_{m, m'}$ for some~$m'$.

    Suppose that $R \neq R_n$ for all~$n$.
    Then $T_n \neq \varnothing$ for all~$n$,
    so $T_{n, m} \neq \varnothing$ for all~$m \geq n \geq 0$,
    and $T_{n, \infty} \neq \varnothing$ for all~$n$.
    Let $\beta_0 \colon E \to X = T_0$ be a non-degenerate face
    that is maximal in~$T_{0, \infty}$.
    Then $\beta_0$ lifts to a sequence $(\beta_n)_{n \geq 0}$,
    with $\beta_n \colon E \to T_{n, \infty}$
    and~$\beta_{n+1}$ mapping to~$\beta_n$ for all~$n$.
    Write $\gamma_n \colon E \to X$
    for the image of~$\beta_m$ in the $n$-th factor~$X$, where $m \geq n$
    and we label the factors by $0, \dotsc, m$.
    Then the faces~$\gamma_n$ are all distinct,
    while having the same image~$\gamma$ in~$Y$.
    Writing $\alpha \colon F \to Y$ for the image of~$\alpha_1$ and~$\alpha_2$,
    which is non-degenerate by assumption,
    we obtain infinitely many distinct morphisms $\gamma \to \alpha$
    in $\mathsf{Face} (Y)$, a contradiction.
\end{proof}

\begin{lemma}
    \label{lemma-quotient-stack-special-face}
    Let~$\mathcal{X} = U / \mathrm{GL} (n)$ be a quotient stack,
    with~$U$ an algebraic space finitely presented over a base $S$ as in \cref{assumption-stack-basic}.
    Then~$\mathcal{X}$ has only finitely many special faces,
    and $\mathrm{CL}_{\mathbb{Q}} (\mathcal{X})$ is coherent.
\end{lemma}

\begin{proof}
    Consider the projection $\mathcal{X} \to \mathrm{BGL} (n)$.
    To show that~$\mathcal{X}$ has only finitely many special faces,
    it is enough to show that there are only finitely many
    non-degenerate faces in~$\mathrm{CL}_\mathbb{Q} (\mathrm{BGL} (n))$
    that are images of special faces of~$\mathcal{X}$,
    since by \cref{lemma-grad-preserves-qc-morphisms},
    every face of~$\mathrm{CL}_\mathbb{Q} (\mathrm{BGL} (n))$
    has only finitely many preimages in~$\mathrm{CL}_\mathbb{Q} (\mathcal{X})$.

    Let $T \subset \mathrm{GL} (n)$ be a split maximal torus,
    and let~$S_T$ be the set of subtori of~$T$,
    identified with the set of linear subspaces of $\Lambda_T \otimes \mathbb{Q}$.
    By the explicit description of $\Grad^n (\mathcal{X})$
    in \cref{eg-grad-quotient-stack},
    a subtorus $T' \in S_T$ corresponds to the image
    of a special face of~$\mathcal{X}$ in~$\Lambda_T \otimes \mathbb{Q}$,
    if and only if it is of the form~$P \cap Q$,
    where~$Q$ is an intersection of hyperplanes dual to the roots of~$\mathrm{GL} (n)$,
    and~$P$ is such that it is the maximal subtorus of~$T$
    fixing the fixed locus~$U^P \subset U$.
    Therefore, it is enough to show that only finitely many elements of~$S_T$
    can act as~$P$ here.
    For a $T$-equivariant étale map $V\to U$, we have $U^P \times_U V \simeq V^P$ by \cite[Corollary~1.1.7]{halpern-leistner-instability}.
    Thus, after passing to a $T$-equivariant étale cover of~$U$
    by an affine scheme using \textcite[Corollary~10.2]{alper-hall-rydh-etale-local},
    we may assume that $U = \Spec R$ is affine.
    Since $R$ is finitely generated over $R^T$, we may embed~$U$ into a $T$-representation over~$R^T$,
    and the claim follows.
    This shows that~$\mathcal{X}$ has only finitely many special faces.

    Next, we show that $\mathrm{CL}_{\mathbb{Q}} (\mathcal{X})$ is coherent.
    By \cref{thm-special-face-closure},
    the disjoint union of all special faces
    surjects onto~$\mathrm{CL}_{\mathbb{Q}} (\mathcal{X})$,
    so the latter is quasi-compact.
    To see that it is quasi-separated,
    it is enough to show that for special faces
    $(F_1, \alpha_1)$ and $(F_2, \alpha_2)$,
    the fibre product $Y = F_1 \times_{\mathrm{CL}_{\mathbb{Q}} (\mathcal{X})} F_2$
    is quasi-compact.
    By the adjoint property of the special face closure,
    any face of~$Y$ extends to one whose image in
    $\mathrm{CL}_{\mathbb{Q}} (\mathcal{X})$
    is special, so~$Y$ is covered by preimages of special faces of~$\mathcal{X}$.
    Finally, by \cref{thm-cl-loc-finite},
    each special face of~$\mathcal{X}$ has finitely many preimages in~$Y$.
\end{proof}

\begin{theorem}
    \label{thm-cl-coherent}
    Let~$\mathcal{X}$ be a quasi-compact stack
    with quasi-compact graded points.
    Then $\mathrm{CL}_{\mathbb{Q}} (\mathcal{X})$ is coherent.
\end{theorem}

\begin{proof}
    By \textcite[Lemma~4.4.6]{halpern-leistner-instability},
    there is a representable smooth cover $\mathcal{U} \to \mathcal{X}$,
    such that~$\mathcal{U}$
    is a quasi-compact quotient stack,
    and for any $F \in \mathsf{Lat} (\mathbb{Q})$,
    the morphism $\Grad^F (\mathcal{U}) \to \Grad^F (\mathcal{X})$
    is also a representable smooth cover.
    Write $\mathcal{V} = \mathcal{U} \times_\mathcal{X} \mathcal{U}$,
    which is also a quasi-compact quotient stack.
    Since~$\Grad^F$ preserves limits,
    we have a coequalizer diagram
    $\mathrm{CL}_\mathbb{Q} (\mathcal{V}) \rightrightarrows
    \mathrm{CL}_\mathbb{Q} (\mathcal{U}) \to
    \mathrm{CL}_\mathbb{Q} (\mathcal{X})$.
    Here, we are using the fact that if
    $\mathcal S \to \mathcal R$
    is a surjective universally open morphism of algebraic stacks, then
    $\uppi_0 (\mathcal{S} \times_\mathcal{R} \mathcal{S}) \rightrightarrows
    \uppi_0 (\mathcal{S}) \to \uppi_0 (\mathcal{R})$
    is a coequalizer diagram.

    By \cref{lemma-quotient-stack-special-face},
    $\mathrm{CL}_\mathbb{Q} (\mathcal{V})$ is coherent,
    so there is a coequalizer diagram
    $\coprod_{1 \leq i \leq r} F_i \rightrightarrows
    \mathrm{CL}_\mathbb{Q} (\mathcal{U}) \to
    \mathrm{CL}_\mathbb{Q} (\mathcal{X})$,
    where $F_i \in \mathsf{Lat} (\mathbb{Q})$ for $i = 1, \dotsc, r$.
    By \cref{lemma-cl-unramified},
    we may assume that the faces $F_i \rightrightarrows \mathrm{CL}_\mathbb{Q} (\mathcal{U})$
    are non-degenerate.

    Write $X_i = \mathrm{coeq} \bigl(\coprod_{1 \leq j \leq i} F_j
    \rightrightarrows \mathrm{CL}_\mathbb{Q} (\mathcal{U}) \bigr)$
    for $0 \leq i \leq r$,
    so $X_0 = \mathrm{CL}_\mathbb{Q} (\mathcal{U})$,
    $X_r = \mathrm{CL}_\mathbb{Q} (\mathcal{X})$,
    and we have coequalizer diagrams
    $F_i \rightrightarrows X_{i-1} \to X_i$.
    Since the map~$X_0 \to X_r$ is unramified,
    so are all the maps $X_{i-1} \to X_i$,
    and the faces $F_i \rightrightarrows X_{i-1}$ are non-degenerate.
    By \cref{thm-cl-loc-finite},
    $X_r$ is locally finite, and so are all the~$X_i$.
    It now follows from \cref{lemma-lattice-coherent-gluing}
    that all the~$X_i$ are coherent.
\end{proof}

\begin{para}[Proof of \texorpdfstring{\cref{thm-finiteness}}{Theorem \ref{thm-finiteness}}]
    \label{proof-finiteness}
    By \cref{thm-cl-loc-finite},
    it is enough to show that~$\mathcal{X}$
    has only finitely many special faces and special cones.
    By \cref{thm-cotangent-arrangement},
    it is enough to show that~$\mathcal{X}$
    has only finitely many special faces.

    By \textcite[Proposition~8.2]{hall-rydh-2019-tannaka},
    we may choose a stratification of~$\mathcal{X}$
    by locally closed substacks~$\mathcal{U}_a = U_a / \mathrm{GL} (n)$,
    where each~$U_a$ is a quasi-affine scheme.
    Write $\mathcal{U} = \coprod_a \mathcal{U}_a$,
    with the stratification morphism $\mathcal{U} \to \mathcal{X}$,
    and write $X = \mathrm{CL}_{\mathbb{Q}} (\mathcal{X})$ for short.

    Let~$(F, \alpha)$ be a special face of~$\mathcal{X}$.
    Then $\Grad^F (\mathcal{U}) \to \Grad^F (\mathcal{X})$
    is also a locally closed stratification,
    and~$\mathcal{X}_\alpha$ is stratified by stacks~$\mathcal{U}_{\alpha_i}$
    corresponding to the preimages
    $\alpha_i$ of~$\alpha$ in $\mathrm{CL}_{\mathbb{Q}} (\mathcal{U})$.
    Let~$(F_i, \tilde{\alpha}_i)$ be the special face closure of~$\alpha_i$.
    The faces~$\tilde{\alpha}_i$ are pairwise distinct,
    since the images of~$\mathcal{U}_{\tilde{\alpha}_i}$
    in~$\mathcal{X}_\alpha$ are pairwise disjoint.

    It is enough to show that for any given collection
    $(\tilde{\alpha}_i)_{i \in I}$ of special faces of~$\mathcal{U}$,
    there are only finitely many special faces~$\alpha$ of~$\mathcal{X}$
    that can give rise to this collection,
    since by \cref{lemma-quotient-stack-special-face},
    there are only finitely many such collections.

    Let~$(F_i, \alpha'_i)$ be the image of~$\tilde{\alpha}_i$ in~$X$.
    The face~$\alpha$ lifts to a face
    $\bar{\alpha} \colon F \to Y = F_1 \times_X \cdots \times_X F_m$.
    Let $\bar{\alpha} \to \bar{\beta}$ be a morphism in~$\mathsf{Face} (Y)$,
    which induces subfaces $\beta_i \subset \tilde{\alpha}_i$ containing~$\alpha_i$,
    and let~$\beta$ be the image of~$\bar{\beta}$ in~$X$.
    Then the stacks $\mathcal{U}_{\beta_i} \simeq \mathcal{U}_{\alpha_i}$
    stratify~$\mathcal{X}_\beta$,
    so the morphism $\mathcal{X}_\beta \to \mathcal{X}_\alpha$
    is universally bijective.
    If it is also étale, 
    then it will be an isomorphism,
    and we will have $\alpha = \beta$
    since~$\alpha$ is special.
    By \cref{thm-tangent-filt-sigma},
    this happens when~$\alpha$ does not lie in any hyperplane
    in the cotangent arrangement $\Phi^\mathrm{cot} (\mathcal{X}, \beta)$
    defined in \cref{para-cotangent-arrangement},
    in which case the above argument gives $\alpha = \beta$.
    We conclude that either $\alpha = \beta$,
    or~$\alpha$ is contained in a hyperplane in
    $\Phi^\mathrm{cot} (\mathcal{X}, \beta)$.
    Now, $Y$ is coherent since~$X$ is,
    so it has only finitely many maximal non-degenerate faces.
    Then~$\bar{\alpha}$ must be obtained by choosing a maximal non-degenerate face of~$Y$,
    then repeatedly taking hyperplanes in the face
    belonging to the cotangent arrangement.
    This shows that there are only finitely many special faces~$\alpha$
    that give rise to a given collection~$(\tilde{\alpha}_i)_{i \in I}$,
    completing the proof.
    \qed
\end{para}

\begin{para}[A derived version]
    Finally, we remark that a derived version of
    \cref{thm-finiteness} is also true.
    More precisely, for a derived stack~$\mathcal{X}$
    as in \cref{para-constancy-derived},
    such that its classical truncation~$\mathcal{X}_\mathrm{cl}$
    is quasi-compact and has quasi-compact graded points,
    then~$\mathcal{X}$ has only finitely many special faces and cones in the derived sense
    defined in \cref{para-constancy-derived},
    and there are only finitely many morphisms between them.

    Indeed, this follows from the classical version as follows.
    Suppose that~$\alpha$ is a special face of~$\mathcal{X}$
    in the derived sense,
    and let~$\alpha'$ be the special face closure of~$\alpha$
    as a special face of~$\mathcal{X}_\mathrm{cl}$.
    Consider the derived cotangent arrangement~$\Phi^\mathrm{cot} (\mathcal{X}, \alpha')$
    as defined in~\cref{para-constancy-derived},
    which is finite since~$\mathcal{X}_{\alpha'}$ is quasi-compact.
    Then~$\alpha$ must be an intersection of
    hyperplanes in~$\Phi^\mathrm{cot} (\mathcal{X}, \alpha')$,
    since if $\alpha'' \subset \alpha'$
    is the intersection of all hyperplanes in~$\Phi^\mathrm{cot} (\mathcal{X}, \alpha')$
    containing~$\alpha$,
    then the induced morphism
    $\dmathcal{X}_{\alpha''} \to \dmathcal{X}_\alpha$
    is classically an isomorphism and is étale,
    so it is an isomorphism of derived stacks, and $\alpha'' = \alpha$.
\end{para}

\subsection{Hall categories and associativity}

\begin{para}
    In this section, we prove the \emph{associativity theorem},
    \cref{thm-associativity},
    which states that for a stack~$\mathcal{X}$
    with finite cotangent weights as in
    \cref{para-finite-cotangent-weights},
    we have a canonical isomorphism
    \begin{equation*}
        \mathcal{X}_{\smash{\sigma_1 \uparrow \sigma_2}}^+ \longrightarrow
        \mathcal{X}_{\sigma_1}^+ \underset{\mathcal{X}_{\alpha_1}}{\times}
        \mathcal{X}_{\sigma_2}^+ \ ,
    \end{equation*}
    where $\sigma_1 \subset \alpha_1 \subset \sigma_2$
    are special cones in~$\mathcal{X}$,
    with $\alpha_1$ a special face,
    and $\sigma_1 \uparrow \sigma_2$ denotes another special cone
    obtained from~$\sigma_1$ and~$\sigma_2$.

    It is perhaps helpful
    to view $\mathcal{X}_{\sigma_2}^+$ as
    $(\mathcal{X}_{\alpha_1})_{\tilde{\sigma}_2}^+$,
    where $\tilde{\sigma}_2$ is the canonical lift
    of~$\sigma_2$ to a cone in $\mathrm{CL}_\mathbb{Q} (\mathcal{X}_{\alpha_1})$,
    as in \cref{para-cone-filt-cone-filt}.
    The morphism $(\mathcal{X}_{\alpha_1})_{\tilde{\sigma}_2}^+ \to \mathcal{X}_{\alpha_1}$
    is then the morphism~$\mathrm{ev}_{1, \tilde{\sigma}_2}$ for~$\mathcal{X}_{\alpha_1}$,
    in the notation of \cref{para-x-sigma-plus}.
    See also \cref{para-intro-associativity}
    for the rough ideas.

    This is called the \emph{associativity theorem},
    because when~$\mathcal{X}$ is a moduli stack of objects
    in an abelian category,
    this property is crucial in proving associativity
    for various types of Hall algebras obtained from~$\mathcal{X}$.
    See \cref{para-intro-associativity} for the main ideas,
    and \cref{subsec-lms} below for details.

    This result can be interpreted as a functorial property,
    where we define a \emph{Hall category} associated to each stack,
    and construct functors from the Hall category.
    Such functors can be seen as analogues of Hall algebras
    for general algebraic stacks.

    See also \textcite[Proposition~5.22]{kinjo-park-safronov-coha}
    for a related result.
\end{para}

\begin{para}[The extension operation]
    \label{para-extension-operation}
    Let~$F$ be a finite-dimensional $\mathbb{Q}$-vector space,
    and~$\Psi$ a cone arrangement in~$F$,
    which we identify with a set of cones in~$F$.
    Assume that~$\Psi$ does not contain
    infinite strictly descending chains of cones.

    We define a binary operation~$\uparrow$
    on the set~$\Psi$ as follows.
    For cones $C, C' \in \Psi$,
    define $C \uparrow C'$
    to be the minimal cone in~$\Psi$
    such that for any~$v \in C^\circ$ and~$v' \in C'$,
    where $C^\circ \subset C$ is the interior,
    there exists~$\varepsilon > 0$ such that~$v + \varepsilon v' \in C \uparrow C'$.

    The operation~$\uparrow$ is often associative.
    For example, this is the case when~$\Psi$
    comes from a hyperplane arrangement as in
    \cref{para-face-arrangement-examples}~\cref{item-eg-hyperplane-arrangement}. When restricted to chambers, the resulting structure is known as the \emph{Tits monoid} of the hyperplane arrangement \cite[\S1.4]{aguiar-topics-hyperplane}. Associativity also holds in the case of composable cones in the Hall category below.
\end{para}

\begin{para}[The Hall category]
    \label{para-hall-category}
    We now globalize the operation~$\uparrow$
    to define the \emph{Hall category} of a stack,
    whose composition is given by this operation.

    For a stack~$\mathcal{X}$ with finite cotangent weights,
    as in \cref{para-finite-cotangent-weights},
    define the \emph{extended Hall category} of~$\mathcal{X}$,
    denoted by $\mathsf{Hall}^+ (\mathcal{X})$, as follows:

    \begin{itemize}
        \item
            An object of~$\mathsf{Hall}^+ (\mathcal{X})$ is an object
            $(F, \alpha) \in \mathsf{Face}^\mathrm{sp} (\mathcal{X})$.

        \item
            A morphism $(F, \alpha) \to (F', \alpha')$
            in~$\mathsf{Hall}^+ (\mathcal{X})$ is a pair~$(f, C)$,
            where $f \colon F \to F'$ is a map of faces,
            and $C \subset F'$ is a cone of full dimension
            such that $f (F) \subset C$
            and $(C, \sigma) \in \mathsf{Cone}^\mathrm{sp} (\mathcal{X})$,
            where $\sigma = \alpha' |_{C}$.
            We often denote such a morphism as
            $\alpha \overset{\sigma}{\to} \alpha'$.

        \item
            The composition of
            $(f, C) \colon (F, \alpha) \to (F', \alpha')$
            and
            $(f', C') \colon (F', \alpha') \to (F'', \alpha'')$,
            is the morphism
            $(f' \circ f, C \uparrow C')$,
            where~$C \uparrow C'$
            is taken in~$F''$.
            We often denote this by $\sigma \uparrow \sigma'$,
            where $\sigma = \alpha' |_{C}$ and $\sigma' = \alpha'' |_{C'}$.
            One can verify using
            \cref{thm-cotangent-arrangement}~\cref{item-cotangent-arrangement-2}
            that composition is associative.

        \item
            The identity morphism of~$(F, \alpha)$
            is the pair~$(\mathrm{id}_F, F)$.
    \end{itemize}
    Define the \emph{Hall category} of~$\mathcal{X}$
    as the subcategory
    $\mathsf{Hall} (\mathcal{X}) \subset \mathsf{Hall}^+ (\mathcal{X})$
    consisting of the same objects,
    but only those morphisms
    $(f, C) \colon (F, \alpha) \to (F', \alpha')$
    satisfying the following condition:

    \begin{itemize}
        \item
            There exists $\lambda \in F'$,
            such that~$C$ is the minimal cone in~$F'$
            that is a special cone and contains both~$f (F)$ and~$\lambda$.
    \end{itemize}
    See \cref{para-hall-category-chambers}
    for an alternative, more explicit description.

    We call these \emph{Hall categories}
    because they are related to the construction of
    various versions of \emph{Hall algebras}
    for moduli stacks of objects in abelian categories,
    which we explain in \cref{subsec-lms} below.
    We will consider functors from the Hall category
    as an analogue of Hall algebras for general algebraic stacks,
    which we explain in \cref{para-gen-hall-alg}.

    The Hall category $\mathsf{Hall} (\mathcal{X})$ generalizes the \emph{category of lunes} of a hyperplane arrangement, as in \textcite[Chapter~4]{aguiar-topics-hyperplane}.
.
\end{para}

\begin{para}[Hall category of~\texorpdfstring{$\mathcal{X}_\alpha$}{X\_α}]
    For a special face $\alpha$ of ~$\mathcal{X}$, we have equivalences
    \[\mathsf{Hall}^+(\mathcal X_\alpha)=(\alpha\downarrow \mathsf{Hall}^+(\mathcal{X}))\quad \mathsf{Hall}(\mathcal X_\alpha)=(\alpha\downarrow \mathsf{Hall}(\mathcal{X}))\ .\]
    This follows directly from \cref{para-cone-filt-cone-filt}
\end{para}

\begin{para}[An alternative description]
    \label{para-hall-category-chambers}
    The category $\mathsf{Hall} (\mathcal{X})$
    can be alternatively described as follows,
    using \cref{thm-cotangent-arrangement}
    and the discussion in \cref{para-special-closure-rays}:

    \begin{itemize}
        \item
            A morphism $\alpha \overset{\sigma}{\to} \alpha'$
            is a morphism of faces $\alpha \to \alpha'$ together with
            a chamber $\sigma \subset \alpha'$
            in the hyperplane arrangement
            $\Phi^\alpha_{\smash{\alpha'}} \subset \Phi^\mathrm{cot} (\mathcal{X}, \alpha')$
            consisting of hyperplanes containing the image of~$\alpha$.

        \item
            For morphisms
            $\alpha \overset{\sigma}{\to} \alpha' \overset{\sigma'}{\to} \alpha''$
            as above,
            the chamber $\sigma \uparrow \sigma' \subset \alpha''$
            is the unique chamber in the hyperplane arrangement
            $\Phi^\alpha_{\smash{\alpha''}}$ defined as above satisfying
            $\sigma\subset \sigma \uparrow \sigma'\subset \sigma'$.
    \end{itemize}
    Here, it is the special closures of the chambers~$\sigma$, etc.,
    that correspond to the cones~$\sigma$, etc., used in \cref{para-hall-category}.
    Although this might cause some ambiguity in notation,
    we try to clarify it whenever possible,
    and in any case, the notation~$\mathcal{X}_\sigma^+$ is unambiguous
    since $\mathcal{X}_\sigma^+ \simeq \mathcal{X}_{\smash{\sigma^\mathrm{sp}}}^+$.
\end{para}

\begin{theorem}
    \label{thm-associativity}
    Let~$\mathcal{X}$ be a stack
    with finite cotangent weights.
    Then for any chain of morphisms
    \begin{equation*}
        \alpha_0 \overset{\sigma_1}{\longrightarrow}
        \alpha_1 \overset{\sigma_2}{\longrightarrow}
        \cdots \overset{\sigma_n}{\longrightarrow} \alpha_n
    \end{equation*}
    in the Hall category~$\mathsf{Hall}^+ (\mathcal{X})$
    defined in \cref{para-hall-category},
    we have a canonical isomorphism
    \begin{equation*}
        \mathcal{X}_{\smash{\sigma_1 \uparrow \cdots \uparrow \sigma_n}}^+
        \longsimto
        \mathcal{X}_{\sigma_1}^+ \underset{\mathcal{X}_{\alpha_1}}{\times}
        \cdots \underset{\mathcal{X}_{\alpha_{n - 1}}}{\times}
        \mathcal{X}_{\sigma_n}^+ \ ,
    \end{equation*}
    where $\sigma_1 \uparrow \cdots \uparrow \sigma_n$
    denotes the composition in the Hall category.
    The fibre product is defined by the morphisms
    $\mathrm{gr}_{\sigma_i} \colon \mathcal{X}_{\sigma_i}^+ \to \mathcal{X}_{\alpha_i}$
    and morphisms
    $\mathcal{X}_{\sigma_i}^+ \to \mathcal{X}_{\alpha_{i-1}}$
    induced by $\alpha_{i-1} \hookrightarrow \sigma_i$.
\end{theorem}

\begin{proof}
    We use the following consequence of
    \textcite[Proposition~3.4.4]{halpern-leistner-instability}:
    \begin{itemize}
        \item
            Let $\alpha \colon \mathbb{Q}^n \to \mathrm{CL}_\mathbb{Q} (\mathcal{X})$
            be a face, and let $(e_1, \dotsc, e_n)$
            be the standard basis of~$\mathbb{Q}^n$.
            If no hyperplanes in the cotangent arrangement
            $\Phi^\mathrm{cot} (\mathcal{X}, \alpha)$
            defined in~\cref{para-cotangent-arrangement}
            intersect with~$(\mathbb{Q}_{> 0})^n \subset F$,
            then the natural morphism
            \begin{equation*}
                \mathcal{X}_{\sigma}^+ \longrightarrow
                \mathcal{X}_{\sigma_1}^+
                \underset{\mathcal{X}_{\alpha_1}}{\times}
                \mathcal{X}_{\sigma_2}^+
                \underset{\mathcal{X}_{\alpha_2}}{\times}
                \cdots
                \underset{\mathcal{X}_{\alpha_{n - 1}}}{\times}
                \mathcal{X}_{\sigma_n}^+
            \end{equation*}
            is universally bijective, where
            $\sigma = (\mathbb{Q}_{\geq 0})^n \subset F$,
            and $\sigma_i \subset F$ is the cone generated by
            $\pm e_1, \dotsc, \pm e_{i-1}, e_i$,
            and~$\alpha_i$ is the face spanned by~$e_1, \dotsc, e_i$.
    \end{itemize}
    This can be deduced from the cited proposition by induction on~$n$.

    To prove the theorem,
    it is enough to consider the case when $n = 2$,
    since the general case follows by induction.
    Let~$F$ be the underlying vector space of~$\alpha_2$,
    and choose a set of generating vectors
    $v_1, \dotsc, v_{m'} \in F$ for~$\sigma_1 \uparrow \sigma_2$,
    such that~$v_1, \dotsc, v_m$ generates~$\sigma_1$ for some~$m \leq m'$.
    Then, applying the above result to the possibly degenerate face
    $\mathbb{Q}^{\smash{m'}} \to F \to \mathrm{CL}_\mathbb{Q} (\mathcal{X})$,
    where the first map sends the standard basis to~$v_1, \dotsc, v_{m'}$,
    we see that the natural morphism
    $\mathcal{X}_{\smash{\sigma_1 \uparrow \sigma_2}}^+ \to
    \mathcal{X}_{\sigma_1}^+ \times_{\mathcal{X}_{\alpha_1}}
    \mathcal{X}_{\sigma_2}^+$
    is universally bijective.
    Therefore, it is enough to show that it is also étale.

    We use an argument similar to the proof of
    \cref{lemma-cotangent-arrangement-cones}.
    By \cref{lemma-theta-retract},
    it suffices to show that the morphism is étale at split filtrations,
    since any open substack of~$\mathcal{X}_{\smash{\sigma_1 \uparrow \sigma_2}}^+$
    containing the split filtrations is the whole stack.
    By \cref{thm-tangent-filt-sigma}, we have
    $\mathbb{L}_{\smash{\dmathcal{X}_{\smash{\sigma_1 \uparrow \sigma_2}}^+}}
    |_{\smash{\mathcal{X}_{\alpha_2}}}
    \simeq (\mathbb{L}_\mathcal{X}
    |_{\smash{\mathcal{X}_{\alpha_2}}})
    _{\sigma_1 \uparrow \sigma_2, -}$,
    with notations as in the theorem,
    while the cotangent complex of the fibre product
    $\dmathcal{X}_{\sigma_1}^+ \times_{\smash{\dmathcal{X}_{\alpha \mathrlap{_1}}}}
    \dmathcal{X}_{\sigma_2}^+$
    at split filtrations is the cofibre of the morphism
    $(\mathbb{L}_\mathcal{X} |_{\smash{\mathcal{X}_{\alpha_1}}})_0
    \to (\mathbb{L}_\mathcal{X} |_{\smash{\mathcal{X}_{\alpha_1}}})_{\sigma_1, -}
    \oplus (\mathbb{L}_\mathcal{X} |_{\smash{\mathcal{X}_{\alpha_2}}})_{\sigma_2, -}$,
    with each term pulled back to the fibre product.
    This cofibre contains weights
    in~$\mathbb{L}_\mathcal{X} |_{\smash{\mathcal{X}_{\alpha_2}}}$
    that either pair non-positively with~$\sigma_2$,
    or have non-positive but not constantly zero pairing with~$\sigma_1$.
    This set of weights contains the set of weights
    that pair non-positively with~$\sigma_1 \uparrow \sigma_2$,
    and the two sets of weights coincide
    in cohomological degrees~$0, 1$.
    Therefore, the morphism
    $\dmathcal{X}_{\smash{\sigma_1 \uparrow \sigma_2}}^+ \to
    \dmathcal{X}_{\sigma_1}^+ \times_{\smash{\dmathcal{X}_{\alpha \mathrlap{_1}}}}
    \dmathcal{X}_{\sigma_2}^+$
    induces a morphism of cotangent complexes
    $\mathrm{H}^i (\mathbb{L}_{(-)})$
    which is an isomorphism for $i = 0, 1$,
    and a surjection for $i = -1$.
    Its relative cotangent complex thus has $\mathrm{H}^i = 0$
    for $i = -1, 0, 1$,
    and so does the relative cotangent complex of the morphism
    of classical truncations
    $\mathcal{X}_{\smash{\sigma_1 \uparrow \sigma_2}}^+ \to
    \mathcal{X}_{\sigma_1}^+ \times_{\mathcal{X}_{\alpha_1}}
    \mathcal{X}_{\sigma_2}^+$,
    which is hence étale.
\end{proof}

\begin{para}[The category of spans]
    \label{para-cat-of-spans}
    We now explain how the associativity theorem,
    \cref{thm-associativity},
    can be reformulated as a functorial property.

    Consider the \emph{category of spans}
    $\mathsf{Span} (\mathsf{St}_S)$
    described as follows:

    \begin{itemize}
        \item
            An object of~$\mathsf{Span} (\mathsf{St}_S)$
            is an object~$\mathcal{Y} \in \mathsf{St}_S$,
            as in \cref{assumption-stack-basic}.

        \item
            A morphism from~$\mathcal{Y}$ to~$\mathcal{Y}'$
            is a span
            $\mathcal{Y} \leftarrow \mathcal{Z} \to \mathcal{Y}'$
            in~$\mathsf{St}_S$.

        \item
            Composition of spans is defined by the fibre product
            \begin{equation*}
                \begin{tikzcd}[sep={2.5em,between origins}]
                    && \mathcal{W} \ar[dl] \ar[dr]
                    \ar[dd, phantom, "\ulcorner" {pos=.2, shift={(-.1em, 0)}, rotate=-45}]
                    && \\
                    & \mathcal{Z} \ar[dl] \ar[dr]
                    && \mathcal{Z}' \ar[dl] \ar[dr]
                    \\
                    \mathcal{Y} && \mathcal{Y}' && \mathcal{Y}'' \rlap{ ,}
                \end{tikzcd}
            \end{equation*}
            where the outer span $\mathcal{Y} \leftarrow \mathcal{W} \to \mathcal{Y}''$
            is the composition.

        \item
            The identity morphism of~$\mathcal{Y}$
            is the span $\mathcal{Y} \leftarrow \mathcal{Y} \to \mathcal{Y}$,
            given by identity morphisms.
    \end{itemize}
    Strictly speaking, since $\mathsf{St}_S$ is a $2$-category,
    $\mathsf{Span} (\mathsf{St}_S)$ is a $3$-category,
    and a precise definition can be found in
    \textcite[\S6.1]{liu-zheng-six-operations}.
    However, for our purposes, it is often enough
    to consider its truncation as a $1$-category.

    The category of spans classifies
    \emph{three-functor formalisms} in the sense of
    \textcite{scholze-six-functor};
    see \cref{para-three-functor} for more discussions.

    We can now restate
    \cref{thm-associativity}
    as the following result.
\end{para}

\begin{theorem}
    \label{thm-functorial-assoc}
    Let~$\mathcal{X}$ be a stack
    with finite cotangent weights,
    and let~$\mathsf{Hall}^+ (\mathcal{X})$ be the Hall category
    defined in \cref{para-hall-category}.
    Then there is a functor
    \begin{equation*}
        \mathcal{X}_{(-)} \colon
        \mathsf{Hall}^+ (\mathcal{X})^\mathrm{op} \longrightarrow
        \mathsf{Span} (\mathsf{St}_S)
    \end{equation*}
    sending an object~$\alpha$ to the stack~$\mathcal{X}_\alpha$,
    and a morphism~$\alpha \overset{\sigma}{\to} \alpha'$
    to the span
    $\mathcal{X}_{\alpha'} \leftarrow \mathcal{X}_\sigma^+ \to \mathcal{X}_\alpha$,
    where the left morphism is $\mathrm{gr}_\sigma$,
    and the right morphism is induced by $\alpha \hookrightarrow \sigma$.
    \qed
\end{theorem}

\section{Applications}

\subsection{Linear moduli stacks}
\label{subsec-lms}

\begin{para}
    As an application of the formalism developed in this paper,
    we formulate a notion of \emph{linear moduli stacks},
    which are stacks that behave like
    moduli stacks of objects in abelian categories,
    such as moduli stacks of quiver representations or
    coherent sheaves on a projective scheme.
    We have already informally discussed them in
    \cref{eg-component-lattice-linear-moduli,eg-special-face-lms}.

    We hope that our axiomatization provides a convenient uniform framework
    to work with various constructions of such moduli stacks,
    such as those of
    \textcite{artin-zhang-2001-hilbert,toen-vaquie-2007-moduli},
    as we discuss in \cref{eg-lms}.
    Our approach is inspired by
    \textcite[Definition~8.1]{kinjo-park-safronov-coha}.

    Using the constancy and associativity theorems,
    we show that our axioms imply many expected properties
    of a moduli stack of objects in an abelian category.
    We also explain the construction of various types of Hall algebras
    using our framework.
\end{para}

\begin{para}[Linear moduli stacks]
    \label{para-linear-moduli-stacks}
    Let~$K$ be an algebraically closed field.
    Define a \emph{linear moduli stack} over~$K$
    to be a stack~$\mathcal{X}$ as in \cref{assumption-stack-basic},
    where $S = \operatorname{Spec} K$,
    together with the following additional structures:

    \begin{itemize}
        \item
            A commutative monoid structure
            $\oplus \colon \mathcal{X} \times_S \mathcal{X} \to \mathcal{X}$,
            with unit~$0 \colon S \hookrightarrow \mathcal{X}$
            an open and closed immersion.

        \item
            A $\mathrm{B} \mathbb{G}_\mathrm{m}$-action
            $\odot \colon \mathrm{B} \mathbb{G}_\mathrm{m} \times \mathcal{X} \to \mathcal{X}$
            compatible with the monoid structure.
    \end{itemize}
    Note that these structures come with extra coherence data.
    We require the following additional property:

    \begin{itemize}
        \item
            There is an isomorphism
            \begin{equation}
                \label{eq-lms-grad}
                \coprod_{\gamma \colon \mathbb{Z} \to \uppi_0 (\mathcal{X})} {}
                \prod_{n \in \mathrm{supp} (\gamma)}
                \mathcal{X}_{\gamma (n)}
                \longsimto \mathrm{Grad} (\mathcal{X}) \ ,
            \end{equation}
            where~$\gamma$ runs through maps of sets
            $\mathbb{Z} \to \uppi_0 (\mathcal{X})$ such that
            $\mathrm{supp} (\gamma) = \mathbb{Z} \setminus \gamma^{-1} (0)$
            is finite,
            the product is over~$S$,
            and the morphism is defined by the composition
            \begin{equation*}
                \mathrm{B} \Gm \times
                \prod_{n \in \mathrm{supp} (\gamma)} \mathcal{X}_{\gamma (n)}
                \overset{(-)^n}{\longrightarrow}
                \prod_{n \in \mathrm{supp} (\gamma)} {}
                (\mathrm{B} \Gm \times \mathcal{X}_{\gamma (n)})
                \overset{\odot}{\longrightarrow}
                \prod_{n \in \mathrm{supp} (\gamma)} \mathcal{X}_{\gamma (n)}
                \overset{\oplus}{\longrightarrow}
                \mathcal{X}
            \end{equation*}
            on the component corresponding to~$\gamma$,
            where the first morphism is given by the
            $n$-th power map $(-)^n \colon \mathrm{B} \Gm \to \mathrm{B} \Gm$
            on the factor corresponding to~$\mathcal{X}_{\gamma (n)}$.
    \end{itemize}
    One could think of~\cref{eq-lms-grad}
    roughly as an isomorphism $\mathrm{Grad} (\mathcal{X}) \simeq \mathcal{X}^\mathbb{Z}$,
    where we only consider components of~$\mathcal{X}^\mathbb{Z}$
    involving finitely many non-zero classes in~$\uppi_0 (\mathcal{X})$.

    One can deduce from this axiom that
    for a $\mathbb{Z}$-lattice~$\Lambda$
    and a finite-dimensional $\mathbb{Q}$-vector space~$F$,
    we have similar isomorphisms
    \begin{align}
        \label{eq-lms-grad-n}
        \coprod_{\gamma \colon \Lambda^\vee \to \uppi_0 (\mathcal{X})} {}
        \prod_{\lambda \in \mathrm{supp} (\gamma)}
        \mathcal{X}_{\gamma (\lambda)}
        & \longsimto \Grad^\Lambda (\mathcal{X}) \ ,
        \\
        \label{eq-lms-grad-q-n}
        \coprod_{\gamma \colon F^\vee \to \uppi_0 (\mathcal{X})} {}
        \prod_{\lambda \in \mathrm{supp} (\gamma)}
        \mathcal{X}_{\gamma (\lambda)}
        & \longsimto \Grad^F (\mathcal{X}) \ ,
    \end{align}
    where~$\gamma$ is assumed of finite support in both cases.
\end{para}

\begin{para}[Examples]
    \label[example]{eg-lms}
    We list here some examples of linear moduli stacks over~$K$.

    \begin{enumerate}
        \item
            \label{item-eg-lms-abelian}
            \emph{Moduli of objects in abelian categories}.
            Let~$\mathcal{A}$ be a $K$-linear abelian category
            which is locally noetherian and cocomplete,
            in the sense of
            \textcite[\S7]{alper-halpern-leistner-heinloth-2023-moduli}.
            Consider the moduli stack~$\mathcal{M}_\mathcal{A}$
            of finitely presented objects in~$\mathcal{A}$
            in the sense of
            \textcite{artin-zhang-2001-hilbert},
            defined by the moduli functor
            \begin{equation*}
                \mathcal{M}_\mathcal{A} (R) = \Bigl(
                    \text{flat, finitely presented objects in }
                    \mathcal{A}_R = \bigl\{
                        (E, \xi_E) \bigm|
                        E \in \mathcal{A}, \ 
                        \xi_E \colon R \to \mathrm{End} (E)
                    \bigr\}
                \Bigr)
            \end{equation*}
            for commutative $K$-algebras~$R$. See
            \cite[\S7]{alper-halpern-leistner-heinloth-2023-moduli}
            for details.
            Then, if~$\mathcal{M}_\mathcal{A}$ is an algebraic stack
            satisfying the assumptions in \cref{assumption-stack-basic},
            then it is a linear moduli stack over~$K$,
            which can be deduced from
            \cite[Proposition~7.12 and Lemma~7.20]{alper-halpern-leistner-heinloth-2023-moduli}.
            Moreover, any open substack of~$\mathcal{M}_\mathcal{A}$
            closed under direct sums and direct summands
            is a linear moduli stack over~$K$.
            For instance, the moduli stack of coherent sheaves
            on a projective $K$-scheme is a linear moduli stack over~$K$.

        \item
            \label{item-eg-lms-dg}
            \emph{Moduli of objects in dg-categories}.
            Let~$\mathcal{C}$ be a $K$-linear dg-category
            of finite type, in the sense of
            \textcite[Definition~2.4]{toen-vaquie-2007-moduli}.
            Consider the moduli stack~$\mathcal{M}_\mathcal{C}$
            of right proper objects in~$\mathcal{C}$,
            as in \cite{toen-vaquie-2007-moduli},
            which is a derived stack over~$K$ defined by the moduli functor
            \begin{equation*}
                \mathcal{M}_\mathcal{C} (R) =
                \mathsf{Map} (\mathcal{C}^\mathrm{op}, \mathsf{Perf} (R))
            \end{equation*}
            for simplicial commutative $K$-algebras~$R$,
            where $\mathsf{Map} (-, -)$ denotes the mapping space of dg-categories.
            See \cite{toen-vaquie-2007-moduli} for details.
            By the argument in
            \textcite[Proposition~8.29]{kinjo-park-safronov-coha},
            if we are given an open substack
            $\mathcal{M} \subset \mathcal{M}_\mathcal{C}$,
            closed under direct sums and direct summands,
            such that it contains the zero object as an open and closed substack,
            and its classical truncation~$\mathcal{M}_\mathrm{cl}$
            is a $1$-stack satisfying the assumptions in \cref{assumption-stack-basic},
            then~$\mathcal{M}_\mathrm{cl}$ is a linear moduli stack over~$K$.

        \item
            \label{item-eg-lms-rep}
            \emph{Moduli of representations of algebras}.
            Let~$A$ be a finitely generated associative $K$-algebra.
            Consider the moduli stack~$\mathcal{M}_A$
            of representations of~$A$,
            defined by the moduli functor
            \begin{equation*}
                \mathcal{M}_A (R) = \bigl(
                    \text{left $(A \otimes_K R)$-modules,
                    flat and finitely presented over~$R$}
                \bigr)
            \end{equation*}
            for commutative $K$-algebras~$R$.
            In fact, choosing a set of generators for~$A$,
            one can deduce that each connected component of~$\mathcal{M}_A$
            admits an affine morphism of finite presentation to
            $\Spec (K) / \mathrm{GL} (n)$ for some~$n$,
            so that~$\mathcal{M}_A$ is an algebraic stack
            locally of finite presentation over~$K$, with affine diagonal.
            One can deduce as in
            \textcite[Proposition~7.12]{alper-halpern-leistner-heinloth-2023-moduli}
            that~$\mathcal{M}_A$ is a linear moduli stack over~$K$.
            In particular, the moduli stack of representations
            of a (finite) quiver, possibly with relations, over~$K$,
            is a linear moduli stack over~$K$.
    \end{enumerate}
\end{para}

\begin{para}[Special faces]
    \label{para-lms-special-faces}
    Let~$\mathcal{X}$ be a linear moduli stack.
    Recall the explicit description of the component lattice of~$\mathcal{X}$
    in \cref{eg-component-lattice-linear-moduli}.

    For classes
    $\gamma_1, \ldots, \gamma_n \in \uppi_0 (\mathcal{X})$,
    consider the face
    \begin{equation*}
        \alpha (\gamma_1, \dotsc, \gamma_n) \in \mathsf{Face} (\mathcal{X})
    \end{equation*}
    corresponding to the map
    $\mathbb{Q}^n \to \uppi_0 (\mathcal{X})$
    sending the $i$-th standard basis vector to~$\gamma_i$,
    and all other elements to~$0$,
    as in \cref{eg-special-face-lms}.

    It follows from \cref{eq-lms-grad-q-n} that
    all special faces of~$\mathcal{X}$ are of this form,
    with $\gamma_1, \dotsc, \gamma_n$ non-zero.
    In particular, by
    \cref{lemma-cot-arr-special-faces},
    every linear moduli stack
    has finite cotangent weights.

    Conversely, such a face $\alpha (\gamma_1, \dotsc, \gamma_n)$
    is special if and only if all~$\gamma_i$ belong to the set
    $M \subset \uppi_0 (\mathcal{X}) \setminus \{ 0 \}$
    consisting of classes~$\gamma$ such that
    for any two non-zero classes
    $\gamma', \gamma'' \in \uppi_0 (\mathcal{X}) \setminus \{ 0 \}$
    with $\gamma = \gamma' + \gamma''$,
    the direct sum morphism
    $\oplus \colon \mathcal{X}_{\gamma'} \times \mathcal{X}_{\gamma''}
    \to \mathcal{X}_\gamma$
    is not an isomorphism.
\end{para}

\begin{para}[Special cones]
    \label{para-lms-special-cones}
    For classes
    $\gamma_1, \ldots, \gamma_n \in \uppi_0 (\mathcal{X})$,
    we also consider the cone
    \begin{equation*}
        \sigma (\gamma_1, \dotsc, \gamma_n)
        \subset \alpha (\gamma_1, \dotsc, \gamma_n)
    \end{equation*}
    defined by the inequalities $x_1 \geq \cdots \geq x_n$,
    where the~$x_i$ are the standard coordinates on~$\mathbb{Q}^n$.

    More generally, for any partial order~$\preceq$
    on the set $\{ 1, \dotsc, n \}$,
    we consider the cone
    \begin{equation*}
        \sigma (\gamma_1, \dotsc, \gamma_n; \preceq)
        \subset \alpha (\gamma_1, \dotsc, \gamma_n)
    \end{equation*}
    defined by the inequalities $x_i \geq x_j$ for $i \preceq j$.

    It follows from \cref{thm-cotangent-arrangement} that
    special cones of~$\mathcal{X}$ are necessarily of the form
    $\sigma (\gamma_1, \dotsc, \gamma_n; \preceq)$,
    with $\gamma_1, \dotsc, \gamma_n$ non-zero.
\end{para}

\begin{para}[Stacks of filtrations]
    \label{para-lms-filt}
    The constancy theorem, \cref{thm-constancy},
    can be used to construct canonical stacks of filtrations
    for a linear moduli stack~$\mathcal{X}$.

    For classes $\gamma_1, \ldots, \gamma_n \in \uppi_0 (\mathcal{X})$,
    define the \emph{stack of filtrations}
    \begin{equation*}
        \mathcal{X}_{\gamma_1, \dotsc, \gamma_n}^+
        = \mathcal{X}_{\sigma (\gamma_1, \dotsc, \gamma_n)}^+ \ .
    \end{equation*}
    Then, by
    \cref{thm-constancy},
    for any map $\gamma \colon \mathbb{Z} \to \uppi_0 (\mathcal{X})$
    whose non-zero values agree with the non-zero elements in the sequence
    $\gamma_n, \dotsc, \gamma_1$, preserving order,
    we have a canonical isomorphism
    $\mathcal{X}_{\gamma_1, \dotsc, \gamma_n}^+
    \simeq \mathcal{X}_\gamma^+$,
    where we identify~$\gamma$
    with an element of $|\mathrm{CL} (\mathcal{X})|$.
    The stack~$\mathcal{X}_{\gamma_1, \dotsc, \gamma_n}^+$
    can thus be thought of as parametrizing filtrations
    with stepwise quotients of classes
    $\gamma_1, \dotsc, \gamma_n$, in that order.

    More generally, for classes
    $\gamma_1, \dotsc, \gamma_n$ as above
    and any partial order~$\preceq$ on the set $\{ 1, \dotsc, n \}$,
    we define a \emph{stack of filtrations}
    \begin{equation*}
        \mathcal{X}_{\smash{\gamma_1, \dotsc, \gamma_n}}^\preceq
        = \mathcal{X}_{\sigma (\gamma_1, \dotsc, \gamma_n; \preceq)}^+ \ .
    \end{equation*}
    When~$\preceq$ is the trivial partial order~$=$,
    this gives
    $\mathcal{X}_{\gamma_1} \times \cdots \times \mathcal{X}_{\gamma_n}$;
    when~$\preceq$ is the usual total order, this gives
    $\mathcal{X}_{\gamma_1, \dotsc, \gamma_n}^+$.
    The case of general partial orders was considered by
    \textcite[\S7]{joyce-2006-configurations-i}
    under a different set of axioms,
    called \emph{moduli stacks of configurations} there.
\end{para}

\begin{para}[Associativity]
    \label{para-lms-assoc}
    For a linear moduli stack~$\mathcal{X}$,
    the associativity theorem, \cref{thm-associativity},
    implies that the stacks of filtrations
    $\mathcal{X}_{\gamma_1, \dotsc, \gamma_n}^+$
    defined above fit into cartesian diagrams
    \begin{equation*}
        \begin{tikzcd}[column sep=2em, row sep={3.5em,between origins}]
            \mathcal{X}_{\gamma_1, \gamma_2, \gamma_3}^+
            \ar[r] \ar[d]
            \ar[dr, phantom, "\ulcorner", pos=.2]
            & \mathcal{X}_{\gamma_1} \times \mathcal{X}_{\gamma_2, \gamma_3}^+
            \ar[d]
            \\
            \mathcal{X}_{\gamma_1, \gamma_2 + \gamma_3}^+
            \ar[r]
            & \mathcal{X}_{\gamma_1} \times \mathcal{X}_{\gamma_2 + \gamma_3}
            \vphantom{^+} \rlap{\ ,}
        \end{tikzcd}
        \qquad
        \begin{tikzcd}[column sep=2em, row sep={3.5em,between origins}]
            \mathcal{X}_{\gamma_1, \gamma_2, \gamma_3}^+
            \ar[r] \ar[d]
            \ar[dr, phantom, "\ulcorner", pos=.2]
            & \mathcal{X}_{\gamma_1, \gamma_2}^+ \times \mathcal{X}_{\gamma_3}
            \ar[d]
            \\
            \mathcal{X}_{\gamma_1 + \gamma_2, \gamma_3}^+
            \ar[r]
            & \mathcal{X}_{\gamma_1 + \gamma_2} \times \mathcal{X}_{\gamma_3}
            \vphantom{^+} \rlap{\ ,}
        \end{tikzcd}
    \end{equation*}
    and similarly, $\mathcal{X}_{\gamma_1, \dotsc, \gamma_n}^+$
    can be expressed as $(n - 1)$-fold fibre products,
    using only stacks of the form~$\mathcal{X}_\gamma$, $\mathcal{X}_{\smash{\gamma, \gamma'}}^+$,
    and their products, where $\gamma, \gamma' \in \uppi_0 (\mathcal{X})$.
    There is also a version for general partial orders,
    analogous to \textcite[Theorem~7.10]{joyce-2006-configurations-i}.
\end{para}

\begin{para}[Hall categories]
    \label{para-lms-hall-cat}
    For a linear moduli stack~$\mathcal{X}$,
    the category $\mathsf{Hall} (\mathcal{X})^\mathrm{op}$
    can be explicitly described as follows:
    Its objects are tuples $(\gamma_1, \dotsc, \gamma_n)$
    of elements of the set
    $M \subset \uppi_0 (\mathcal{X}) \setminus \{ 0 \}$
    defined in \cref{para-lms-special-faces},
    where $n \geq 0$,
    and a morphism
    $(\gamma_1, \dotsc, \gamma_n) \to (\gamma'_1, \dotsc, \gamma'_m)$
    is a partition $\{ 1, \dotsc, n \} = I_1 \sqcup \cdots \sqcup I_m$
    such that $\gamma_i = \sum_{j \in I_i} \gamma'_j$ for all~$i$,
    together with a total order on each~$I_i$,
    but two choices of total orders are identified
    if they only differ when comparing elements
    $\gamma'_j, \gamma'_k$ for $j, k$ in some~$I_i$
    such that the morphism
    $\oplus \colon \mathcal{X}_{\gamma'_j} \times \mathcal{X}_{\gamma'_k}
    \to \mathcal{X}_{\gamma'_j + \gamma'_k}$
    is étale.
    Composition of morphisms is given by concatenating the total orders.
\end{para}

\begin{para}[Remark]
    In \cref{para-lms-hall-cat},
    the description of $\mathsf{Hall} (\mathcal{X})$
    is made complicated by the fact that
    not all faces $\alpha (\gamma_1, \dotsc, \gamma_n)$
    for $\gamma_1, \dotsc, \gamma_n \in \uppi_0 (\mathcal{X}) \setminus \{ 0 \}$
    are special faces,
    and that not all cones
    $\sigma (\gamma_1, \dotsc, \gamma_n; \preceq)$
    defined in \cref{para-lms-special-cones}
    are special cones.

    One could circumvent this by considering the cone arrangement
    $\mathsf{Cone}^\mathrm{lin} (\mathcal{X}) \subset \mathsf{Cone} (\mathcal{X})$
    consisting of all cones
    $\sigma (\gamma_1, \dotsc, \gamma_n; \preceq)$
    defined in \cref{para-lms-special-cones},
    and then defining a Hall category
    $\mathsf{Hall}^\mathrm{lin} (\mathcal{X})$
    similarly to \cref{para-hall-category},
    but using cones in $\mathsf{Cone}^\mathrm{lin} (\mathcal{X})$
    instead of special cones.
    This category is still well-defined,
    and \cref{thm-associativity}
    holds with $\mathsf{Hall}^\mathrm{lin} (\mathcal{X})$
    in place of~$\mathsf{Hall} (\mathcal{X})$.
    Taking special closures gives a functor
    $\mathsf{Hall}^\mathrm{lin} (\mathcal{X}) \to \mathsf{Hall} (\mathcal{X})$.

    Now, $\mathsf{Hall}^\mathrm{lin} (\mathcal{X})^\mathrm{op}$
    can be explicitly described as follows:
    Its objects are tuples $(\gamma_1, \dotsc, \gamma_n)$
    of elements of $\uppi_0 (\mathcal{X}) \setminus \{ 0 \}$,
    where $n \geq 0$,
    and morphisms are partitions with choices of total orders
    as in \cref{para-lms-hall-cat},
    with no extra identifications.
    In fact, $\mathsf{Hall}^\mathrm{lin} (\mathcal{X})^\mathrm{op}$
    is equivalent to the total category
    of a coloured operad over the category~$\mathsf{Fin}$ of finite sets,
    similarly to \textcite[Construction~2.1.1.7]{lurie-ha},
    whose algebras are $\uppi_0 (\mathcal{X})$-graded associative algebras
    $(A_\gamma)_{\gamma \in \uppi_0 (\mathcal{X})}$
    such that~$A_0$ is the unit object.
\end{para}

\begin{para}[Hall algebras]
    \label{para-three-functor}
    We now explain how to construct various types of
    \emph{Hall algebras} for linear moduli stacks.
    These include \emph{Hall algebras} studied by
    \textcite{ringel-1990-hall-1,ringel-1990-hall-2},
    \emph{motivic Hall algebras}
    of \textcite{joyce-2007-configurations-ii},
    and \emph{cohomological Hall algebras}
    of \textcite{kontsevich-soibelman-2011-coha}.

    Let~$(\mathcal{C}, \otimes)$ be a symmetric monoidal category
    with coproducts compatible with the tensor product.
    Assume that we are given a
    \emph{three-functor formalism} valued in~$\mathcal{C}$,
    as in \textcite{scholze-six-functor},
    that is, a lax symmetric monoidal functor
    \begin{equation*}
        \mathrm{H} \colon
        \mathcal{S} \longrightarrow \mathcal{C} \ ,
    \end{equation*}
    where $\mathcal{S} \subset \mathsf{Span} (\mathsf{St}_S)$
    is an isomorphism-closed symmetric monoidal subcategory
    of the category of spans defined in \cref{para-cat-of-spans}.

    More explicitly, this assigns to each stack
    $\mathcal{Y} \in \mathcal{S}$
    an object~$\mathrm{H} (\mathcal{Y}) \in \mathcal{C}$,
    together with the data of
    pullback and pushforward operations~$f^*, g_!$
    for spans
    $\smash{\mathcal{Y}' \overset{f}{\gets} \mathcal{Y} \overset{g}{\to} \mathcal{Y}''}$
    in~$\mathcal{S}$, and an exterior product operation
    $\boxtimes \colon \mathrm{H} (\mathcal{Y}') \otimes \mathrm{H} (\mathcal{Y}'')
    \to \mathrm{H} (\mathcal{Y}' \times \mathcal{Y}'')$,
    satisfying extra properties including the base change property.

    We explain how various types of Hall algebras
    fit into this framework:
    \begin{enumerate}
        \item
            To define Hall algebras in the style of
            \textcite{ringel-1990-hall-1,ringel-1990-hall-2},
            we assume $S = \Spec (K)$ for a finite field~$K$,
            take~$\mathrm{H} (\mathcal{X})$ to be the space of
            linear combinations of $K$-points of~$\mathcal{X}$,
            and require the legs of a span
            to be of finite type,
            with the right leg representable.

        \item
            To define motivic Hall algebras
            as in \textcite{joyce-2007-configurations-ii},
            we take~$\mathrm{H} (\mathcal{X})$ to be a ring of motives on~$\mathcal{X}$,
            where we require the legs of a span
            to be of finite type,
            with the right leg representable.

        \item
            To define cohomological Hall algebras
            in the style of \textcite{kontsevich-soibelman-2011-coha},
            we take~$\mathrm{H} (\mathcal{X})$ to be the cohomology of a smooth stack~$\mathcal{X}$,
            where we require the stacks in a span to be smooth,
            and the right leg to be proper.
    \end{enumerate}
    The following result is not difficult to verify directly
    using the associativity theorem
    in its form stated in \cref{para-lms-assoc}.
\end{para}

\begin{theorem}
    Let~$\mathcal{X}$ be a linear moduli stack,
    and let~$\mathrm{H}$ be a three-functor formalism as above,
    such that the functor
    $\mathcal{X}_{(-)} \colon \mathsf{Hall} (\mathcal{X})^\mathrm{op} \to \mathsf{Span} (\mathsf{St}_S)$
    in \cref{thm-functorial-assoc} factors through~$\mathcal{S}$.

    Then the object
    $\coprod_{\gamma \in \uppi_0 (\mathcal{X})} \mathrm{H} (\mathcal{X}_\gamma)$
    has the structure of a $\uppi_0 (\mathcal{X})$-graded
    associative algebra object in~$\mathcal{C}$,
    called the \emph{Hall algebra} of~$\mathcal{X}$,
    where the multiplication is given by the composition
    \begin{equation*}
        \mathrm{H} (\mathcal{X}_{\gamma_1}) \otimes \mathrm{H} (\mathcal{X}_{\gamma_2})
        \overset{\boxtimes}{\longrightarrow}
        \mathrm{H} (\mathcal{X}_{\gamma_1} \times \mathcal{X}_{\gamma_2})
        \overset{\mathrm{gr}^*}{\longrightarrow}
        \mathrm{H} (\mathcal{X}_{\gamma_1, \gamma_2}^+)
        \overset{(\mathrm{ev}_1)_!}{\longrightarrow}
        \mathrm{H} (\mathcal{X}_{\gamma_1 + \gamma_2}) \ ,
    \end{equation*}
    where $\mathcal{X}_{\gamma_1, \gamma_2}^+$
    is defined in \cref{para-lms-filt}.
    \qed
\end{theorem}

\begin{para}
    \label{para-gen-hall-alg}
    For a general stack~$\mathcal{X}$,
    assuming that the functor
    $\mathcal{X}_{(-)}$ in \cref{thm-functorial-assoc}
    factors through~$\mathcal{S}$,
    there is also a well-defined composition functor
    $\mathrm{H} (\mathcal{X}_{(-)}) \colon
    \mathsf{Hall} (\mathcal{X})^\mathrm{op} \to \mathcal{S} \to \mathcal{C}$.
    This should be seen as a generalized version of Hall algebras for~$\mathcal{X}$,
    where the algebra structure is replaced with
    the functorial property with respect to the Hall category.
\end{para}

\subsection{Real filtrations}\label{S:real_filtrations}

\begin{para}
    \label{para-real-filt-intro}
    As an application of the constancy theorem,
    \cref{thm-constancy}, we introduce a generalization
    of the stacks $\mathrm{Filt}^\Sigma (\mathcal{X})$
    and $\mathrm{Filt}^C (\mathcal{X})$
    of cone filtrations constructed in \cref{subsec-filt-cones}
    to the case of real cones $K \subset \mathbb{R}^n$,
    not necessarily polyhedral.
    In particular, taking $K = \mathbb{R}_{\geq 0} \subset \mathbb{R}$,
    this gives a stack $\mathrm{Filt}_\mathbb{R} (\mathcal{X})$
    of real-weighted filtrations in~$\mathcal{X}$.

    Filtrations arise in many moduli problems as optimizers of functionals
    on the space of all filtrations, but these
    optimal filtrations often have real weights. For instance, in \textcite[Theorem 1.2]{blum-liu-xu-zhuang-2023-kahler-ricci},
    the authors use a functional on the space of valuations to identify a canonical finitely generated real-weighted filtration on the anti-canonical ring of any klt log Fano pair. This gives an algebro-geometric generalization of the degeneration to a Kähler--Ricci soliton under K\"{a}hler--Ricci flow. The issue of real-weighted filtrations also arises in the study of Bridgeland stability conditions \cite[Remark~6.5.7]{halpern-leistner-instability}, where Harder--Narasimhan filtrations can be naturally indexed by either slopes or phases, both of which can be arbitrary real numbers.

    Instead of previous ad hoc workarounds to deal with real-weighted filtrations,
    we provide a general framework to work with them.

    Throughout, we fix an algebraic stack~$\mathcal{X}$
    as in \cref{assumption-stack-basic}.
    The definition of stacks of real filtrations
    will be given in \cref{para-real-filt}.
    A main result of this section is the algebraicity of the stack
    $\mathrm{Filt}^K (\mathcal{X})$,
    and in particular of $\mathrm{Filt}_\mathbb{R} (\mathcal{X})$,
    under a mild assumption, which we prove in
    \cref{thm-algebraicity-real-filtrations}.
    We also provide a criterion for a connected component of
    $\mathrm{Filt}_\mathbb{R} (\mathcal{X})$
    to be isomorphic to a component of $\mathrm{Filt} (\mathcal{X})$
    in \cref{prop-connected-components-real-filtrations}.
\end{para}

\begin{para}[Generalized graded points]
    \label{para-grad-real}
    Let $A$ be a commutative $\mathbb Q$-algebra,
    such as~$\mathbb R$ or~$\mathbb C$.
    Recall from \cref{para-formal-lattices}
    that an \emph{$A$-lattice} means a free $A$-module of finite rank.

    For an $A$-lattice $V$,
    define the \emph{stack of\/ $V^\vee$-graded points} of $\mathcal X$ as the colimit
    \begin{equation}
        \Grad^V(\mathcal X) =
        \colim_{(\Lambda,f)\in (V\downarrow \mathsf{Lat} (\Z))^{\mathrm{op}}}
        \Grad^\Lambda(\mathcal X)\ ,
    \end{equation}
    where $(V\downarrow \mathsf{Lat} (\Z))$
    is the comma category whose objects are pairs $(\Lambda, f)$
    with $\Lambda \in \mathsf{Lat} (\Z)$ and
    $f \colon V \to \Lambda \otimes_\mathbb{Z} A$ an $A$-linear map.
    We also write
    $\Grad_A (\mathcal X) = \Grad^A (\mathcal X)$.

    The assignment $(V,\cX)\mapsto \mathrm{Grad}^V(\cX)$ is contravariant on $V$ and covariant on $\cX$.
    There is a morphism $\mathrm{tot}\colon \Grad^V(\cX)\to \cX$
    defined as the colimit of the morphisms
    $\mathrm{tot} \colon \Grad^\Lambda (\cX)\to \cX$.

    We will see in \cref{thm-rational-closure}
    that $\mathrm{Grad}^V (\mathcal{X})$ is algebraic,
    with each connected component isomorphic to a component of
    $\mathrm{Grad}^n (\mathcal{X})$ for some~$n$.

    Alternatively, we also have
    \begin{equation}\label{eq-grad-v}
        \Grad^V(\mathcal X) \simeq
        \colim_{(F,f)\in (V\downarrow \mathsf{Lat} (\Q))^{\mathrm{op}}} \Grad^F(\mathcal X) \ ,
    \end{equation}
    where the comma category $(V \downarrow \mathsf{Lat} (\Q))$
    is defined analogously.
\end{para}

\begin{para}[Componentwise description]
    Denote
    $\mathrm{CL}_A(\mathcal X) = \mathrm{CL}_\Q(\mathcal X) \otimes_\mathbb{Q} A$,
    defined as in \cref{para-restriction}.
    We abbreviate $\mathsf{Face}_A(\mathcal X)=\mathsf{Face}(\mathrm{CL}_A(\mathcal X))$ and, as usual, $\mathsf{Face}(\mathcal X)=\mathsf{Face}_\Q(\mathcal X)$.
    There is an extension of scalars functor
    $(-)_A \colon \mathsf{Face}(\mathcal X)\to \mathsf{Face}_A(\mathcal X)$.

    The colimit \cref{eq-grad-v} can be rewritten as
    \begin{equation*}
        \Grad^V(\cX) \simeq
        \coprod_{\beta \colon V\to \mathrm{CL}_A(\mathcal X) } \cX_\beta \ ,
        \quad \text{where} \quad
        \cX_\beta = \colim_{(F,\alpha,f)\in (\beta\downarrow \mathsf{Face}(\mathcal X))^{\mathrm{op}}} \mathcal X_\alpha \ ,
    \end{equation*}
    where $(\beta\downarrow \mathsf{Face}(\mathcal X))$
    is the comma category whose objects are triples $(F,\alpha,f)$
    with $(F, \alpha) \in \mathsf{Face}(\mathcal{X})$
    and $f\colon (V,\beta)\to (F_A,\alpha_A)$ a morphism in $\mathsf{Face}_A(\mathcal X)$.
\end{para}

\begin{theorem}
    \label{thm-rational-closure}
    The extension of scalars functor
    $\mathsf{Face}^{\mathrm{nd}}(\mathcal X)\to \mathsf{Face}_A(\mathcal X)$
    has a left adjoint
    \begin{align*}
        (-)^\mathrm{rat} \colon
        \mathsf{Face}_A(\mathcal X)
        &\longrightarrow \mathsf{Face}^{\mathrm{nd}}(\mathcal X) \ ,
    \end{align*}
    called the \emph{rational closure} functor.
    For any $(V,\beta)\in \mathsf{Face}_A(\mathcal X)$,
    we have a canonical isomorphism
    \[
        \mathcal{X}_{\smash{\beta^\mathrm{rat}}} \longsimto
        \mathcal{X}_{\beta} \ .
    \]
    In particular, for an $A$-lattice $V$, the stack $\Grad^V(\mathcal X)$ is an algebraic stack as in \cref{assumption-stack-basic}.
\end{theorem}

\begin{proof}
    For the existence of the left adjoint,
    it is enough to show that for any face $(V,\beta)\in \mathsf{Face}_A(\cX)$,
    the comma category $(\beta~\downarrow~\mathsf{Face}^{\mathrm{nd}}(\cX))$
    has an initial object. 

    Fix $(V,\beta)\in \mathsf{Face}_A(\cX)$ and let $(F,\alpha,f)\in (\beta \downarrow \mathsf{Face}^{\mathrm{nd}}(\cX))$ be an object with $\dim F$ minimal. We claim that $(F, \alpha,f)$ is initial in $(\beta~\downarrow~\mathsf{Face}^{\mathrm{nd}}(\cX))$.

    First we note that the extension of scalars functor $\mathsf{fLat} (\Q)\to \mathsf{fLat} (A)$ preserves finite limits. 
    This follows formally from the fact that, for every $A$-lattice $V$, the category $(V\downarrow\mathsf{Lat} (\Q))$ is cofiltered, which is in turn a consequence of the functor $\mathsf{Lat} (\Q)\to \mathsf{Lat} (A)$ preserving finite limits.

    Now let $(F_1,\alpha_1, f_1)$ be another object in $(\beta \downarrow \mathsf{Face}^{\mathrm{nd}}(\cX))$. 
    Both $f$ and $f_1$ factor through $F_A\times_{\mathrm{CL}_A(\cX)}F_{1,A}=(F\times_{\mathrm{CL}_\Q(\cX)}F_1)_A$, and thus there is an object $(F_2,\alpha_2,f_2)$ of $(\beta\downarrow \mathsf{Face}(\cX))$ mapping both to $(F,\alpha,f)$ and $(F_1,\alpha_1,f_1)$. After taking the non-degenerate quotient as in \cref{theorem-non-degenerate-quotient-functor}, we may assume that $\alpha_2$ is non-degenerate. Then, by minimality of $\dim F$, the map $F_2\to F$ is an isomorphism, so after inverting it we obtain the sought-after morphism from $(F,\alpha,f)\simeq (F_2,\alpha_2,f_2)$ to $(F_1,\alpha_1,f_2)$.

    If we have two maps from $(F,\alpha,f)$ to $(F_1,\alpha_1,f_2)$, the corresponding equalizer $\mathrm{eq}(F \rightrightarrows F_1)$ must be the whole of $F$ by minimality, so the two maps are equal.

    For the final statement,
    it follows formally from \cref{theorem-non-degenerate-quotient-functor}
    that $\mathcal X_\beta$ is the colimit
    \[
        \mathcal X_\beta \simeq
        \colim_{(F,\alpha,f)\in (\beta\downarrow \mathsf{Face}^{\mathrm{nd}}(\mathcal X))^{\mathrm{op}}} \mathcal X_\alpha\ .
    \]
    Since $(\beta\downarrow \mathsf{Face}^{\mathrm{nd}}(\mathcal X))$
    has an initial object given by $(V^{\mathrm{rat}},\beta^{\mathrm{rat}},f)$,
    with $f\colon (V,\beta)\to (V^{\mathrm{rat}}_A,\beta^{\mathrm{rat}}_A)$
    the unit map, we get the desired isomorphism.
\end{proof}

\begin{para}[The special face closure]\label{para-real-special-face-closure}
    Composing the rational closure from \cref{thm-rational-closure} with the special face closure from \cref{thm-special-face-closure}, we obtain a generalized \emph{special face closure} functor
    \begin{align*}
        (-)^{\mathrm{sp}} \colon
        \mathsf{Face}_A(\mathcal X)&\longrightarrow \mathsf{Face}^{\mathrm{sp}}(\mathcal X) \ ,
    \end{align*}
    which is left adjoint to the extension of scalars functor $\mathsf{Face}^{\mathrm{sp}}(\mathcal X)\to \mathsf{Face}_A(\mathcal X)$.
    For any $(V,\beta)\in \mathsf{Face}_A(\mathcal X)$,
    we have a canonical isomorphism
    $\mathcal{X}_{\smash{\beta^{\mathrm{sp}}}} \simto \mathcal{X}_{\beta}$.
\end{para}

\begin{para}[Real cones]
    We now start to discuss a generalization of
    the stacks of cone filtrations defined in \cref{subsec-filt-cones}
    to real cones.

    A \emph{real cone} $K$ is a commutative monoid
    with an $\bR_{\geq 0}$-module structure,
    such that the groupification $V_K$ of $K$
    is a finite-dimensional real vector space,
    the canonical map $K\to V_K$ is injective,
    and $K$ is closed in $V_K$ in the euclidean topology.
    Morphisms between real cones are morphisms of $\bR_{\geq 0}$-modules.
    Equivalently, a morphism $K_1\to K_2$ between real cones is a
    linear map $f\colon V_{K_1}\to V_{K_2}$ such that $f(K_1)\subset K_2$.
    Real cones form a category, denoted $\mathsf{Cone}_{\mathbb R}$.
\end{para}

\begin{para}[Real filtrations]
    \label{para-real-filt}
    For a real cone~$K$,
    define the \emph{stack of\/ $K$-filtrations} of $\mathcal X$ as the colimit
    \vspace{-6pt}
    \begin{equation}
        \Filt^K(\mathcal X) =
        \colim_{(\Sigma,f)\in (K\downarrow \mathsf{Cone}_\Z)^{\mathrm{op}}} \Filt^\Sigma(\mathcal X)\ ,
    \end{equation}
    where~$\mathsf{Cone}_\mathbb{Z}$ is as in \cref{para-integral-cones},
    and $(K \downarrow \mathsf{Cone}_\Z)$ is the comma category
    of pairs $(\Sigma,f)$ with
    $\Sigma \in \mathsf{Cone}_\Z$ and
    $f\colon K\to \Sigma_\bR$ a morphism of real cones,
    where
    $\Sigma_\bR = \mathbb{R}_{\geq 0} \cdot \Sigma \subset \Lambda_\Sigma \otimes_\mathbb{Z} \mathbb{R}$.

    In particular, the stack $\mathrm{Filt}_\mathbb{R} (\mathcal{X})
    = \Filt^{\mathbb{R}_{\geq 0}} (\mathcal{X})$
    is called the \emph{stack of real filtrations} of~$\mathcal{X}$.

    A~condition for the stack $\Filt^K(\mathcal X)$ to be algebraic
    is given in \cref{thm-algebraicity-real-filtrations} below.

    The assignment $(K,\cX)\mapsto \mathrm{Filt}^K(\cX)$ is contravariant on $K$ and covariant on $\cX$.
    We have the induced morphisms
    \vspace{-12pt}
    \begin{equation*}
        \begin{tikzcd}
            \Grad^{V_K} (\mathcal{X})
            \ar[rr, bend left, start anchor=north east, end anchor=north west, looseness=.6, "\smash{\mathrm{tot}}"]
            \ar[r, shift right=0.5ex, "\mathrm{sf}"']
            &
            \Filt^K (\mathcal{X})
            \ar[l, shift right=0.5ex, "\mathrm{gr}"']
            \ar[r, shift left=0.5ex, "\mathrm{ev}_0"]
            \ar[r, shift right=0.5ex, "\mathrm{ev}_1"']
            &
            \mathcal{X} \ , \vphantom{^{V_K}}
        \end{tikzcd}
    \end{equation*}
    defined by taking the colimits of the corresponding maps for each
    $\Sigma \in \mathsf{Cone}_\Z$.

    Alternatively, we also have
    \vspace{-6pt}
    \begin{equation}\label{eq-real-filt}
        \Filt^K(\mathcal X) \simeq
        \colim_{(C,f)\in (K\downarrow \mathsf{Cone}_\Q)^\mathrm{op}} \Filt^C(\mathcal X) \ ,
    \end{equation}
    where the comma category $(K\downarrow \mathsf{Cone}_\Q)$ is defined analogously.
\end{para}

\begin{para}[As mapping stacks]
    Analogously to \cref{para-rational-mapping}, for a real cone $K$,
    we can see $\Filt^K(\mathcal X)$ as a mapping stack from a non-algebraic source. Set
    $R_K=\Spec \bZ[K^\vee]$, 
    $T_K=\Spec \bZ[V_K^\vee]$, and
    $\Theta_K=R_K/T_K$.
    Then we have
    \vspace{-6pt}
    \[
        \Filt^K(\mathcal X)
        \simeq \Map(\Theta_K,\mathcal X)
        \simeq \Map_S(\Theta_K\times S,\mathcal X) \ .
    \]
    This follows by the same argument as \cref{para-rational-mapping}, observing that the category
    $(K\downarrow \mathsf{Cone}_\bZ)^\mathrm{op}$
    is filtered, as it has finite colimits, and that we have
    $R_{K,i}=\lim_{(\Sigma,f)\in (K\downarrow \mathsf{Cone}_\bZ)^\mathrm{op}} R_{\Sigma,i}$,
    where $R_{K,i}$ is defined in the same way as $R_{C,i}$ in \cref{para-rational-mapping}.
\end{para}

\begin{para}[Componentwise description]
    \label{para-x-kappa}
    Let $\mathsf{Cone}_\bR(\mathcal X)$ be the category of real cones
    in $\mathrm{CL}_\mathbb{R} (\mathcal X)$,
    whose objects are pairs $(K, \kappa)$
    with $K\in \mathsf{Cone}_\bR$ and
    $\kappa \colon V_K \to \mathrm{CL}_\bR(\mathcal X)$ a morphism,
    seen as a morphism
    $\kappa \colon K \to \mathrm{CL}_\bR(\mathcal X)$, similarly to \cref{para-cones}.

    Then, the colimit \cref{eq-real-filt} can be rewritten as
    \vspace{-6pt}
    \begin{equation*}
        \Filt^K(\cX) \simeq \coprod_{\kappa\colon K\to \mathrm{CL}_\bR(\cX)} \cX_\kappa^+ \ ,
        \quad \text{where} \quad
        \cX_\kappa^+ = \colim_{(C,\sigma,f)\in (\kappa\downarrow \mathsf{Cone}(\cX))^{\mathrm{op}}} \mathcal X_\sigma^+ \ ,
    \end{equation*}
    where $(\kappa\downarrow \mathsf{Cone}(\cX))$
    is the comma category
    whose objects are triples $(C,\sigma,f)$ with $(C,\sigma)\in \mathsf{Cone}(\cX)$ and $f\colon (K,\kappa)\to (C_\bR,\sigma_\bR)$ a morphism in $\mathsf{Cone}_\bR(\cX)$.
\end{para}

\begin{theorem}
    \label{thm-real-special-cone-closure}
    Suppose that $\mathcal X$ has finite cotangent weights.
    Then the extension of scalars functor
    $\mathsf{Cone}^{\mathrm{sp}}(\mathcal X)\to \mathsf{Cone}_\bR(\mathcal X)$
    has a left adjoint
    \begin{align*}
        (-)^{\mathrm{sp}} \colon
        \mathsf{Cone}_\bR(\mathcal X)
        &\longrightarrow \mathsf{Cone}^{\mathrm{sp}}(\mathcal X) \ ,
    \end{align*}
    referred to as the \emph{special cone closure} functor.
    Moreover, for any $(K,\kappa)\in \mathsf{Cone}_\bR(\mathcal X)$,
    there is a canonical isomorphism
    \begin{equation*}
        \mathcal{X}_{\smash{\kappa^{\mathrm{sp}}}}^+ \longsimto
        \mathcal{X}_\kappa^+ \ .
    \end{equation*}
    In particular, for any real cone $K\in \mathsf{Cone}_\bR$,
    the stack $\Filt^K(\mathcal X)$ is an algebraic stack as in \cref{assumption-stack-basic}.
\end{theorem}

\begin{proof}
    For the first statement,
    it is enough to show that for any $(K,\kappa)\in \mathsf{Cone}_\bR(\cX)$,
    the comma category
    $(\kappa\downarrow \mathsf{Cone}^{\mathrm{sp}}(\cX))$
    has an initial object.
    Let $V=V_K$ and let $\kappa^\mathrm{sp}\colon V^{\mathrm{sp}}\to \mathrm{CL}_\Q(\cX)$ be the special face closure of the face $\kappa\colon V\to \mathrm{CL}_\bR(\cX)$. 
    Denote $f\colon V\to V^{\mathrm{sp}}_\bR$ the unit map.  
    Let $K^{\mathrm{sp}}$ be the smallest special cone in $V^\mathrm{sp}$ such that $K^{\mathrm{sp}}_\bR$ contains $f(K)$. 
    It follows from \cref{para-real-special-face-closure} and \cref{thm-special-cone-closure} that $(K^{\mathrm{sp}},\kappa^{\mathrm{sp}},f)$ is initial in $(\kappa\downarrow \mathsf{Cone}^{\mathsf{sp}}(\mathcal X))$.

    For the second statement,
    it follows formally from \cref{thm-special-cone-closure} that $\mathcal X_\kappa^+$ can be expressed as the colimit
    \[
        \mathcal X_\kappa^+ =
        \colim_{(C,\sigma,f)\in (\kappa\downarrow \mathsf{Cone}^{\mathrm{sp}}(\mathcal X))^{\mathrm{op}}} \mathcal X_\sigma^+\ .
    \]
    By \cref{thm-real-special-cone-closure},
    the category $(\kappa\downarrow \mathsf{Cone}^{\mathrm{sp}}(\cX))$
    has an initial object $(K^{\mathrm{sp}},\kappa^{\mathrm{sp}},f)$,
    where $f\colon (K,\kappa)\to (K^{\mathrm{sp}},\kappa^{\mathrm{sp}})_{\bR}$ is the unit map,
    so we obtain the desired isomorphism.
\end{proof}

\begin{para}[Stacks with locally finite cotangent weights]
    \label{para-locally-finite-cotangent-weights}
    We now extend the algebraicity statement in
    \cref{thm-real-special-cone-closure}
    to a more general class of stacks.

    The stack $\mathcal X$ is said to have \emph{locally finite cotangent weights} if every quasi-compact open substack $\mathcal U$ of $\mathcal X$ has finite cotangent weights in the sense of \cref{para-finite-cotangent-weights}.

    For example, the stack $\mathcal X$ has locally finite cotangent weights if either
    \begin{enumerate}
        \item $\mathcal X$ has finite cotangent weights; or
        \item $\mathcal X$ has quasi-compact graded points, as in \cref{para-quasi-compact-graded-points}.
    \end{enumerate}
\end{para}

\begin{theorem}[Algebraicity]\label{thm-algebraicity-real-filtrations}
    Suppose that $\mathcal X$ has locally finite cotangent weights,
    and let~$K$ be a real cone.
    Then the stack $\Filt^K(\mathcal X)$ is an algebraic stack satisfying the assumptions in \cref{assumption-stack-basic}.
\end{theorem}
\begin{proof}
    \allowdisplaybreaks
    For every rational cone $C\in \mathsf{Cone}_\Q$, we may express
    \[\Filt^C(\mathcal X) \simeq \colim_{\cU\subset \mathcal X} \Filt^C(\cU) \ ,\]
    where the colimit runs over all quasi-compact open substacks $\cU$ of $\cX$. By commutativity of colimits with colimits, we have
    \begin{align*}
        \Filt^K(\mathcal X)
        & \simeq \colim_{(C,f)\in(K\downarrow \mathsf{Cone}_\Q)^{\mathrm{op}}} \Filt^C(\mathcal X)
        \\
        & \simeq \colim_{(C,f)\in (K\downarrow \mathsf{Cone}_\Q)^{\mathrm{op}}} \
        \Bigl( {} \colim_{\cU\subset \mathcal X} \Filt^C(\cU) \Bigr)
        \\
        & \simeq \colim_{\cU\subset \mathcal X} \
        \Bigl( {} \colim_{(C,f)\in(K\downarrow \mathsf{Cone}_\Q)^{\mathrm{op}}} \Filt^C(\cU) \Bigr)
        \simeq \colim_{\cU\subset \mathcal X} \Filt^K(\cU) \ .
    \end{align*}
    The last colimit is filtered and all maps involved are open immersions. Moreover, by \cref{thm-real-special-cone-closure}, each stack $\mathrm{Filt}^K(\mathcal U)$ satisfies the assumptions in \cref{assumption-stack-basic}. The result follows.
\end{proof}

\begin{proposition}\label{prop-connected-components-real-filtrations}
    If\/ $\mathcal X$ is quasi-compact and has quasi-compact graded points,
    then every connected component of\/
    $\mathrm{Filt}_\bR(\mathcal X)$
    is non-canonically isomorphic to a component of\/ $\mathrm{Filt}(\mathcal X)$.

    Moreover, the corresponding components of\/
    $\mathrm{Grad}_\bR(\mathcal X)$ and $\mathrm{Grad}(\mathcal X)$
    are also isomorphic,
    and these isomorphisms are compatible with the morphisms
    $\mathrm{tot}$, $\mathrm{gr}$, $\mathrm{sf}$, $\mathrm{ev}_0$, and $\mathrm{ev}_1$.
\end{proposition}
\begin{proof}
    Let $(\bR_{\geq 0},\kappa)\in \mathsf{Cone}_\bR(\mathcal X)$
    and let $(C,\sigma)\in \mathsf{Cone}^{\mathrm{sp}}(\mathcal X)$
    be its special closure,
    so that we have isomorphisms
    $\mathcal X_\sigma^+ \simto \mathcal X_{\kappa}^+$
    and $\mathcal X_{\mathrm{span}(\sigma)} \simto \mathcal X_{\mathrm{span}(\kappa)}$.
    Then~$\sigma$ is also the special closure
    of the cotangent chamber in the special face
    $\alpha = \mathrm{span} (\sigma)$ containing~$\kappa$;
    see also the discussion in \cref{para-special-closure-rays}.
    By the finiteness theorem, \cref{thm-finiteness},
    the face~$\alpha$ only contains finitely many special subfaces,
    so there exists a rational ray~$\mathbb{Q}_{\geq 0} \cdot \lambda \subset \alpha$
    which does not lie in any of these special subfaces,
    so that $(\mathbb{Q}_{\geq 0} \cdot \lambda)^\mathrm{sp} = \sigma^\mathrm{sp}$,
    and the result follows.
\end{proof}

\subsection{Semistable reduction}
\label{subsec-semistable-reduction}

\begin{para}
    We apply our main results
    to establish a generalization of the
    \emph{semistable reduction theorem} of
    \textcite[Theorem~B]{alper-halpern-leistner-heinloth-2023-moduli}
    to the cases of \emph{real $\Theta$-stratifications}
    and \emph{cone $\Theta$-stratifications},
    in \cref{thm-semistable-reduction} below. A similar result was previously
    described in \textcite[Corollary~3.10, Theorem~3.18]{odaka2024stability}, but we feel
    it is useful to present a different approach to its proof here.
    See \cite[p.~28--29]{odaka2024stability} for the detailed summary comparing these two approaches.

    We expect that the real filtrations that arise in the context
    of klt Fano varieties or Bridgeland stability conditions,
    as discussed in \cref{para-real-filt-intro},
    define a real $\Theta$-stratification of the corresponding moduli stack. In this case the
    semistable reduction theorem below would provide a more canonical approach to proving properness of the moduli spaces of semistable objects. (See \textcite{blum-halpern-leistner-liu-xu-2021-properness} for example.)
\end{para}

\begin{para}[Real \texorpdfstring{$\Theta$}{Θ}-stratifications]
    \label{D:theta_stratification}
    Let $\mathcal{X}$ be a stack as in \cref{assumption-stack-basic},
    with locally finite cotangent weights
    as in \cref{para-locally-finite-cotangent-weights}.

    A \emph{well-ordered real $\Theta$-stratification} of~$\mathcal{X}$
    consists of an open substack $\mathcal{S} \subset \Filt_\mathbb{R} (\mathcal{X})$
    and a well-ordering of $\uppi_0(\mathcal{S})$,
    such that for all $c \in \uppi_0(\mathcal{S})$:
    \begin{enumerate}
        \item $\mathrm{ev}_1 \colon \mathcal{S} \to \mathcal{X}$ is universally bijective.
        \item The union $\bigcup_{d>c} \mathrm{ev}_1(\mathcal{S}_d) \subset |\mathcal{X}|$ is closed, where $\mathcal{S}_d \subset \mathcal{S}$ denotes the connected component indexed by $d \in \uppi_0(\mathcal{S})$.
        \item Letting $\mathcal{X}_{\leq c} = \mathcal{X} \setminus \bigcup_{d>c} \mathcal{S}_d$, the component $\mathcal{S}_c$ is a closed substack of the open substack $\Filt_\mathbb{R} (\mathcal{X}_{\leq c}) \subset \Filt_\mathbb{R} (\mathcal{X})$.
        \item
            \label{D:theta_stratification-closed-immersion}
            $\mathrm{ev}_1 \colon \mathcal{S}_c \to \mathcal{X}_{\leq c}$ is a closed immersion.
    \end{enumerate}
    This generalizes \emph{well-ordered $\Theta$-stratifications} in
    \textcite[Definition~2.1.2]{halpern-leistner-instability},
    which are the special case where the strata live in
    the open and closed substack
    $\mathrm{Filt} (\mathcal{X}) \subset \mathrm{Filt}_\mathbb{R} (\mathcal{X})$.

    A \emph{real $\Theta$-stratum} of~$\mathcal{X}$
    is an open and closed substack $\mathcal{S}_1 \subset \mathrm{Filt}_\mathbb{R} (\mathcal{X})$
    such that $\mathrm{ev}_1 \colon \mathcal{S}_1 \to \mathcal{X}$ is a closed immersion.
    This can be seen as a special case of the definition above,
    where $\uppi_0 (\mathcal{S}) = \{ 0, 1 \}$
    and $\mathcal{S}_0$ consists of trivial filtrations,
    that is,
    $\mathrm{ev}_1 \colon \mathcal{S}_0 \to \mathcal{X} \setminus \mathcal{S}_1$
    is an isomorphism.

    These notions will be investigated more thoroughly in future work of
    D.~Halpern-Leistner and A.~Ibáñez Núñez.
\end{para}

\begin{para}[Cone \texorpdfstring{$\Theta$}{Θ}-stratifications]
    We also introduce variants of $\Theta$-stratifications
    for cone filtrations defined in \cref{subsec-filt-cones}
    and real cone filtrations defined in \cref{S:real_filtrations}.

    Let~$\mathcal{X}$ be a stack as in \cref{assumption-stack-basic}.
    A \emph{well-ordered cone $\Theta$-stratification} of $\mathcal{X}$
    is an open substack
    $\mathcal{S} \subset \coprod_C \Filt^C (\mathcal{\cX})$
    and a well-ordering of $\uppi_0(\mathcal{S})$
    satisfying the conditions in \cref{D:theta_stratification},
    where we range over all isomorphism types of rational cones~$C$,
    and $\Filt^C(\mathcal{X})$ is defined in \cref{para-rational-filt}.

    Similarly, if~$\mathcal{X}$ has locally finite cotangent weights,
    a \emph{well-ordered real cone $\Theta$-stratification} of $\mathcal{X}$
    is an open substack
    $\mathcal{S} \subset \coprod_K \Filt^K (\mathcal{\cX})$
    and a well-ordering of $\uppi_0(\mathcal{S})$
    satisfying the conditions in \cref{D:theta_stratification},
    where we range over all isomorphism types of real cones~$K$,
    and $\Filt^K(\mathcal{X})$ is defined in \cref{para-real-filt}.

    We can also define \emph{cone $\Theta$-strata}
    and \emph{real cone $\Theta$-strata} similarly.

    Real cone $\Theta$-stratifications generalize
    both real $\Theta$-stratifications and
    cone $\Theta$-stratifications,
    where the latter fact is because for any rational cone~$C$, we have
    $\Filt^C(\mathcal{X}) \subset \Filt^{C_\mathbb{R}} (\mathcal{X})$
    as an open and closed substack.
\end{para}

\begin{para}[A condition on cotangent weights]
    \label{para-real-cotangent-weights}
    For a face $\beta \colon V \to \mathrm{CL}_\mathbb{R} (\mathcal{X})$,
    where~$V$ is a finite-dimensional $\mathbb{R}$-vector space, let
    $\beta^\mathrm{rat} \in \mathsf{Face} (\mathcal{X})$
    be the rational closure of~$\beta$,
    as in \cref{thm-rational-closure}.

    Define the \emph{cotangent arrangement}
    $\Phi^\mathrm{cot} (\mathcal{X}, \beta)$
    as the set of hyperplanes in~$V$
    which are restrictions of hyperplanes in
    $\Phi (\mathcal{X}, \beta^\mathrm{rat})$
    defined in \cref{para-cotangent-arrangement}.
    Note that this set need not be finite in general.

    Recall the notation $\mathsf{Cone}_\mathbb{R} (\mathcal X)$
    from \cref{para-x-kappa}.
    We say that a real cone
    $(K, \kappa) \in \mathsf{Cone}_\mathbb{R} (\mathcal X)$
    is \emph{contained in a cotangent chamber},
    if, writing $(V, \beta) = \mathrm{span} (K, \kappa)$,
    the interior of~$K$ in~$V$ does not intersect with
    any hyperplane in $\Phi^\mathrm{cot} (\mathcal{X}, \beta)$.

    Note that in this case, the notion of \emph{cotangent chambers}
    as in \cref{para-special-closure-rays}
    may not be well-defined unless~$\mathcal{X}$ has finite cotangent weights,
    but the above condition is always meaningful.
\end{para}

\begin{theorem}[Semistable reduction]
    \label{thm-semistable-reduction}
    Let~$\mathcal{X}$ be a stack as in \cref{assumption-stack-basic},
    with quasi-compact graded points as in
    \cref{para-quasi-compact-graded-points},
    and suppose we are given one of the following:

    \begin{enumerate}
        \item
            \label{item-semistable-reduction-real}
            A well-ordered real $\Theta$-stratification~$\mathcal{S}$
            of\/~$\mathcal{X}$, or
        \item
            \label{item-semistable-reduction-real-cone}
            A well-ordered real cone $\Theta$-stratification~$\mathcal{S}$
            of\/~$\mathcal{X}$, such that~$\mathcal{S}$
            lies in the union of connected components
            $\mathcal{X}_\kappa^+ \subset \Filt^K(\mathcal{X})$
            corresponding to cones $(K, \kappa) \in \mathsf{Cone}_\bR(\mathcal X)$
            contained in cotangent chambers
            in the sense of\/ \cref{para-real-cotangent-weights}.
    \end{enumerate}
    Then for any map from a discrete valuation ring
    $\xi \colon {\Spec(R)} \to \mathcal{X}$
    whose generic point lies in~$\mathcal{X}_{\leq c}$ for some~$c$,
    there is a finite extension $R \subset R'$
    and a map $\xi' \colon {\Spec(R')} \to \mathcal{X}_{\leq c}$
    that agrees with $\xi$ at the generic point.
\end{theorem}

This generalizes the semistable reduction theorem
of \textcite[Theorem~B]{alper-halpern-leistner-heinloth-2023-moduli},
which is the special case of usual $\Theta$-stratifications.

This theorem is mainly applied to show that when $\mathcal{X}$ satisfies the existence part of the valuative criterion for properness, then so does $\mathcal{X}_{\leq c}$ for any $c$. In particular, if $\mathcal{X}_{\leq 0}$ has a separated good moduli space for some $0 \in \uppi_0(\mathcal{S})$, then the moduli space is proper.

\begin{proof}
    The case~\cref{item-semistable-reduction-real}
    is a special case of~\cref{item-semistable-reduction-real-cone},
    since the condition on cotangent weights
    is automatically satisfied if $K \simeq \mathbb{R}_{\geq 0}$.
    It is thus enough to prove the theorem in the case~\cref{item-semistable-reduction-real-cone}.

    By induction, it suffices to show that for a morphism
    $\Spec(R) \to \mathcal{X}_{\leq c}$
    whose generic point does not lie in
    $\mathrm{ev}_1(\mathcal{S}_c)$,
    there is a morphism
    $\Spec(R') \to \mathcal{X}_{\leq c} \setminus \mathrm{ev}_1(\mathcal{S}_c)$
    that agrees with the original morphism at the generic point.
    Therefore, it is enough to prove the theorem
    in the case of a single real cone $\Theta$-stratum,
    that is, when $\uppi_0(\mathcal{S}) \simeq \{0,1\}$.

    Suppose that~$\mathcal{S}_1$ lies in the component
    $\mathcal{X}_\kappa^+ \subset \Filt^K(\mathcal{X})$
    for some $(K, \kappa) \in \mathsf{Cone}_\bR(\mathcal X)$,
    and let $(V, \beta) = \mathrm{span} (K, \kappa)$.
    Choose a quasi-compact open substack
    $\mathcal{S}'_1 \subset \mathcal{S}_1$ saturated under $\mathrm{gr}$ and
    containing the preimage of $\Spec (R)$ in~$\mathcal{S}_1$,
    and then choose a quasi-compact open substack
    $\mathcal{X}' \subset
    \mathcal{X} \setminus \mathrm{ev}_1 (\mathcal{S}_1 \setminus \mathcal{S}'_1)$
    containing the images of $\Spec (R)$ and $\mathcal{S}'_1$.
    The induced morphism
    $\Filt^K (\mathcal{X}') \to \Filt^K (\mathcal{X})$ is an open immersion,
    and $\mathcal{S}'_1$ defines a real cone $\Theta$-stratum of~$\mathcal{X}'$,
    which satisfies the condition on cotangent weights,
    since for any $(V, \beta') \in \mathsf{Face}_\mathbb{R} (\mathcal{X}')$
    lying over $(V, \beta)$, we have
    $\Phi^\mathrm{cot} (\mathcal{X}', \beta')
    \subset \Phi^\mathrm{cot} (\mathcal{X}, \beta)$
    as sets of hyperplanes in~$V$.
    It is thus enough to prove the theorem
    for $\mathcal{X}'$ and $\mathcal{S}'_1$,
    so we may assume that~$\mathcal{X}$ is quasi-compact.

    In this case, $\mathcal{X}$ has finite cotangent weights.
    By assumption, the cone $K \subset V$
    is contained in a chamber of the cotangent arrangement
    $\Phi^\mathrm{cot} (\mathcal{X}, \beta)$.
    We may now replace~$\beta$ by its special closure~$\alpha = \beta^\mathrm{sp}$,
    as in \cref{para-real-special-face-closure},
    and~$\kappa$ by a cotangent chamber~$\sigma \subset \alpha$
    such that $\kappa \subset \sigma_\bR$,
    so that $\kappa^\mathrm{sp} = \sigma^\mathrm{sp}$,
    and hence $\mathcal{X}_\kappa^+ \simeq \mathcal{X}_\sigma^+$
    by \cref{thm-real-special-cone-closure}.
    As in the proof of \cref{prop-connected-components-real-filtrations},
    the finiteness theorem ensures that there exists a rational ray
    $\mathbb{Q}_{\geq 0} \cdot \lambda \subset \alpha$
    such that $(\mathbb{Q}_{\geq 0} \cdot \lambda)^\mathrm{sp} = \sigma^\mathrm{sp}$,
    so that $\mathcal{X}_\sigma^+ \simeq \mathcal{X}_\lambda^+$
    for some $\lambda \in |\mathrm{CL}_\mathbb{Q} (\mathcal{X})|$.
    We can now identify~$\mathcal{S}_1$
    with a usual $\Theta$-stratum lying in
    $\mathcal{X}_\lambda^+ \subset \Filt_\mathbb{Q} (\mathcal{X})$,
    and hence one lying in $\Filt (\mathcal{X})$.
    The theorem now follows from the case of a usual $\Theta$-stratum
    proved in \textcite[Theorem~B]{alper-halpern-leistner-heinloth-2023-moduli}.
\end{proof}

\clearpage
\preparebibliography
\printbibliography

@article{bu-self-dual-1,
    author       = {Bu, Chenjing},
    title        = {Enumerative invariants in self-dual categories. {I}: {Motivic} invariants},
    version      = {preprint v4},
    eprinttype   = {arxiv},
    eprint       = {2302.00038v4},
    date         = {2025-03-28},
    sortyear     = {2023},
}

@article{bu-self-dual-2,
    author       = {Bu, Chenjing},
    title        = {Enumerative invariants in self-dual categories. {II}: {Homological} invariants},
    version      = {preprint v2},
    eprinttype   = {arxiv},
    eprint       = {2309.00056v2},
    date         = {2023-09-11},
}

@article{ibanez-nunez-filtrations,
    author       = {Ibáñez Núñez, Andrés},
    title        = {Refined Harder–Narasimhan filtrations in moduli theory},
    version      = {preprint v1},
    eprinttype   = {arxiv},
    eprint       = {2311.18050v1},
    date         = {2023-11-29},
}

@article{halpern-leistner-instability,
    author       = {Halpern-Leistner, Daniel},
    title        = {On the structure of instability in moduli theory},
    version      = {preprint v5},
    eprinttype   = {arxiv},
    eprint       = {1411.0627v5},
    date         = {2022-02-04},
    sortyear     = {2014},
}

@article{halpern-leistner-derived,
    author       = {Halpern-Leistner, Daniel},
    title        = {Derived $\Theta$-stratifications and the $D$-equivalence conjecture},
    version      = {preprint v2},
    eprinttype   = {arxiv},
    eprint       = {2010.01127v2},
    date         = {2021-06-18},
}

@article{part-ii,
    author       = {Bu, Chenjing and Ibáñez Núñez, Andrés and Kinjo, Tasuki},
    title        = {Intrinsic Donaldson--Thomas theory. {II}: {Stability} measures and invariants},
    version      = {preprint v1},
    eprinttype   = {arxiv},
    eprint       = {2502.20515v1},
    date         = {2025-02-27},
}

@article{part-iii,
    author       = {Bu, Chenjing and Ibáñez Núñez, Andrés and Kinjo, Tasuki},
    title        = {Intrinsic Donaldson--Thomas theory. {III}: {Wall}-crossing and applications},
    note         = {in preparation},
}

@article{kontsevich-soibelman-motivic-dt,
    author       = {Kontsevich, Maxim and Soibelman, Yan},
    title        = {Stability structures, motivic {Donaldson}–{Thomas} invariants and cluster transformations},
    version      = {preprint v1},
    eprinttype   = {arxiv},
    eprint       = {0811.2435v1},
    date         = {2008-11-16},
}

@article{hall-rydh-2019-tannaka,
    author       = {Hall, Jack and Rydh, David},
    title        = {Coherent {Tannaka} duality and algebraicity of {Hom}-stacks},
    volume       = {13},
    doi          = {10.2140/ant.2019.13.1633},
    pages        = {1633--1675},
    number       = {7},
    journaltitle = {Algebra Number Theory},
    eprinttype   = {arxiv},
    eprint       = {1405.7680},
    year         = {2019},
}

@book{joyce-song-2012-dt,
    author       = {Joyce, Dominic and Song, Yinan},
    title        = {A theory of generalized {Donaldson}–{Thomas} invariants},
    series       = {Mem. Am. Math. Soc.},
    number       = {1020},
    publisher    = {American Mathematical Society},
    eprinttype   = {arxiv},
    eprint       = {0810.5645},
    year         = {2012},
    doi          = {10.1090/s0065-9266-2011-00630-1},
}

@article{alper-halpern-leistner-heinloth-2023-moduli,
    author       = {Alper, Jarod and Halpern-Leistner, Daniel and Heinloth, Jochen},
    title        = {Existence of moduli spaces for algebraic stacks},
    volume       = {234},
    doi          = {10.1007/s00222-023-01214-4},
    pages        = {949--1038},
    number       = {3},
    journaltitle = {Invent. Math.},
    eprinttype   = {arxiv},
    eprint       = {1812.01128},
    year         = {2023},
}

@article{artin-zhang-2001-hilbert,
    author       = {Artin, M. and Zhang, J. J.},
    title        = {Abstract {Hilbert} schemes},
    volume       = {4},
    doi          = {10.1023/a:1012006112261},
    pages        = {305--394},
    number       = {4},
    journaltitle = {Algebr. Represent. Theory},
    year         = {2001},
}

@article{joyce-2006-configurations-i,
    author       = {Joyce, Dominic},
    title        = {Configurations in abelian categories. {I}: {Basic} properties and moduli stacks},
    volume       = {203},
    doi          = {10.1016/j.aim.2005.04.008},
    pages        = {194--255},
    number       = {1},
    journaltitle = {Adv. Math.},
    eprinttype   = {arxiv},
    eprint       = {math/0312190},
    year         = {2006},
}

@article{joyce-2007-configurations-ii,
    author       = {Joyce, Dominic},
    title        = {Configurations in abelian categories. {II}: {Ringel}–{Hall} algebras},
    volume       = {210},
    doi          = {10.1016/j.aim.2006.07.006},
    pages        = {635--706},
    number       = {2},
    journaltitle = {Adv. Math.},
    eprinttype   = {arxiv},
    eprint       = {math/0503029},
    year         = {2007},
}

@article{joyce-2007-configurations-iii,
    author       = {Joyce, Dominic},
    title        = {Configurations in abelian categories. {III}: {Stability} conditions and identities},
    volume       = {215},
    doi          = {10.1016/j.aim.2007.04.002},
    pages        = {153--219},
    number       = {1},
    journaltitle = {Adv. Math.},
    eprinttype   = {arxiv},
    eprint       = {math/0410267},
    year         = {2007},
}

@article{joyce-2008-configurations-iv,
    author       = {Joyce, Dominic},
    title        = {Configurations in abelian categories. {IV}: {Invariants} and changing stability conditions},
    volume       = {217},
    doi          = {10.1016/j.aim.2007.06.011},
    pages        = {125--204},
    number       = {1},
    journaltitle = {Adv. Math.},
    eprinttype   = {arxiv},
    eprint       = {math/0410268},
    year         = {2008},
}

@article{joyce-2007-stack-functions,
    author       = {Joyce, Dominic},
    title        = {Motivic invariants of {Artin} stacks and ``stack functions''},
    volume       = {58},
    doi          = {10.1093/qmath/ham019},
    pages        = {345--392},
    number       = {3},
    journaltitle = {Q. J. Math.},
    eprinttype   = {arxiv},
    eprint       = {math/0509722},
    year         = {2007},
}

@article{drinfeld-gaitsgory-2014-braden,
    author       = {Drinfeld, V. and Gaitsgory, D.},
    title        = {On a theorem of {Braden}},
    volume       = {19},
    doi          = {10.1007/s00031-014-9267-8},
    pages        = {313--358},
    number       = {2},
    journaltitle = {Transform. Groups},
    eprinttype   = {arxiv},
    eprint       = {1308.3786},
    year         = {2014},
}

@article{halpern-leistner-preygel-2023-mapping-stacks,
    author       = {Halpern-Leistner, Daniel and Preygel, Anatoly},
    title        = {Mapping stacks and categorical notions of properness},
    volume       = {159},
    doi          = {10.1112/s0010437x22007667},
    pages        = {530--589},
    number       = {3},
    journaltitle = {Compos. Math.},
    eprinttype   = {arxiv},
    eprint       = {1402.3204},
    year         = {2023},
}

@article{alper-2013-good-moduli,
    author       = {Alper, Jarod},
    title        = {Good moduli spaces for {Artin} stacks},
    volume       = {63},
    doi          = {10.5802/aif.2833},
    pages        = {2349--2402},
    number       = {6},
    journaltitle = {Ann. Inst. Fourier},
    eprinttype   = {arxiv},
    eprint       = {0804.2242},
    year         = {2013},
}

@article{kinjo-park-safronov-coha,
    author       = {Kinjo, Tasuki and Park, Hyeonjun and Safronov, Pavel},
    title        = {Cohomological {Hall} algebras for $3$-{Calabi}–{Yau} categories},
    version      = {preprint v1},
    eprinttype   = {arxiv},
    eprint       = {2406.12838v1},
    date         = {2024-06-18},
}

@book{lurie-2009-htt,
    author       = {Lurie, Jacob},
    title        = {Higher topos theory},
    series       = {Ann. Math. Stud.},
    number       = {170},
    publisher    = {Princeton University Press},
    year         = {2009},
    doi          = {10.1515/9781400830558},
}

@book{lurie-ha,
    author       = {Lurie, Jacob},
    title        = {Higher algebra},
    date         = {2017-09-18},
    url          = {https://www.math.ias.edu/~lurie/papers/HA.pdf},
}

@thesis{lurie-2004-dag,
    author       = {Lurie, Jacob},
    title        = {Derived algebraic geometry},
    institution  = {Massachusetts Institute of Technology},
    type         = {phdthesis},
    year         = {2004},
}

@book{lurie-sag,
    author       = {Lurie, Jacob},
    title        = {Spectral algebraic geometry},
    date         = {2018-02-03},
    url          = {https://www.math.ias.edu/~lurie/papers/SAG-rootfile.pdf},
}

@book{toen-vezzosi-2008-hag-ii,
    author       = {Toën, Bertrand and Vezzosi, Gabriele},
    title        = {Homotopical algebraic geometry {II}: Geometric stacks and applications},
    series       = {Mem. Am. Math. Soc.},
    number       = {902},
    publisher    = {American Mathematical Society},
    eprinttype   = {arxiv},
    eprint       = {math/0404373},
    year         = {2008},
    doi          = {10.1090/memo/0902},
}

@article{toen-vaquie-2007-moduli,
    author       = {Toën, Bertrand and Vaquié, Michel},
    title        = {Moduli of objects in dg-categories},
    volume       = {40},
    doi          = {10.1016/j.ansens.2007.05.001},
    pages        = {387--444},
    number       = {3},
    journaltitle = {Ann. Sci. Éc. Norm. Supér.},
    eprinttype   = {arxiv},
    eprint       = {math/0503269},
    year         = {2007},
}

@incollection{donaldson-thomas-1998,
    author       = {Donaldson, S. K. and Thomas, R. P.},
    editor       = {Huggett, S. A. and others},
    title        = {Gauge theory in higher dimensions},
    pages        = {31--47},
    booktitle    = {The geometric universe: Science, geometry, and the work of Roger Penrose},
    publisher    = {Oxford University Press},
    year         = {1998},
    doi          = {10.1093/oso/9780198500599.003.0003},
}

@article{thomas-2000-dt,
    author       = {Thomas, R. P.},
    title        = {A holomorphic {Casson} invariant for {Calabi}–{Yau} $3$-folds, and bundles on {K3} fibrations},
    volume       = {54},
    doi          = {10.4310/jdg/1214341649},
    pages        = {367--438},
    number       = {2},
    journaltitle = {J. Differ. Geom.},
    eprinttype   = {arxiv},
    eprint       = {math/9806111},
    year         = {2000},
}

@article{joyce-wall-crossing,
    author       = {Joyce, Dominic},
    title        = {Enumerative invariants and wall-crossing formulae in abelian categories},
    version      = {preprint v1},
    eprinttype   = {arxiv},
    eprint       = {2111.04694v1},
    date         = {2021-11-08},
}

@article{alper-hall-rydh-etale-local,
    author       = {Alper, Jarod and Hall, Jack and Rydh, David},
    title        = {The étale local structure of algebraic stacks},
    version      = {preprint v4},
    eprinttype   = {arxiv},
    eprint       = {1912.06162v4},
    date         = {2025-04-03},
}

@book{gubeladze-winfried-2009-polytopes,
    author       = {Gubeladze, Joseph and Bruns, Winfried},
    title        = {Polytopes, rings, and $K$-theory},
    series       = {Springer Monogr. Math.},
    publisher    = {Springer},
    year         = {2009},
    isbn         = {9780387763569},
    doi          = {10.1007/b105283},
}

@article{scholze-six-functor,
    author       = {Scholze, Peter},
    title        = {Six-functor formalisms},
    url          = {https://people.mpim-bonn.mpg.de/scholze/SixFunctors.pdf},
    date         = {2022-10},
}

@incollection{ringel-1990-hall-1,
    author       = {Ringel, Claus Michael},
    title        = {Hall algebras},
    series       = {Banach Cent. Publ.},
    pages        = {433--447},
    number       = {26, Part 1},
    booktitle    = {Topics in algebra. 1: Rings and representations of algebras},
    year         = {1990},
}

@article{ringel-1990-hall-2,
    author       = {Ringel, Claus Michael},
    title        = {Hall algebras and quantum groups},
    volume       = {101},
    doi          = {10.1007/bf01231516},
    pages        = {583--591},
    number       = {1},
    journaltitle = {Invent. Math.},
    year         = {1990},
    langid       = {english},
}

@article{kontsevich-soibelman-2011-coha,
    author       = {Kontsevich, Maxim and Soibelman, Yan},
    title        = {Cohomological Hall algebra, exponential Hodge structures and motivic Donaldson–Thomas invariants},
    volume       = {5},
    doi          = {10.4310/cntp.2011.v5.n2.a1},
    pages        = {231--352},
    number       = {2},
    journaltitle = {Commun. Number Theory Phys.},
    eprinttype   = {arxiv},
    eprint       = {1006.2706},
    year         = {2011},
}

@book{johnstone-2002-elephant-ii,
    author       = {Johnstone, Peter},
    title        = {Sketches of an elephant: A topos theory compendium. Volume~{II}},
    series       = {Oxf. Logic Guides},
    number       = {44},
    publisher    = {Oxford University Press},
    year         = {2002},
    doi          = {10.1093/oso/9780198515982.001.0001},
}

@article{borisov-joyce-2017,
    author       = {Borisov, Dennis and Joyce, Dominic},
    title        = {Virtual fundamental classes for moduli spaces of sheaves on {Calabi}–{Yau} four-folds},
    volume       = {21},
    doi          = {10.2140/gt.2017.21.3231},
    pages        = {3231--3311},
    number       = {6},
    journaltitle = {Geom. Topol.},
    eprinttype   = {arxiv},
    eprint       = {1504.00690},
    date         = {2017},
}

@article{oh-thomas-2023-i,
    author       = {Oh, Jeongseok and Thomas, Richard P.},
    title        = {Counting sheaves on {Calabi}–{Yau} 4-folds. {I}},
    doi          = {10.1215/00127094-2022-0059},
    journaltitle = {Duke Math. J.},
    volume       = {172},
    number       = {7},
    pages        = {1333--1409},
    eprinttype   = {arxiv},
    eprint       = {2009.05542},
    date         = {2023},
}

@article{oh-thomas-ii,
    author       = {Oh, Jeongseok and Thomas, Richard P.},
    title        = {Complex Kuranishi structures and counting sheaves on {Calabi}–{Yau} 4-folds, {II}},
    version      = {preprint v2},
    eprinttype   = {arxiv},
    eprint       = {2305.16441v2},
    date         = {2024-03-01},
}

@article{cao-leung-dt4,
    author       = {Cao, Yalong and Leung, Naichung Conan},
    title        = {Donaldson–Thomas theory for {Calabi}–{Yau} 4-folds},
    version      = {preprint v2},
    eprinttype   = {arxiv},
    eprint       = {1407.7659v2},
    date         = {2015-09-24},
}

@article{bu-davison-ibanez-nunez-kinjo-padurariu,
    author       = {Bu, Chenjing and Davison, Ben and Ibáñez Núñez, Andrés and Kinjo, Tasuki and Pădurariu, Tudor},
    title        = {Cohomology of symmetric stacks},
    version      = {preprint v2},
    eprinttype   = {arxiv},
    eprint       = {2502.04253v2},
    date         = {2025-05-31},
}

@article{davison-meinhardt-2020,
    author       = {Davison, Ben and Meinhardt, Sven},
    title        = {Cohomological {Donaldson}–{Thomas} theory of a quiver with potential and quantum enveloping algebras},
    volume       = {221},
    doi          = {10.1007/s00222-020-00961-y},
    pages        = {777--871},
    number       = {3},
    journaltitle = {Invent. Math.},
    eprinttype   = {arxiv},
    eprint       = {1601.02479},
    year         = {2020},
}

@article{tanaka-thomas-2020-i,
    author       = {Tanaka, Yuuji and Thomas, Richard P.},
    title        = {Vafa–Witten invariants for projective surfaces. I: Stable case},
    volume       = {29},
    doi          = {10.1090/jag/738},
    pages        = {603--668},
    number       = {4},
    journaltitle = {J. Algebr. Geom.},
    eprinttype   = {arxiv},
    eprint       = {1702.08487},
    year         = {2020},
    sortyear     = {2017},
}

@article{tanaka-thomas-2018-ii,
    author       = {Tanaka, Yuuji and Thomas, Richard P.},
    title        = {Vafa–Witten invariants for projective surfaces. {II}: Semistable case},
    volume       = {13},
    doi          = {10.4310/pamq.2017.v13.n3.a6},
    pages        = {517--562},
    number       = {3},
    journaltitle = {Pure Appl. Math. Q.},
    eprinttype   = {arxiv},
    eprint       = {1702.08488},
    year         = {2018},
}

@article{odaka2024stability,
    title        = {Stability theory over toroidal or {N}ovikov type base and canonical modifications},
    author       = {Odaka, Yuji},
    version      = {preprint v3},
    eprinttype   = {arxiv},
    eprint       = {2406.02489v3},
    date         = {2025-02-19},
}

@article{blum-liu-xu-zhuang-2023-kahler-ricci,
    author     = {Blum, Harold and Liu, Yuchen and Xu, Chenyang and Zhuang, Ziquan},
    journal    = {Forum Math. Pi},
    title      = {The existence of the {K}ähler–{R}icci soliton degeneration},
    year       = {2023},
    issn       = {2050-5086},
    volume     = {11},
    eid        = {e9},
    doi        = {10.1017/fmp.2023.5},
    eprinttype = {arxiv},
    eprint     = {2103.15278},
}

@article{blum-halpern-leistner-liu-xu-2021-properness,
    author = {Blum, Harold and Halpern-Leistner, Daniel and Liu, Yuchen and Xu, Chenyang},
    doi = {10.1007/s00029-021-00694-7},
    issn = {1022-1824},
    journal = {Sel. Math.},
    number = {4},
    title = {On properness of {K}-moduli spaces and optimal degenerations of {F}ano varieties},
    volume = {27},
    year = {2021},
    eid = {73},
    eprinttype = {arxiv},
    eprint = {2011.01895},
}

@article{chen-sun-wang-2018,
    author       = {Chen, Xiuxiong and Sun, Song and Wang, Bing},
    title        = {Kähler–Ricci flow, Kähler–Einstein metric, and K-stability},
    volume       = {22},
    doi          = {10.2140/gt.2018.22.3145},
    pages        = {3145--3173},
    number       = {6},
    journaltitle = {Geom. Topol.},
    eprinttype   = {arxiv},
    eprint       = {1508.04397},
    year         = {2018},
}

@article{liu-zheng-six-operations,
    author       = {Liu, Yifeng and Zheng, Weizhe},
    title        = {Enhanced six operations and base change theorem for higher Artin stacks},
    version      = {preprint v4},
    eprinttype   = {arxiv},
    eprint       = {1211.5948v4},
    date         = {2024-12-17},
}

@Book{laumon-champs,
    author       = {Laumon, Gérard and Moret-Bailly, Laurent},
    title        = {Champs algébriques},
    series       = {Ergeb. Math. Grenzgeb., 3. Folge},
    number       = {39},
    publisher    = {Springer},
    year         = {2000},
    doi          = {10.1007/978-3-540-24899-6},
}

@Article{antieau-gepner-brauer-groups,
    author       = {Antieau, Benjamin and Gepner, David},
    journaltitle = {Geom. Topol.},
    title        = {Brauer groups and {\'e}tale cohomology in derived algebraic geometry},
    year         = {2014},
    issn         = {1465-3060},
    number       = {2},
    pages        = {1149--1244},
    volume       = {18},
    doi          = {10.2140/gt.2014.18.1149},
    eprinttype   = {arxiv},
    eprint       = {1210.0290},
}

@book{aguiar-topics-hyperplane,
    author       = {Aguiar, Marcelo and Mahajan, Swapneel},
    title        = {Topics in hyperplane arrangements},
    series       = {Math. Surv. Monogr.},
    number       = {226},
    publisher    = {American Mathematical Society},
    year         = {2017},
    doi          = {10.1090/surv/226},
}

\authorinforule

\authorinfo{Chenjing Bu}
    {bucj@mailbox.org}
    {Mathematical Institute, University of Oxford, Oxford OX2 6GG, United Kingdom}

\authorinfo{Daniel Halpern-Leistner}
    {daniel.hl@cornell.edu}
    {Mathematics Department, Cornell University, Ithaca, NY 14850, USA}

\authorinfo{Andrés Ibáñez Núñez}
    {andres.ibaneznunez@columbia.edu}
    {Department of Mathematics, Columbia University, New York, NY 10027, USA}

\authorinfo{Tasuki Kinjo}
    {tkinjo@kurims.kyoto-u.ac.jp}
    {Research Institute for Mathematical Sciences, Kyoto University, Kyoto 606-8502, Japan}

\end{document}